\RequirePackage{fix-cm}
\documentclass[11pt]{article}
\usepackage[english]{babel}
\usepackage{url}
\usepackage{authblk}
\usepackage{comment}

\usepackage{abstract}
\usepackage{amsmath,amssymb}
\usepackage{mathtools}
\usepackage{graphicx}
\usepackage{caption}
\usepackage{subfigure}
\usepackage{algorithm}
\usepackage{algorithmic}
\usepackage{amsthm}
\usepackage{bm}
\usepackage{epstopdf}
\usepackage{float}
\usepackage{xcolor,colortbl}
\definecolor{MyBlueUrl}{rgb}{0.12,0.35,1}
\definecolor{MyredLink}{rgb}{0.9,0.11,0.17}
\definecolor{greenTcite}{rgb}{0.25,0.65,0.38}
\usepackage{enumitem}
\usepackage{hyperref}
\usepackage{cite}
\hypersetup{colorlinks=true,
linkcolor=MyredLink,
citecolor=greenTcite,
urlcolor=MyBlueUrl}
\setenumerate[0]{label=\roman*)}
\setenumerate[1]{label=\arabic*.}
\setenumerate[2]{label=(\alph*)}
\usepackage[verbose=true]{geometry}
\AtBeginDocument{
  \newgeometry{
    textheight=9in,
    textwidth=7.2in,
    top=1in,
    headheight=14pt,
    headsep=25pt,
    footskip=30pt
  }
}
\newcommand{\headeright}{}
\usepackage{fancyhdr}
\fancyheadoffset{0pt}
\newcommand{\R}{\mathbb{R}}

\newcommand{\adj}{\operatorname{adj}}
\newcommand{\dd}{\operatorname{d}\!}
\DeclareMathOperator{\sign}{\mathrm{sign}}
\newtheorem{theorem}{Theorem}

\newtheorem{lemma}[theorem]{Lemma}
\newtheorem{definition}{Definition}[section]
\newtheorem{remark}[definition]{Remark}

\newtheorem{assumption}{Assumption}

\title{On Rescuing Diverging Algorithms:\\ Crossing Infinity in Continuous Newton Methods \\ via a Self-Similarity Detection Framework}
\newcommand{\shorttitle}{Crossing infinity in continuos Newton methods}
\newcommand{\myand}{$\cdot \text{ }$}

\author[1,2,*]{Gianluca Fabiani}
\author[3]{Ralf Hannemann-Tam\'{a}s}
\author[4]{Alexander Mitsos}
\author[5]{Jonathan D.\ Hauenstein}
\author[6]{\\ Constantinos Siettos}
\author[7,*]{Panayotis G. Kevrekidis}
\author[2,8,*]{Ioannis G.  Kevrekidis}

\affil[1]{Hopkins Extreme Materials Institute, \emph{Johns Hopkins University}, Baltimore, 21218, MD, USA}
\affil[2]{Dept. of Chemical and Biomolecular Engineering, \emph{Johns Hopkins University}, Baltimore, 21218, MD, USA}
\affil[3]{Institute of Climate and Energy Systems: Energy Systems Engineering,
Forschungszentrum Jülich GmbH, Jülich, Germany.}
\affil[4]{RWTH Aachen University, Process Systems Engineering (AVT.SVT), Aachen 52074, Germany}
\affil[5]{Dept. of Applied and Computational Mathematics and Statistics, University of Notre Dame, Notre Dame, IN USA.}
\affil[6]{Dipartimento di Matematica e Applicazioni ‘‘Renato Caccioppoli", \emph{Universit\`a degli Studi di Napoli} \emph{Federico II}, Naples, Italy}
\affil[7]{Department of Mathematics and Statistics and Department of Physics, \emph{University of Massachusetts}, Amherst, 01003, MA, USA}
\affil[8]{Dept. of Applied Mathematics and Statistics, \emph{Johns Hopkins University}, Baltimore, 21218, MD, USA}
\affil[*]{Corresponding authors. Emails: \texttt{gfabian2@jh.edu}, \texttt{kevrekid@umass.edu}, \texttt{yannisk@jhu.edu}}

\date{\vspace{-3ex}\vspace{-3ex}}

\begin{document}

\maketitle

% REQUIRED
\begin{abstract}
We address the solution of systems of nonlinear algebraic equations via continuous Newton methods; these methods reformulate the discrete Newton update as an initial value problem for a system of ordinary differential equations (ODEs) whose stationary points correspond to the desired roots. A fundamental challenge arises when the Jacobian of these ODEs becomes singular along the solution trajectory, inducing blow-up in the vector field --or in the state itself-- and causing the numerical integration to break down before reaching the desired root.
To overcome this difficulty, we propose a directional self-similar blow-up detection and crossing procedure, wrapped around an adaptive step-size ODE solver with event detection. The method assumes local self-similar behavior near the blow-up; the self-similarity exponent, blow-up time, and amplitude are estimated on the fly from a short sequence of threshold crossings via linear regression, requiring no prior knowledge of the singularity location, its order, or the blow-up regime---in contrast to approaches that rely on explicit analytical continuation or problem-specific coordinate transformations.
The estimated exponent distinguishes two qualitatively different scenarios
at finite time: one in which only the vector field diverges while the state remains finite, admitting an explicit integral correction, and one in which the state itself blows-up.
In addition, infinite-time divergence of the state (the so-called ``grow-up" regime) is handled through a suitable time reparameterization.
%
%For scalar problems, explicit formulas are derived for both cases.
In higher dimensions, we track the norm of the vector field along its instantaneous direction; near a {\em simple} singularity, the flow asymptotically aligns with a single dominant direction (see text), reducing detection and estimation to a scalar problem regardless of ambient space dimension.
The approach is demonstrated on problems of increasing complexity, from scalar nonlinear equations to the steady-state solution of the Liouville-Bratu-Gelfand partial differential equation discretized on a fine spatial grid. %,
\end{abstract}
\hspace{0.5cm} {\footnotesize
%%%%%%%%%%%%%%%%%KEYWORDS
\myand Detecting Singularities
\myand Continuous Newton Method
\myand Roots of Nonlinear Algebraic Systems
}
-----
%%%%%%%%%%%%%%%%%%%%%%%%%%%% MSC Codes
%\hspace{0.2cm}
{\small
\myand 	65H10 %Numerical solutions to inverse eigenvalue problems
\myand 	65H20 %Ill-posedness and regularization problems in numerical linear algebra
\myand 	65L08  	%Numerical investigation of stability of solutions to ordinary differential equations
}
\vskip 0.05in
\hrule height 1pt
\vskip 0.1in

\sloppy
\section{Introduction}
\label{sec:intro}
Solving systems of nonlinear algebraic equations, of the form
\begin{equation}
    \bm{g}(\bm{x}) = \bm{0}, \qquad \bm{g}:\mathbb{R}^n \to \mathbb{R}^n,
    \label{eq:system_nonlin_eqs}
\end{equation}
is a central and fundamental task in computational science and engineering.
Such systems arise directly in numerous applications and, perhaps more significantly, as the inner core of many numerical methods.
The spatial discretization of steady-state problems for nonlinear partial differential equations (PDEs) -- whether by classical numerical analysis approaches, or via modern physics-informed machine learning parametrizations -- leads to large-scale nonlinear algebraic systems~\cite{quarteroni2008numerical, fabiani2021numerical}. Similarly, implicit time-stepping schemes for time-dependent PDEs and differential-algebraic equations (DAEs) require the solution of a nonlinear system at every time step~\cite{hairer2000solving, shampine1999solving, fabiani2023parsimonious, kavousanakis2026going}.
Furthermore, the numerical study of bifurcations and solution branches via continuation methods -- where one tracks how solutions vary with a parameter -- repeatedly solves nonlinear systems along the path~\cite{doedel2012numerical,  galaris2022numerical, fabiani2024task, fabiani2025enabling, fabiani2026stability}.
%Beyond direct discretizations, nonlinear systems are central to parameter estimation and inverse problems, where one seeks model parameters that minimize a misfit between predictions and observations~\cite{isakov2006inverse}.
%
%In control theory, the analysis of feedback systems often reduces to solving algebraic Riccati equations or finding equilibrium points~\cite{armaou2004time, fabiani2025equation}.

Among the numerical methods available for tackling such systems~\cite{hirsch1979algorithms, ortega2000iterative}, the Newton scheme and its variants stand as one of the most widely used, owing to their simplicity and, when convergent, their characteristically quadratic local convergence rate~\cite{kelley1995iterative, rheinboldt1998methods}. 
The classical Newton iteration,
\begin{equation}
    \bm{x}_{k+1} = \bm{x}_k - \big(\nabla \bm{g}(\bm{x}_k)\big)^{-1} \bm{g}(\bm{x}_k), \qquad k = 0,1,2,\ldots,
    \label{eq:discrete_newton}
\end{equation}
where \(\nabla \bm{g}(\bm{x}_k) \in \mathbb{R}^{n\times n}\) denotes the Jacobian matrix evaluated at \(\bm{x}_k\), transforms the nonlinear problem into a sequence of linear systems, each solved to obtain the next iterate.

However, its strength comes with a well-known limitation: convergence is guaranteed only if the Jacobian $\nabla \bm{g}$ is locally  Lipschitz continuous and the initial guess lies sufficiently close to a root~\cite{kelley1995iterative, rheinboldt1998methods}.
In practice, without prior knowledge of the root's location -- knowledge that is often unavailable for complex problems -- one cannot guarantee that the iteration will succeed. As a consequence, Newton's method may diverge, cycle, converge to an unintended root, or, as is the focus of this work, encounter numerical breakdown due to singularities encountered along the iteration path~\cite{branin1972widely}.
Furthermore, the sensitivity of Newton-type methods to the initial condition is illustrated by the classical Cayley's problem for the complex polynomial $z^3 - 1 = 0$, where the basins of attraction of the three roots form a fractal set~\cite{jacobsen2007approximations}.

One way to mitigate such difficulties is to embed the discrete Newton iteration into a continuous-time framework, giving rise to \emph{continuous Newton methods} (CNM). The underlying idea, tracing back to Davidenko~\cite{davidenko1953new} and later developed by Branin~\cite{branin1972widely} and others~\cite{diener1986trajectory,diener1987global, zufiria2002role}, is to reformulate the root-finding problem as an initial value problem for a system of ODEs of the form
\begin{equation}
\dot{\bm{x}}(t) = \bm{f}(t,\bm{x}(t)) := -\big(\nabla \bm{g}(\bm{x}(t))\big)^{-1} \bm{g}(\bm{x}(t)), \qquad \bm{x}(0) = \bm{x}_0,
\label{eq:vector_field_continuousNewton}
\end{equation}
whose stationary points correspond to roots of Eq.~\eqref{eq:system_nonlin_eqs}.
Viewing the discrete Newton iteration as a forward Euler discretization of this ODE with unit step size renders it as a particular numerical scheme within this broader continuous formulation.
%
%Several related approaches have been developed to address difficulties arising in continuous-time iterative methods and root-finding more broadly.
Continuum-time variants of discrete iterative algorithms have been developed for a variety of methods, including Nesterov's accelerated gradient method~\cite{nesterov1983method}, leading to efficient optimization algorithms~\cite{weijie2016differential,wibisono2016variational,subramanian2021exploring,ward2017toolkitsteadystatesnonlinear}, and have been generalized to a wide range of first-order optimization algorithms interpreted as physical systems with the potential for unified convergence analysis~\cite{yang2018physical}.
%, encompassing, among others, also the notion of deflation~\cite{farrell2015deflation} which may produce singularities in the vector field, due to its exclusion of already-identified roots, in order to converge to new ones.

%
Despite their appeal, the Davidenko formulation is not without drawbacks and difficulties. 
First, this formulation comes at a cost: a numerical method of choice used to simulate the CNM flow loses the characteristic Newton quadratic convergence. Indeed, any damped Newton step with $\Delta t < 1$ yields at best linear convergence with rate $C = 1-\Delta t$ and requires $1/\Delta t$ Jacobian evaluations per unit time. %Higher-order integrators approximate the continuous flow and thus inherit its linear convergence near the root with rate $C = e^{-\Delta t}$, requiring multiple evaluations and linear solves per step without improving the rate.
The Davidenko CNM is therefore best viewed as a robust exploratory tool for navigating the solution landscape, 
%particularly across singularities, 
rather than a computationally efficient alternative to Newton's method near a root.
Besides, the most fundamental challenge arises when the Jacobian $\nabla \bm{g}(\bm{x}(t))$ becomes singular at some point along its solution trajectory, that is, when $\det\big(\nabla \bm{g}(\bm{x}(t^\ast))\big)=0$, for some finite time $t^\ast$. 
At such points, the right-hand side of the ODE ceases to be well-defined, inducing a singularity in the vector field $\bm{f}(t,\bm{x})$ of Eq.~\eqref{eq:vector_field_continuousNewton}.
Depending on the structure of the singularity, two qualitatively different scenarios may arise: $\bm{f}(t,\bm{x})$ may diverge while the state $\bm{x}(t)$ remains bounded, or the state itself may diverge in finite time.
In either case, standard adaptive ODE solvers will fail to continue the integration beyond such a singularity, reducing the step size indefinitely without being able to continue past the singularity. Conversely, forcing the integrator to use artificially large step sizes---as in Branin's heuristic ``sloppy integration''---may jump past the singularity but can introduce unphysical trajectory deviations and is highly sensitive to the chosen step size, as documented by Zufiria and Guttalu~\cite{zufiria2002role}.
A third, distinct scenario may also occur independently of any Jacobian singularity: the state $\bm{x}(t)$ may diverge in infinite time, a behavior referred to as \emph{grow-up} in the dynamical systems and PDE literature to distinguish it from finite-time blow-up \cite{kobayashi1977growing, graf2026blow}.

In these cases, rather than halting the integration, we accelerate the trajectory toward infinity via a time reparameterization, reset the state at the ``opposite end" of the real line, and continue the integration while preserving the orientation of the flow -- treating $+\infty$ and $-\infty$ as the same point (in the spirit of a projective compactification of the domain~\cite{kevrekidis2017infinity}).
This projective identification is prototypical for blow-up of the form $\dot{x} = x^2$, with solution of type $x(t)=1/(t^*-t)$. Indeed, with the (singular!) change of variable $y = 1/x$, the ODE becomes $\dot{y} = -1$, whose solution $y(t) = t^* - t$ passes linearly through $y = 0$ (i.e., $x = \infty$) and continues to negative values, corresponding to $x$ re-emerging from $-\infty$.
However, for general blow-up $\dot{x} = x^p$ with odd $p > 1$, the real-valued solution ceases to exist after blow-up; its analytical continuation becomes complex, and no real trajectory continues past the singularity. The projective reset $x \to -x$ is therefore a heuristic that creates a new real trajectory on the opposite side of infinity, rather than a continuation of the original solution.
%
%In higher dimensions, such a passage through infinity may not be merely incidental but a necessary feature of the flow: the trajectory may be required to traverse unbounded regions of $\mathbb{R}^n$ before reaching a finite solution, and without a mechanism to handle this, convergence would be impossible regardless of the quality of the initial guess.
%

Branin's CNM~\cite{branin1972widely} (sometimes also called the ``global Newton method") represents a significant tweak towards global optimization (see also discussions and related path-following approaches~\cite{merrill1972applications, garcia1980relations}). Branin introduced a switching rule that reverses the sign of the vector field whenever a Jacobian singularity is encountered or a root is discovered, thereby enabling the trajectory to circumvent singularities as well as, after finding a root, continue its search for additional roots. 
This trajectory-based perspective, also pursued in homotopy and continuation methods~\cite{keller1978global, garcia1980relations, kuno1988computing, allgower2012numerical} %, offers a compelling alternative to purely local iterative schemes, with
has the ultimate goal of recovering all solutions of Eq.~\eqref{eq:system_nonlin_eqs}---a task of central importance in global optimization and in the study of multiple equilibria in nonlinear systems.
Subsequently, deflation techniques, $\gamma$-trick method, and other modifications have been incorporated with homotopy methods to generate multiple solutions~\cite{brow1971deflation, garcia1980relations, diener1986trajectory, diener1987global, kuno1988computing, farrell2015deflation}.
In the context of algebraic numerical tools, tool suites like {\em Bertini} have been devised, at the level of dynamical extensions within the complex plane, that enable the bypassing of singularities~\cite{bates2013numerically}.
Indeed, for polynomial systems, such homotopy-continuation methods in numerical algebraic geometry are designed to obtain approximations to {\em all} isolated complex solutions.
Similar singularity-crossing ideas have recently been explored in the context of multi-objective optimization~\cite{manoj2026singularities}.

Initially, Branin~\cite{branin1972widely} conjectured that his CNM with sign switching is globally convergent, provided that the zero sets of the individual components of $\bm{g}$ are homeomorphic to hyperplanes. However, Brent~\cite{brent1972davidenko} constructed a counterexample, where the Newton flow trajectory remains confined to a bounded region containing no root.
In response, Branin~\cite{branin1972widely} attributed this failure to the presence of \emph{extraneous singularities}---points where the Jacobian is singular and the adjugate satisfies $\operatorname{adj}(\nabla \bm{g}(\bm{x}))\,\bm{g}(\bm{x}) = \bm{0}$---and, in particular, to the interaction between saddle-like and center-like equilibria of the desingularized flow
\(
{\dot{\bm{x}}} = \operatorname{adj}(\nabla \bm{g}(\bm{x}))\,\bm{g}(\bm{x}).
\)

In such settings, the trajectory can inadvertently produce spurious closed loops (resembling cyclic orbits) by repeatedly switching between two nearby singularities, thereby preventing convergence to any root.
A general class of such cycling counterexamples was also described by Treccani~\cite{treccani1975convergence}, who constructed a family of twice continuously differentiable functions, topologically equivalent to convex quadratics, for which Branin's method exhibits regions of nonconvergence. %Treccani's analysis revealed that such behavior can arise from the presence of saddle points of an auxiliary function \(\beta(x,y) = \partial\phi/\partial y - \partial\phi/\partial x\), which form the boundaries of nonconvergent regions.
This demonstrated that, even for well-behaved objective functions with a unique global minimum and connected, bounded level sets, the Branin flow can fail to converge.
The conjecture that non-convergence is exclusively attributable to \emph{extraneous singularities} was later disproved in~\cite{zufiria2002role}.

The CNM and related path-following approaches have been shown to exhibit global convergence properties only under certain stricter conditions. For instance, Diener~\cite{diener1987global} proved that if the Jacobian matrix $\nabla \bm{g}$ satisfies a uniform non-degeneracy condition on the subspace orthogonal to a fixed direction, then the so-called Newton flow trajectory--the set of points where $\bm{g}(\bm{x})$ is parallel to $\bm{g}(\bm{x}_0)$---is connected and contains all roots. %This provides a sufficient condition for global convergence and for finding all solutions.

%The Branin's method have been shown to exhibit global convergence properties under certain conditions~\cite{diener1987global}, offering improved robustness over the classical discrete Newton method at the cost of a slower, asymptotic convergence rate. However, the stronger conjecture put forward by Branin~\cite{branin1971solution,branin1972widely} — that the continuous Newton method converges to a root from any initial guess — has been shown to be false in general, with counter-examples constructed in~\cite{brent1972davidenko, zufiria2002role}.
%
%\textcolor{blue!90}{might discuss here about saddle-vortex points interaction.}
%
%In both cases, the failure is caused by singularities of the Jacobian along the trajectory, which disrupt the flow and prevent convergence — the same difficulty this work aims to overcome.
%
From a numerical implementation standpoint, Branin's method is hampered by unbounded vector fields near singular manifolds, where standard ODE integrators stall or produce erratic jumps. By equipping the continuous Newton method with a systematic procedure for detecting and crossing such singularities, the present work aims to recover convergence in situations where standard integration would fail, thereby extending the practical range of applicability of the method.

A key aspect of the proposed approach concerns the orientation of the flow across singular passages. As discussed above, when a trajectory approaches a Jacobian singularity, the sign of the flow must be updated upon crossing to prevent the trajectory from cycling back toward the same non-root minimum,  or from becoming trapped. Maintaining the correct orientation through each obstruction ensures that progress toward a solution is retained and that integration continues to advance.

% Consider, for instance, a trajectory attracted toward a local minimum of Eq.~\eqref{eq:system_nonlin_eqs} that is not a root: the Jacobian vanishes there, inducing a singularity, and resuming the integration without adjusting the flow direction would cause the trajectory to cycle back toward the same point indefinitely. By appropriately flipping the local orientation of the flow upon crossing, the trajectory is instead directed away from the singularity and toward a potential solution. The same rationale applies to passages through infinity: after resetting the state at the opposite end of the real line, the flow must continue in a direction consistent with the ongoing trajectory.
% Maintaining the orientation of the flow through each such obstruction---rather than allowing the trajectory to retrace its path---ensures that progress toward a potential solution is retained and that the integration continues to advance.
%

In practice, however, Branin's original implementation relied on heuristic strategies---most notably ``sloppy integration"---to escape regions of nonconvergence. As Zufiria and Guttalu~\cite{zufiria2002role} discussed, such approaches are highly sensitive to step size and can produce numerical trajectories that deviate substantially from the theoretical flow. Here, we propose an alternative to Branin's sloppy integration by adapting the techniques of~\cite{kevrekidis2017infinity} to the continuous Newton setting. Rather than using a fixed step size to jump across singularities, we employ an adaptive ODE solver to detect the approach of a singularity, estimate its blow-up parameters on the fly, and apply a crossing strategy appropriate for the detected regime. The method operates directly on the real ODE flow without extension to the complex plane or restriction to polynomial systems, and provides a principled and numerically robust means of continuing the integration across singularities without requiring prior knowledge of their location or order.
For passages through infinity of the state itself, however, the projective reset is heuristic and does not necessarily constitute a true continuation of the solution.

Complementary to our local, on-the-fly crossing strategy, Graf et al.~\cite{graf2026blow} recently developed a global framework to classify parameter regions where {\em polynomial ODEs} exhibit finite-time blow-up.
Their approach uses Poincaré compactification to analyze the sphere at infinity and algebraic elimination to partition parameter space into regions with distinct blow-up behavior. In contrast to our method, which tracks individual trajectories of a fixed system and intervenes only when a singularity is encountered along a specific path, their framework performs a parameter-space classification based on the existence of blow-up for at least one initial condition, without simulating  any specific trajectory or attempting to cross singularities.

Beyond root-finding, the challenge of handling singularities and crossing infinity is a recurring theme across computational mathematics, appearing in signal processing~\cite{mallat1992singularity}, dynamical systems~\cite{kevrekidis2017infinity, graf2026blow}, operator theory~\cite{monfared2026algebra}, numerical methods for blow-up problems~\cite{matsue2020numerical, jia2026singularity}, and in transformations of dynamical systems where singularities can arise during learning~\cite{alvarez2023discrete}. In cosmology, transformations between frames have been used to describe the crossing of Big Bang–Big Crunch singularities~\cite{kamenshchik2016transformations}. Recent developments in scientific machine learning have explored rational neural networks for approximating functions with poles~\cite{boulle2020rational} and specialized architectures for singularity detection~\cite{derevianko2025neural}; techniques for handling self-similar blow-up have also been integrated into physics-informed neural networks~\cite{wang2023asymptotic, kavousanakis2026going}. Our work contributes to this broader frontier by providing a systematic, on-the-fly framework for detecting and crossing such singularities within the specific context of continuous Newton methods-- since such a crossing serves the original root-finding intention of the algorithm.

The remainder of the paper is organized as follows. Section~\ref{sec:continewton} introduces the continuous Newton method and its key properties. Section~\ref{sec:scalarodes} develops the infinity-crossing procedure for scalar problems, and Section~\ref{sec:multidimensionalodes} extends it to $n$-dimensional systems. Section~\ref{sec:algorithm} presents the complete algorithm. Section~\ref{sec:casestudies} illustrates the methodology on a set of test problems of increasing complexity. Conclusions are drawn in Section~\ref{sec:conclusions}.

%Finally, Section~\ref{sec:projectivespace} establishes the connection between the proposed infinity-crossing procedure and the extension of ODE flows to projective space~\cite{wright1985finding, kevrekidis2017infinity}.

%%%%%%%%%%%%%%%%%%%%%%%%
\section{Continuous Newton methods}
\label{sec:continewton}
The starting point for the CNM is the nonlinear system in Eq.~\eqref{eq:system_nonlin_eqs}, which we assume to be continuously differentiable, i.e., $\bm{g} \in \mathcal{C}^1(\mathbb{R}^n, \mathbb{R}^n)$.
The CNM generalizes the classical discrete Newton iteration in Eq.~\eqref{eq:discrete_newton} by embedding it into a continuous-time framework. In its most general form, this gives rise to the preconditioned initial value problem
\begin{equation}
    \dot{\bm{x}}(t) = \mathbf{M}(\bm{x}(t))\, \bm{g}(\bm{x}(t)), 
    \qquad \bm{x}(0) = \bm{x}_0,
    \label{eq:continuous_newton_general}
\end{equation}
where $\mathbf{M}(\bm{x}(t)) \in \mathbb{R}^{n \times n}$ is a preconditioning matrix whose choice determines the convergence properties of the flow.
In this work, we adopt the Branin formulation \cite{branin1972widely}, in which the preconditioner matrix takes the form
\begin{equation}
    \mathbf{M}(\bm{x}(t)) = \gamma(t) \big(\nabla \bm{g}(\bm{x}(t))\big)^{-1},
    \qquad \gamma(t) \in \{-1, +1\}.
    \label{eq:branin_precond}
\end{equation}
Substituting this choice into the general continuous framework in Eq.~\eqref{eq:continuous_newton_general} yields the \emph{Branin-Davidenko-Newton flow}
\begin{equation}
    \dot{\bm{x}}(t) = \bm{f}(t, \bm{x}(t)) := \gamma(t) 
    \big(\nabla \bm{g}(\bm{x}(t))\big)^{-1} \bm{g}(\bm{x}(t)), 
    \qquad \bm{x}(0) = \bm{x}_0.
    \label{eq:continuous_newton}
\end{equation}
Note that the stationary points of Eq.~\eqref{eq:continuous_newton} correspond to roots of Eq.~\eqref{eq:system_nonlin_eqs} whenever $\bm{g}$ is defined and continuously differentiable ($\bm{g} \in C^1$) in the domain of interest and the Jacobian $\nabla \bm{g}(\bm{x})$ is nonsingular\footnote{If $\bm{g}$ is not $C^1$, the Jacobian may become unbounded at some point, giving rise to spurious stationary behavior in the Newton flow. For example, the scalar function $g(x) = x^{1/3} + 1$ has a root at $x = -1$, but $g'(x) = \frac{1}{3}x^{-2/3}$ is unbounded at $x = 0$. The CNM flow $\dot{x} = -g(x)/g'(x) = -3x^{2/3}(x^{1/3} + 1)$ vanishes at $x = 0$, making it a spurious stationary point that traps trajectories starting from the right, preventing convergence to the root at $x = -1$.}. Indeed, $\bm{g}(\bm{x}) = \bm{0}$ implies $\dot{\bm{x}} = \bm{0}$. Conversely, if $\dot{\bm{x}} = \bm{0}$, then from Eq.~\eqref{eq:continuous_newton} we have $\big(\nabla \bm{g}(\bm{x})\big)^{-1}\bm{g}(\bm{x}) = \bm{0}$; because $\nabla \bm{g}(\bm{x})$ is invertible, this is equivalent to $\bm{g}(\bm{x}) = \bm{0}$. %Thus, the stationary points of the flow are exactly the roots of $\bm{g}$.

It can be shown that under this Branin-Davidenko formulation, the domain of guaranteed convergence for $\gamma(t) \equiv -1$ coincides with that given by the Newton--Kantorovich theorem~\cite{kelley1995iterative}.
Moreover, if $\gamma(t) \equiv -1$ and Eq.~\eqref{eq:continuous_newton} is discretized by the forward Euler method with unit step size $\dd t = 1$, one recovers the classical Newton iteration in Eq.~\eqref{eq:discrete_newton}, which may thus be viewed as a special case of this broader continuous framework.

The sign factor $\gamma(t)$ is the key mechanism by which the proposed methodology controls the orientation of the flow across singularities and passages through infinity.
As discussed in the Introduction, the orientation of the flow must be carefully maintained across singularities and passages through infinity to ensure that the trajectory continues to advance toward a potential solution rather than retracing its path. The choice $\gamma(t) = -1$ corresponds to the standard convergence direction of the continuous Newton method; departures from this choice are necessary when the trajectory encounters a flow singularity or a passage through 
infinity. In such cases, $\gamma(t)$ is updated locally to preserve the orientation of the flow, as described in detail in Sections~\ref{sec:scalarodes}--\ref{sec:algorithm}.

%\color{blue!90!black}
%
It can be shown (see, e.g.,\cite{branin1972widely}) that %if $\gamma(t) \equiv -1$ and
if $\nabla \bm{g}(\bm{x}(t))$ remains nonsingular for all $t \geq t'$ for some $t' \geq 0$, then the trajectory of Eq.~\eqref{eq:continuous_newton} satisfies
\begin{equation}
\dot{\bm{g}}(\bm{x}(t)) -\gamma(t) \bm{g}(\bm{x}(t)) = \bm{0},
\qquad t \geq t',
\label{eq:branin_ode}
\end{equation}
where $\dot{\bm{g}} = \frac{\dd \bm{g}}{\dd t}$. This follows directly from the chain rule:
\begin{equation}
\dot{\bm{g}}(\bm{x}(t)) = \nabla\bm{g}(\bm{x}(t)) \cdot \dot{\bm{x}}(t)
= \gamma(t)\nabla\bm{g}(\bm{x}(t)) \cdot \big(\nabla\bm{g}(\bm{x}(t))\big)^{-1} \bm{g}(\bm{x}(t))
= \gamma(t)\bm{g}(\bm{x}(t)).
\label{eq:branin_derivation}
\end{equation}
The analytical solution of Eq.~\eqref{eq:branin_ode}, for $\gamma(t)\equiv -1$, is
\begin{equation}
\bm{g}(\bm{x}(t)) = e^{(t' - t)} \, \bm{g}(\bm{x}(t')),
\qquad t \geq t'.
\label{eq:branin_solution}
\end{equation}
Hence $|\bm{g}(\bm{x}(t))| \to 0$ exponentially as $t \to \infty$, establishing that the Newton flow is a contraction in the residual. If, in addition, $\lim_{t \to \infty} \bm{x}(t)$ exists, then this limit is a solution of~\eqref{eq:system_nonlin_eqs}.

Note that the derivation above requires the Jacobian $\mathbf{J}:=\nabla \bm{g}(\bm{x}(t))$ to be invertible along the entire trajectory.
If there exists a time $t^{\ast} > t'$ at which $\det\big(\mathbf{J}(\bm{x}(t^\ast))\big) = 0$, then the right-hand side of Eq.~\eqref{eq:branin_derivation} becomes undefined, and the contraction property in Eq.~\eqref{eq:branin_solution} ceases to hold beyond $t^\ast$.
Nevertheless, on any interval $(t', t^\ast)$ where the Jacobian remains nonsingular, each component of the residual $|\bm{g}(\bm{x}(t))|$ decreases monotonically.  
%\color{black}

Furthermore, as detailed in~\cite{garcia1980relations, zufiria2002role}, Branin's continuous Newton method admits a natural interpretation within the framework of homotopy and continuation methods. Indeed, from the explicit solution in Eq.~\eqref{eq:branin_solution} for $t'=t_0$ and defining the homotopy parameter
\begin{equation}
\lambda(t) := e^{t_0 - t} \in \mathbb{R},
\end{equation}
the trajectory $\bm{x}(t)$ satisfies
\begin{equation}
\bm{G}(\bm{x}(t), \lambda(t)) := \bm{g}(\bm{x}(t)) - \lambda(t) \bm{g}(\bm{x}_0) = 0.
\end{equation}
Thus, in the absence of sign changes, the CNM flow corresponds to tracing the zero set of the global homotopy \(\bm{G}(\bm{x}, \lambda) = 0\) with \(\lambda\) varying monotonically from \(\lambda(0)=1\) toward \(0\) (for \(\gamma \equiv -1\)) or from \(1\) toward \(+\infty\) (for \(\gamma \equiv +1\)). When a singularity is encountered, Branin's method deliberately reverses the sign of \(\gamma(t)\), which in the homotopy picture corresponds to reversing the direction of \(\lambda\) along the path. This is precisely analogous to the pseudo-arclength continuation technique, where the sign of \(\dd\lambda/\dd s\) is reversed at limit points to traverse the solution manifold past critical points~\cite{allgower2012numerical}. Branin's sign-switching rule thus provides a heuristic, yet dynamically motivated, mechanism for navigating singularities by reversing the parametrization direction.
This connection with bifurcation theory is particularly instructive: points where \(\det(\nabla \bm{g}(\bm{x})) = 0\) correspond to turning/fold points of the global homotopy path.
Failure of the method, as for example in the Brent counterexample~\cite{brent1972davidenko}, due to looping between singularities, might correspond to closed, detached components of the solution manifold (sometimes called \emph{isolas}), i.e.\ regions in which solutions exist but are disconnected from the rest of the solution set.

\subsection{Illustrative example: discrete versus continuous Newton}
\label{subsec:example_disc_cont}
\begin{figure}[ht!]
    \centering
    \subfigure[]{\includegraphics[width=0.4\linewidth]{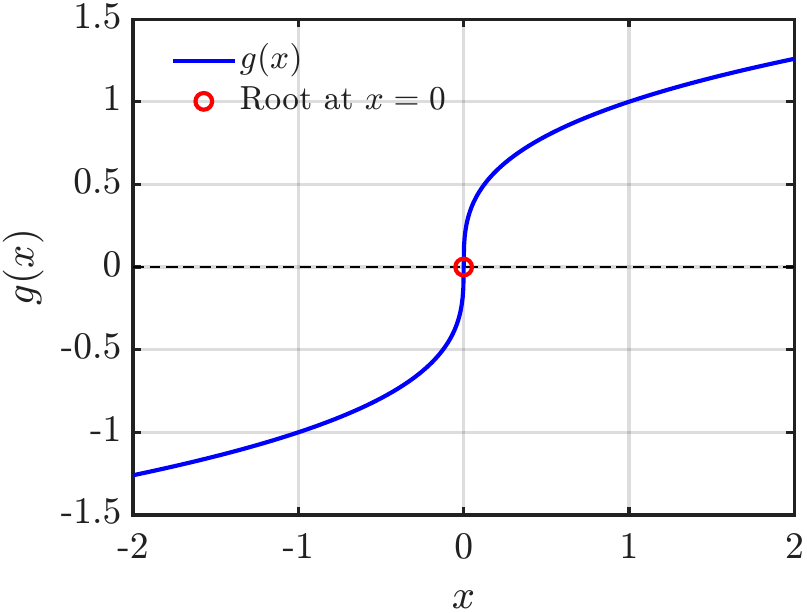}
    }
    \subfigure[]{\includegraphics[width=0.4\linewidth]{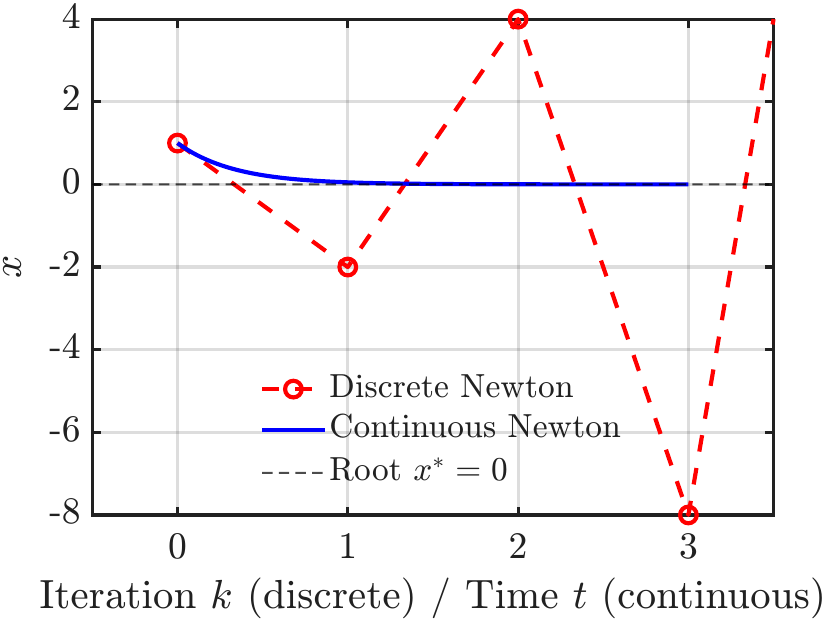}
    }
    \caption{\scriptsize Discrete versus continuous Newton for $g(x) = x^{1/3}$. 
(a) The function $g(x)$ (blue) with root at $x^*=0$ (red dot). 
(b) Evolution from $x_0=1$: discrete Newton (red) gives $x_k=(-2)^k x_0$, 
diverging to infinity, while continuous Newton (blue) gives 
$x(t)=x_0 e^{-3t}$, converging to $0$. \label{fig:xpower1_3}
}
\vspace{0.1cm}\hrule
\end{figure}

The advantage of the continuous formulation over the discrete Newton iteration can be illustrated by the scalar function 
$g(x) = x^{\alpha}$ (taking the real power, defined for all real 
$x$, more precisely $g(x) =  \sign(x)|x|^{\alpha}$), whose unique real root is $x^\ast = 0$ and whose derivative $g'(x) = \alpha |x|^{\alpha-1}$ is unbounded at the root for $\alpha<1$.
For any nonzero initial guess $x_0 \neq 0$, the classical Newton iteration yields
\begin{equation}
    x_{k+1} = x_k - \frac{x_k^{\alpha}}{\alpha x_k^{\alpha-1}} 
    = \left(1-\frac{1}{\alpha}\right)x_k,
    \label{eq:newton_x13}
\end{equation}
so that $x_k = (1-\frac{1}{\alpha})^k x_0$.
The iteration converges if $|1 - 1/\alpha| < 1$, i.e., $\alpha > 1/2$, and diverges otherwise. For $\alpha = 1/3$, the factor is $-2$, giving the oscillatory divergence shown in Fig.~\ref{fig:xpower1_3}.
The iterates oscillate in sign and grow exponentially in magnitude, diverging to infinity regardless of the choice of $x_0$. This failure is a direct consequence of the 
unboundedness of $g'(x)$ at the root, which violates the assumptions of the Newton--Kantorovich theorem.

In contrast, taking $\gamma(t) = -1$ in the continuous formulation in Eq.~\eqref{eq:continuous_newton} yields
\begin{equation}
    \dot{x}(t) = -\frac{(x(t))^{\alpha}}{\alpha (x(t))^{\alpha-1}} 
    = -\frac{1}{\alpha}x(t),
    \label{eq:cont_newton_x13}
\end{equation}
whose unique solution is $x(t) = x_0 e^{-\frac{1}{\alpha} t}$, which converges to $x^\ast = 0$ for all $\alpha > 0$ and $x_0 \in \mathbb{R}$ as $t \to \infty$. This example illustrates that the continuous Newton method may converge globally in situations where the discrete iteration fails entirely (see Fig.~\ref{fig:xpower1_3} for $\alpha = 1/3$).

%%%%%%%%%%%%%%%%%%%%%%%%%%%%%%
\subsection{Classification of singularities along the trajectory}
\label{subsec:classification}
Despite the improved robustness illustrated above, the numerical integration of Eq.~\eqref{eq:continuous_newton} may encounter difficulties whenever the Jacobian $\nabla \bm{g}(\bm{x}(t))$ becomes singular along the trajectory. We say that $\bm{x}(t)$ or $\bm{f}(t, \bm{x}(t))$ \emph{approach infinity} if the absolute value of one or more of their components grows without bound, i.e., $\bm{x}(t)$ approaches infinity if
\begin{equation}
    \exists\, i \in \{1, \ldots, n\} \;\text{and}\; 
    \exists\, t^\ast \in (0, \infty] \;\text{ such that }\; 
    \lim_{t \to t^\ast} |x_i(t)| = +\infty.
    \label{eq:infinity_def}
\end{equation}
If $t^\ast < \infty$ we say that $\bm{x}(t)$ approaches infinity in \emph{finite time} (\emph{blow-up}); if $t^\ast = \infty$ we say that it approaches infinity in \emph{infinite time} (\emph{grow-up}). In particular, one or more of the following situations may arise along the trajectory of Eq.~\eqref{eq:continuous_newton}:
\begin{enumerate}%[label=(\roman*)]
    \item[\texttt{IT1}]
    \emph{(Infinity type 1 - Blow-up)} $f_i(t, \bm{x}(t))$ approaches infinity in finite time, for at least one index $i \in \{1,\ldots,n\}$;
    \item[\texttt{IT2}]
    \emph{(Infinity type 2 - Grow-up)} $x_i(t)$ itself approaches infinity in infinite time, for at least one index $i \in \{1,\ldots,n\}$.
\end{enumerate}
%
%\begin{remark}
%\label{rem:types_prime}
One might also consider the cases where the state $x_i(t)$ approaches infinity in finite time (\emph{infinity type $1'$} (\texttt{IT}--$1'$) or $f_i(t,\bm{x}(t))$ approaches infinity in infinite time (\emph{infinity type $2'$} (\texttt{IT}--$2'$)). 
Under Assumption~\ref{ass:monotonicity} below, \texttt{IT}--$1'$ implies infinity type~\texttt{IT1}, and \texttt{IT}--$2'$ implies infinity type~\texttt{IT2}. Hence, without loss of generality, it suffices to consider only the conditions for~\texttt{IT1} and~\texttt{IT2}.
%\end{remark}
%
\begin{assumption}
\label{ass:monotonicity}
If $x_i(t)$ approaches infinity at $t = t^\ast$, then $\dot{x}_i(t)$ is monotone on $(t^\ast - \varepsilon, t^\ast)$ for some $\varepsilon > 0$.
\end{assumption}
\begin{lemma}
\label{lem:type1prime}
Under Assumption~\ref{ass:monotonicity}, \texttt{IT}--$1'$ implies ~\texttt{IT1}.
\end{lemma}
\begin{proof}
Without loss of generality, suppose $x_i(t) \to +\infty$ as $t \to t^\ast$ with $t^\ast < \infty$. Since $x_i(t)$ is continuous and unbounded on $(0, t^\ast)$, there exists an ascending sequence $t_1 < t_2 < \cdots$ with $t_k \nearrow t^\ast$ such that $x_i(t_{k+1}) - x_i(t_k) = 1$ for all $k \in \mathbb{N}$ and $t_{k+1} - t_k \to 0$. By the mean value theorem, for each $k \in \mathbb{N}$ there exists $\tau_k \in (t_k, t_{k+1})$ such that
\begin{equation}
    f_i(\tau_k, \bm{x}(\tau_k)) = \frac{x_i(t_{k+1}) - x_i(t_k)}
    {t_{k+1} - t_k} = \frac{1}{t_{k+1} - t_k} \to +\infty 
    \quad \text{as } k \to \infty.
    \label{eq:mvt_blowup}
\end{equation}
Since $\tau_k \nearrow t^\ast$ and $\dot{x}_i(t) = f_i(t, \bm{x}(t))$ is monotone on $(t^\ast - \varepsilon, t^\ast)$ by Assumption~\ref{ass:monotonicity}, it follows that
\begin{equation}
    \lim_{t \to t^\ast} f_i(t, \bm{x}(t)) = +\infty,
\end{equation}
establishing~\texttt{IT1}.
\end{proof}
\begin{remark}
\label{rem:converse}
The converse of Lemma~\ref{lem:type1prime} does not hold in general: it is possible that $f_i(t, \bm{x}(t))$ diverges in finite time while $x_i(t)$ remains bounded, as will indeed be observed in the illustrative example of Section~\ref{subsec:blowup_example} and in several of the numerical experiments of Section~\ref{sec:casestudies}.
\end{remark}
Similarly, \texttt{IT2} does not necessarily imply \texttt{IT}--$2'$: for instance, if $f_i \equiv 1$, then $x_i(t) = x_i(0) + t \to \infty$ as $t \to \infty$ (\texttt{IT2}), while $f_i$ remains bounded (so \texttt{IT}--$2'$ does not occur). The direct implication, however, is straightforward. 
\begin{lemma}
\label{lem:type2prime}
\texttt{IT}--$2'$ implies \texttt{IT2}.
\end{lemma}
\begin{proof}
If $f_i(t,\bm{x}(t)) \to \infty$ as $t \to \infty$, then for any $M > 0$ there exists $T > 0$ such that $f_i(t,\bm{x}(t)) > M$ for all $t > T$. Integrating from $T$ to $t$ gives $x_i(t) > x_i(T) + M(t - T) \to \infty$ as $t \to \infty$, establishing \texttt{IT2}.
\end{proof}

\subsection{An illustrative example: finite-time blow-up and multiple singularity crossings}
\label{subsec:blowup_example}
To illustrate the condition~\texttt{IT1} and the role of the orientation factor $\gamma(t)$, consider the scalar function
\begin{equation}
    g(x) = \frac{x^3 - 3x}{2} + 2,
    \label{eq:blowup_example}
\end{equation}
whose derivative $g'(x) = \frac{3(x^2-1)}{2}$ vanishes at $x = \pm 1$, corresponding to a local maximum at $x = -1$ with $g(-1) = 3 > 0$ and a local minimum at $x = 1$ with $g(1) = 1 > 0$. Since both stationary points of $g$ lie strictly above zero, the function has a unique real root $x^\ast \approx -2.1958$, satisfying $x^3 - 3x + 4 = 0$, located to the left of the local maximum.
\begin{figure}
    \centering
    \subfigure[]{\includegraphics[width=0.45\linewidth]{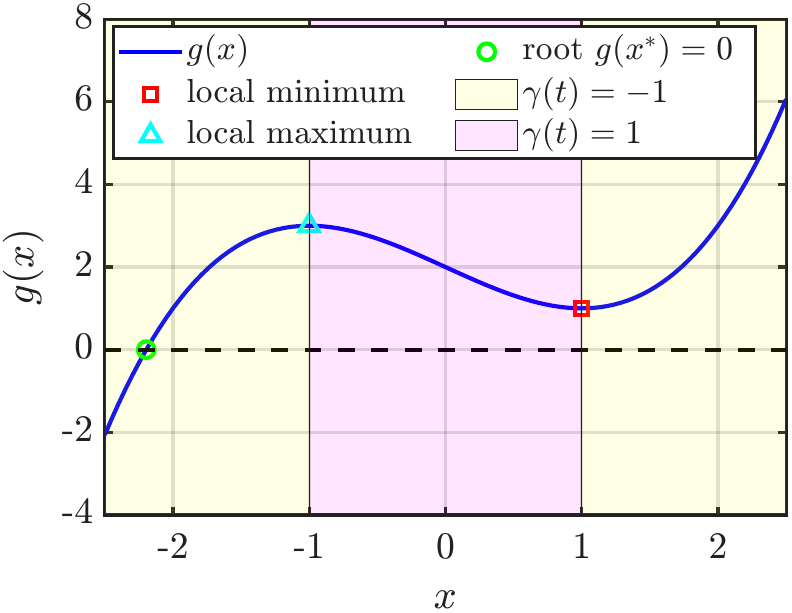}
    }
    \subfigure[]{\includegraphics[width=0.45\linewidth]{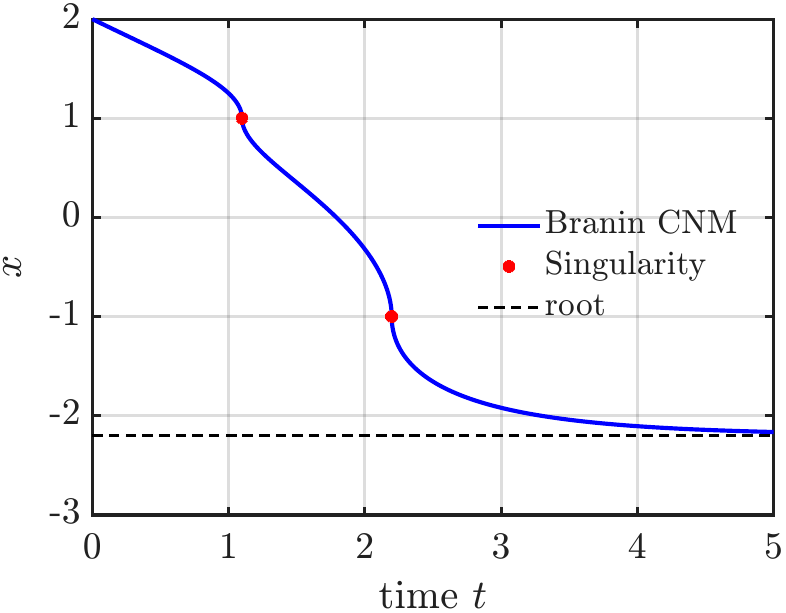}
    }
    \subfigure[]{\includegraphics[width=0.45\linewidth]{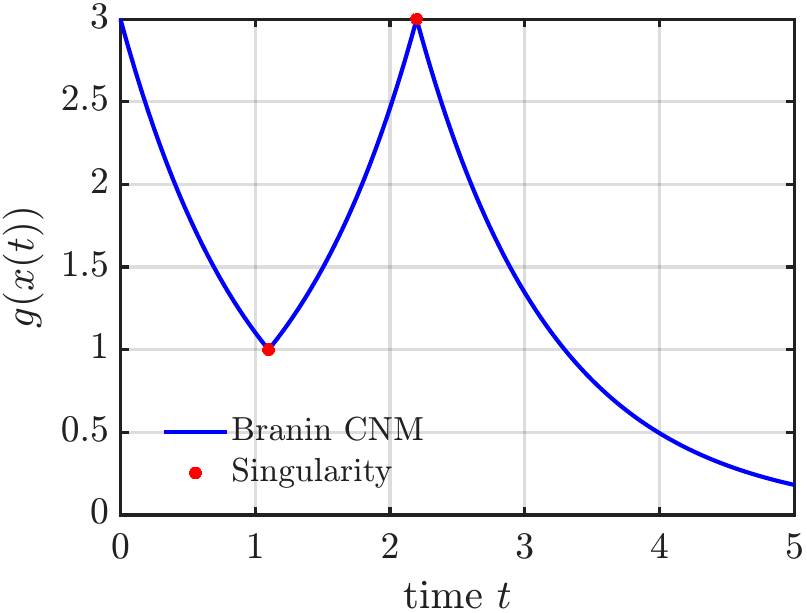}
    }
    \subfigure[]{\includegraphics[width=0.45\linewidth]{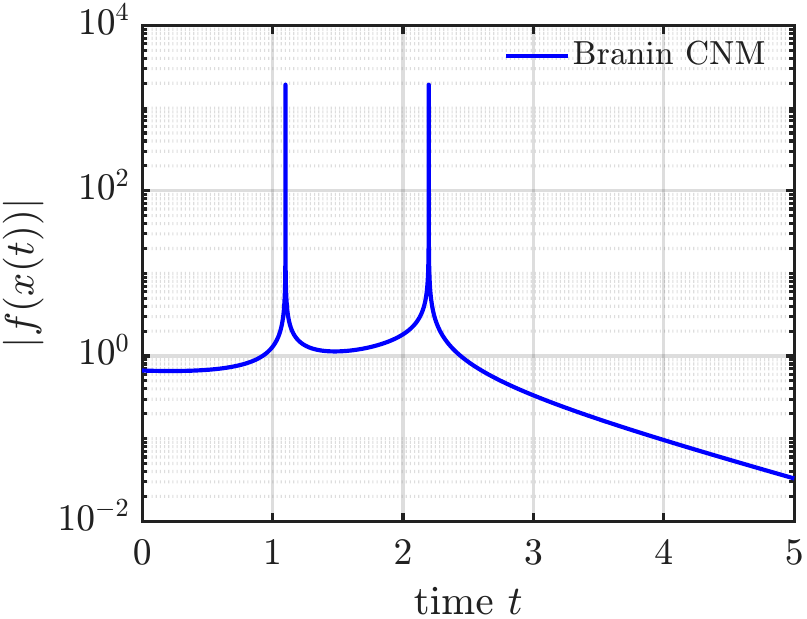}
    }
    
    \caption{\scriptsize Numerical demonstration of finite-time blow-up and multiple Jacobian singularity crossings (infinity type~1 (\texttt{IT1})) for the scalar example in Eq.~(\ref{eq:blowup_example}). Panel (a) displays $g(x)$ with its local maximum, local minimum, and unique root. Panel (b) shows the trajectory $x(t)$ of the Branin CNM, while panels (c) and (d) depict the evolution of $g(x(t))$ and the logarithmic evolution of $f(x(t))$, respectively; the vertical asymptotes in (d) confirm the finite-time blow-up of the vector field at $x = \pm 1$. \label{fig:ex2}
    }
    \vspace{0.1cm}\hrule
\end{figure}

Consider the CNM flow in Eq.~\eqref{eq:continuous_newton} with $\gamma(t) = -1$ and initial condition $x(0) = 2 \in (1, +\infty)$. 
In this region we have $g(x) > 0$ and $g'(x) > 0$, so that $\dot{x} = -g/g' < 0$ and the trajectory decreases toward the Jacobian singularity at $x = 1$. As $x(t) \to 1^+$, we have $g'(x) \to 0^+$ while $g(1) = 1 > 0$, so that
\begin{equation}
    \dot{x}(t) = -\frac{g(x(t))}{g'(x(t))} \to -\infty 
    \quad \text{as } x(t) \to 1^+,
\end{equation}
and the vector field blows-up in finite time, providing a concrete instance of \texttt{IT1}.
Upon crossing this singularity (e.g., as described in Section~\ref{sec:scalarodes}), the orientation is updated to $\gamma = +1$ to preserve the direction of the flow.
In the interval $(-1, 1)$ we have $g'(x) < 0$, so the flow now advances toward the local maximum at $x = -1$, where a second Jacobian singularity is encountered.
After crossing this second singularity with another update of $\gamma$, the trajectory enters the region $(-\infty, -1)$, where $g'(x) > 0$ and $g(x)$ changes sign, and the flow converges to the unique root $x^\ast \approx -2.1958$.
This behavior, including the two finite-time blow-ups and the subsequent convergence to the root, is illustrated in Fig.~\ref{fig:ex2}. As the title of the paper indicates, the diverging solution trajectory has been ``salvaged".

%%%%%%%%%%%%%%%%%%%%%
\section{Infinity crossing for scalar ODEs}
\label{sec:scalarodes}

We now develop the infinity-crossing procedure for the two types of singularities (\texttt{IT1} and \texttt{IT2}) defined in Section~\ref{subsec:classification}. For ease of presentation, we first consider the scalar case $x(t) \in \mathbb{R}$; the extension to $n$-dimensional systems is given in Section~\ref{sec:multidimensionalodes}.

Throughout this section, we assume that Eq.~\eqref{eq:continuous_newton} is scalar, so that the vector field reduces to
\begin{equation}
    f(x(t)) := \gamma(t) \frac{g(x(t))}{g'(x(t))},
    \label{eq:scalar_vf}
\end{equation}
where $\gamma(t) \in \{-1, +1\}$ and $g'$ is the first derivative of $g$. 
%Without loss of generality, we assume throughout that
%
% \begin{equation}
%     \lim_{t \uparrow t^\ast} f(x(t)) = +\infty;
%     \label{eq:finfty}
% \end{equation}
% %
% the case $f(x(t)) \to -\infty$ is handled analogously by symmetry.

\subsection{Infinity crossing in finite time (\texttt{IT1})}
\label{subsec:infinity1}
We investigate the case in which $f(x(t))$ blows-up in finite time at some $t^\ast < \infty$, that is, the condition \texttt{IT1}.

\subsubsection{Self-similar model for the blow-up}
\label{subsubsec:selfsimilar}

The foundation of the proposed crossing procedure is the assumption that the vector field admits a self-similar representation in a neighborhood of the blow‑up time \(t^\ast\).

\begin{assumption}
\label{ass:selfsimilar}
In a neighborhood of the blow-up time \(t^\ast\), the vector field satisfies
\begin{equation}
    |f(x(t))| = a \cdot |t^\ast - t|^q,
    \label{eq:selfsimilar}
\end{equation}
where \(a > 0\) is an amplitude constant and \(q < 0\) is the self-similarity exponent.
\end{assumption}

The self-similar ansatz in Eq.~\eqref{eq:selfsimilar}, motivated by~\cite{kevrekidis2017infinity}, postulates that the vector field behaves as a power law in time along the trajectory near \(t^\ast\).
Further considerations and validations of this assumption are considered in the next subsection.

The \textit{self-similarity exponent} \(q\) plays a central role in distinguishing two qualitatively different scenarios:
\begin{enumerate}[label=(\roman*)]
    \item \label{it:qsmall} \(-1 < q < 0\): the vector field \(f(x(t))\) diverges in finite time while the state \(x(t)\) remains bounded, since the integral
    \begin{equation}
        x(t^\ast + h) - x(t^\ast - h) = 
        \int_{t^\ast-h}^{t^\ast+h} f(x(t))\, \dd t 
        \approx \int_{t^\ast-h}^{t^\ast+h} 
        a|t^\ast-t|^q\, \dd t = \frac{2ah^{q+1}}{q+1}
        \label{eq:integral_finite}
    \end{equation}
    is finite for \(q > -1\);
    \item \label{it:qlarge} \(q \leq -1\): both the vector field and the state \(x(t)\) diverge in finite time, since the integral above diverges. 
\end{enumerate}
This distinction determines the crossing strategy adopted in Subsection~\ref{subsubsec:crossing}.

In practice, the parameters \((a,t^\ast,q)\) are not known a priori; they are estimated on the fly from a short sequence of threshold crossings during the integration (see Section~\ref{subsubsec:detection}). The ansatz therefore serves as a local parametrization model whose parameters are inferred {\em a posteriori}, allowing the method to adapt to the observed blow‑up behavior without requiring prior knowledge of the singularity type.

%%%%%%%%%%%%%%%%%%%%%%%%%%%
\subsubsection{Validation of the self-similarity assumption and implications for the Continuous Newton flow}
\label{subsec:validation}
We now justify the self-similar ansatz in Assumption~\ref{ass:selfsimilar} for the continuous Newton flow~\eqref{eq:scalar_vf} encountering a singularity in finite time $t^\ast$. For the Newton flow, assuming $g$ is bounded, the only points where the vector field can become unbounded are points $x^\ast$ such that $g'(x^\ast)=0$ (since $g$ remains finite).

The following lemma shows that for the scalar continuous Newton flow, the self-similar ansatz is exact whenever the trajectory approaches a finite point where \(g'\) has a finite‑order zero (multiplicity $m$) while \(g\) is nonzero. Such points correspond precisely to critical points of \(g\) — i.e., local minima, local maxima, or inflection points — where the derivative vanishes but the residual does not.

\begin{lemma}[Power-law behavior of the Newton's Flow near a local extremum]
\label{lem:m_zero_taylor}
Let \(g \in C^{m+2}(\mathbb{R})\) for some integer \(m \ge 1\), and suppose \(x^\ast \in \mathbb{R}\) satisfies
\begin{equation}
g'(x^\ast) = g''(x^\ast) = \dots = g^{(m)}(x^\ast) = 0,\qquad g^{(m+1)}(x^\ast) \neq 0,
\end{equation}
and \(g(x^\ast) \neq 0\). Consider the continuous Newton flow
\begin{equation}
\dot{x} = f(x) = \gamma\frac{g(x)}{g'(x)},
\end{equation}
with a constant \(\gamma \in \{-1,1\}\). If the solution \(x(t)\) approaches \(x^\ast\) as \(t \to t^\ast\) for some finite \(t^\ast\), then there exist constants \(a > 0\) and \(q = -\frac{m}{m+1}\) such that
\begin{equation}
|f(x(t))| \sim a\,(t-t^\ast)^q \qquad \text{as } t \to t^\ast.
\end{equation}
\end{lemma}
\begin{proof}
We begin by noting that, although the flow \( f = g / g' \) is singular at \( x^\ast \), the individual functions \( g \) and \( g' \) are smooth in a neighborhood of \( x^\ast \). Hence, each admits a Taylor expansion about \( x^\ast \).
Therefore, let us expand \( g \) and \( g' \) around \( x^\ast \) using Taylor's theorem:
\begin{equation}
g(x) = g(x^\ast) + \frac{g^{(m+1)}(x^\ast)}{(m+1)!}(x-x^\ast)^{m+1} + O\big((x-x^\ast)^{m+2}\big),
\end{equation}
\begin{equation}
g'(x) = \frac{g^{(m+1)}(x^\ast)}{m!}(x-x^\ast)^m + O\big((x-x^\ast)^{m+1}\big).
\end{equation}
Thus
\begin{equation}
f(x) = \gamma\frac{g(x)}{g'(x)} = \gamma\frac{g(x^\ast)}{\frac{g^{(m+1)}(x^\ast)}{m!}(x-x^\ast)^m} + o\big((x-x^\ast)^{-m}\big) = C\,(x-x^\ast)^{-m},
\end{equation}
with \(C = \gamma \, m!\,g(x^\ast) / g^{(m+1)}(x^\ast) \neq 0\). Near \(x^\ast\) the differential equation therefore behaves like \(\dot{x} = C (x-x^\ast)^{-m}\). Separating variables and integrating from \(t\) to \(t^\ast\) gives
\begin{equation}
\int_{x(t)}^{x^\ast} (y-x^\ast)^m\,\dd y = C \int_t^{t^\ast} \dd \tau
\quad\Longrightarrow\quad
\frac{(x^\ast-x(t))^{m+1}}{m+1} = C\,(t^\ast-t).
\end{equation}
Hence
\begin{equation}
x(t)-x^\ast \sim \big((m+1)C\big)^{1/(m+1)}\,(t-t^\ast)^{1/(m+1)}.
\end{equation}
Substituting back,
\begin{equation}
f(x(t)) \sim C \big((m+1)C\big)^{-m/(m+1)}\,(t-t^\ast)^{-m/(m+1)}.
\end{equation}
Therefore \(|f(x(t))| \sim a\,(t-t^\ast)^q\) with \(q = -\frac{m}{m+1}\) and
\(a = |C|\,|(m+1)C|^{-m/(m+1)}\).
\end{proof}
\begin{remark}
As noted in the proof of Lemma~\ref{lem:m_zero_taylor}, the smoothness (or analyticity) of $g$ and $g'$ near $x^\ast$ guarantees their convergent Taylor expansions. Since $g'$ vanishes to order $m$ while $g(x^\ast)\neq0$, the vector field $f = \gamma g/g'$ is meromorphic at $x^\ast$ with a pole of order $m$. Consequently, the expansion $f(x) = C (x-x^\ast)^{-m} + O\big((x-x^\ast)^{-m+1}\big)$ is exact in the sense of a Laurent series, and the leading term dominates as $x\to x^\ast$. %Hence the power-law ansatz in time is asymptotically exact.
\end{remark}

An immediate consequence of Lemma~\ref{lem:m_zero_taylor} is that for every local extremum or inflection point encountered by the Newton flow at a finite $x^\ast$, the self-similarity assumption holds asymptotically in the limit $t\to t^\ast $, and the exponent satisfies $q = -m/(m+1) \in (-1,0)$, i.e., strictly greater than $-1$.
Therefore, whenever estimation yields $q < -1$, the trajectory cannot be approaching a finite point $x^\ast$; it must be escaping to infinity. This occurs, for instance, when $g(x)$ approaches a nonzero constant horizontal asymptote as $|x|\to\infty$.

To model such behavior, consider the family
\begin{equation}
g(x) = c + \frac{d}{x^\beta}, \qquad \beta>0,\; c\neq0,\; d\neq0,
\label{eq:asymptote_family}
\end{equation}
which satisfies $g(x)\to c$ and $g'(x)\to 0$ as $x\to\infty$. For this family, a direct computation (which we summarize in Lemma~\ref{lem:asymptote} below) yields $q = -1 - 1/\beta < -1$, showing that any exponent $q<-1$ can be realized.

\begin{lemma}[Blow-up at a horizontal asymptote]
\label{lem:asymptote}
Let $g(x) = c + d x^{-\beta}$ with $\beta> 0$, $c\neq0$, $d\neq0$, and consider the 
Newton flow $\dot{x} = \gamma g(x)/g'(x)$. Assume $C=\gamma c/(\beta d) > 0$.
Then as $x\to\infty$, the flow satisfies
\begin{equation}
\dot{x} \simeq C\, x^{\beta+1},
\end{equation}
and the solution blows-up in finite time $t^\ast$ with self-similar exponent
\begin{equation}
q = -1 - \frac{1}{\beta} < -1.
\end{equation}
\end{lemma}

\begin{proof}
For $g(x) = c + d \,x^{-\beta}$, we have $g'(x) = -\beta \,d \,x^{-(\beta+1)}$. Hence
\begin{equation}
f(x) = \gamma\frac{g(x)}{g'(x)} = -\gamma\frac{c + d \,x^{-\beta}}{\beta \, d \, x^{-(\beta+1)}} 
= C x^{\beta+1} - \gamma\frac{1}{\beta} x.
\end{equation}
For large $x$, the leading term is $C x^{\beta+1}$.
Solving $\dot{x} = C x^{\beta+1}$ gives $x(t) \sim (t-t^\ast)^{-1/\beta}$ and consequently $f(x(t)) \sim (t-t^\ast)^{-(\beta+1)/\beta}$, so $q = -(\beta+1)/\beta = -1 - 1/\beta$.
\end{proof}

Thus, the self-similar assumption in Assumption~\ref{ass:selfsimilar} is asymptotically valid for the two fundamental scenarios encountered by the continuous Newton flow: passage through a local extremum or inflection point of $g$ (giving $q\in(-1,0)$), and escape to a horizontal asymptote (giving $q<-1$).

%
%From the integral representation in Eq.~\eqref{eq:integral_finite}, the self-similar 
%ansatz directly implies a dichotomy: if $q>-1$ then the state $x(t)$ remains bounded 
%across the singularity, whereas if $q\leq -1$ then $x(t)$ itself blows-up in finite 
%time.
%

\subsubsection{Detection and estimation of blow-up parameters}
\label{subsubsec:detection}
Given Assumption~\ref{ass:selfsimilar}, the parameters $(q, t^\ast, a)$ are estimated on the fly from a sequence of threshold crossings. Let $f_{\max} \gg 1$ be a prescribed detection threshold. Detection is implemented via the event function
\begin{equation}
    E_{\textrm{fun}}(x) = |f(x)| - f_{\max},
    \label{eq:event}
\end{equation}
which triggers when $E_{\textrm{fun}}$ transitions from negative to positive. i.e., when $|f(x(t))|$ first exceeds $f_{\max}$. The sign of $f$ 
at the event time determines whether $f$ is approaching $+\infty$ or $-\infty$. In practice, this is implemented in \textsc{Matlab} via the \texttt{Events} option of the ODE suite (e.g., \texttt{ode45}, \texttt{ode15s}), with \texttt{isterminal = 1} and \texttt{direction = 1}.

To estimate the self-similarity parameters, the integration is continued past 
the first threshold crossing through a sequence of $N \geq 3$ successively larger thresholds:
\begin{equation}
    E^{(k)}_{\textrm{fun}}(x) = |f(x)| - c^{k-1} f_{\max}, \qquad k = 1, \ldots, N,
    \label{eq:event_sequence}
\end{equation}
where $c > 1$ is a fixed growth factor. The integration is advanced from one threshold to the next, recording the crossing times $t_1 < t_2 < \cdots < t_N$ at which $E_{\textrm{fun}}^{(k)}(x(t_k)) = 0$, i.e.,
\begin{equation}
    |f(x(t_k))| = c^{k-1} f_{\max} = a|t^\ast - t_k|^q, 
    \qquad k = 1, \ldots, N.
    \label{eq:threshold_crossings}
\end{equation}

\begin{itemize}
    \item \textit{Estimation of $q$.}
Solving Eq.~\eqref{eq:threshold_crossings} for $t^\ast - t_k$ gives
\begin{equation}
    t^\ast - t_k = c^{\frac{k-1}{q}} a^{-\frac{1}{q}} 
    f_{\max}^{\frac{1}{q}}, \qquad k = 1, \ldots, N.
    \label{eq:tstar_tk}
\end{equation}
Defining the time increments $\Delta t_k := t_{k+1} - t_k$ and subtracting consecutive terms in Eq.~\eqref{eq:tstar_tk}, one obtains
\begin{equation}
    \Delta t_k = \left(c^{\frac{k-1}{q}} - c^{\frac{k}{q}}\right) 
    a^{-\frac{1}{q}} f_{\max}^{\frac{1}{q}}, 
    \qquad k = 1, \ldots, N-1,
    \label{eq:deltat}
\end{equation}
from which it follows that
\begin{equation}
    \frac{\Delta t_{k+1}}{\Delta t_k} = c^{\frac{1}{q}}, 
    \qquad k = 1, \ldots, N-2.
    \label{eq:ratio}
\end{equation}
Taking logarithms, we obtain the linear relation
\begin{equation}
    \log(\Delta t_{k+1}) - \log(\Delta t_k) = 
    \frac{\log c}{q}, \qquad k = 1, \ldots, N-2.
    \label{eq:loglinear}
\end{equation}
%
%Each $k$ provides an estimate $q_k = \frac{\log c}{\log(\Delta t_{k+1}) - \log(\Delta t_k)}$.
Since the right-hand side of Eq.~\eqref{eq:loglinear} is constant, the least squares estimator for $\log (c)/q$ is the sample mean of the left-hand sides:
\begin{equation}
    \label{eq:telescoping}
    \frac{1}{N-2}\sum_{k=1}^{N-2} \log(\Delta t_{k+1}) - \log(\Delta t_k) = 
    \frac{\log c}{\hat q},
\end{equation}

By reducing the telescoping sum on the left-hand-side of Eq.~\eqref{eq:telescoping}, this reduces to
\begin{equation}
    q\simeq\hat{q} = \frac{(N-2)\log c}{\log(\Delta t_{N-1}) - \log(\Delta t_1)}.
    \label{eq:q_estimate}
\end{equation}

%\color{blue!90!black}
\begin{remark}
Note that $q$ is estimated \emph{indirectly}, rather than by directly averaging some individual local estimates, e.g., of type $\hat{q}_{(k)}= \frac{\log c}{\log(\Delta t_{k+1}) - \log(\Delta t_k)}$. Here, the assumption is that the source of error is in the computed time increments $\Delta t_k$, which carry numerical errors from the integration near the singularity.
Indeed, as the trajectory approaches a singularity, the adaptive solver computes increasingly smaller $\Delta t_k$, tending to 0. Although the detection threshold prevents these increments from reaching machine precision, they become increasingly contaminated by numerical errors arising from the growth of the flow and the ill-conditioning of the underlying linear systems. 
%These errors affect the log-ratios $\log(\Delta t_{k+1}) - \log(\Delta t_k)$. By formulating the estimation as a linear least squares problem in the log-domain, the telescoping estimator in Eq.~\eqref{eq:q_estimate} uses all available intervals in an optimal way, reducing the influence of the noisiest measurements without assuming any particular error distribution.
\end{remark}

% which is solved by linear regression to obtain the estimate
% %
% \begin{equation}
%     q = \frac{(N-2)\log c}{\log(\Delta t_{N-1}) - \log(\Delta t_1)}.
%     \label{eq:q_estimate}
% \end{equation}

\item \textit{Estimation of $t^\ast$.}
From Eq.~\eqref{eq:ratio}, the sequence $\Delta t_k$ is geometric with ratio $c^{1/q}$. Since $\hat{q} < 0$ and $c > 1$, we have $c^{1/\hat{q}} \in (0,1)$, so the series converges and the blow-up time is estimated as
\begin{equation}
    t^\ast\simeq \hat{t}^\ast = t_N + \sum_{j=1}^{\infty} c^{\frac{j}{\hat{q}}} \Delta t_{N-1} 
    = t_N + \frac{c^{1/\hat{q}}}{1 - c^{1/\hat{q}}} \Delta t_{N-1}.
    \label{eq:tstar_estimate}
\end{equation}

\item \textit{Estimation of $a$.}
Given $\hat{q}$ and $\hat{t}^\ast$, the amplitude $a$ is recovered from Eq.~\eqref{eq:selfsimilar} evaluated at $t = t_N$:
\begin{equation}
    a\simeq \hat{a} = \frac{c^{N-1} f_{\max}}{(\hat{t}^\ast - t_N)^{\hat{q}}}.
    \label{eq:a_estimate}
\end{equation}
\end{itemize}

\subsubsection{Crossing strategy}
\label{subsubsec:crossing}
With the parameters $(q, t^\ast, a)\simeq (\hat q, \hat t^\ast, \hat a)$ estimated, the crossing of the singularity proceeds differently according to the two cases identified in Subsection~\ref{subsubsec:selfsimilar}.

\begin{itemize}
    \item \textit{Case $-1 < q < 0$: vector field blows-up, state bounded.}
Since the integral in Eq.~\eqref{eq:integral_finite} is finite, the state $x(t)$ remains bounded across the singularity. Letting $h = t^\ast - t_N$ denote the distance from the last computed time point to the estimated blow-up time, we approximate $f(x(t))$ by 
the self-similar model in Eq.~\eqref{eq:selfsimilar} on the interval $[t^\ast - h, t^\ast + h]$ and integrate explicitly to obtain
\begin{equation}
    x(t^\ast + h) = x(t^\ast - h) + \int_{t^\ast-h}^{t^\ast+h} 
    f(x(t))\, \dd t \approx x(t^\ast - h) + \frac{2ah^{q+1}}{q+1}.
    \label{eq:crossing_formula}
\end{equation}
%
% The orientation $\gamma$ is updated at $t^\ast$ to preserve the direction of the flow:
% %
% \begin{equation}
%     \gamma(t^\ast + \varepsilon) = -\gamma(t^\ast - \varepsilon),
%     \label{eq:sigma_update}
% \end{equation}
% %
% and integration resumes from $(t^\ast + h,\, x(t^\ast + h))$.

\item \textit{Case $q \leq -1$: both vector field and state blow-up.}
When $q \leq -1$, the integral in Eq.~\eqref{eq:integral_finite} diverges and $x(t)$ itself blows-up in finite time. In this case, a direct integral correction is no longer available. Instead, motivated by the projective compactification of $\mathbb{R}$ discussed in~\cite{kevrekidis2017infinity}, we reset the state via
\begin{equation}
    x(t^\ast + h) = -x(t^\ast - h),
    \label{eq:sign_flip}
\end{equation}

which corresponds to identifying $+\infty$ and $-\infty$ as the same point and continuing the trajectory on the other side of the compactified real line.
%The orientation $\gamma$ is updated according to Eq.~\eqref{eq:sigma_update} in both cases, \color{blue!90!black} with the goal of preventing the trajectory from returning to the same singularity. Integration then resumes from $(t^\ast + h,\, x(t^\ast + h))$.
%\color{black}
\end{itemize}

\textbf{\textit{Orientation update.}}
In both cases, {in order to prevent the trajectory from returning to the same singularity}, after the crossing the orientation factor $\gamma$ is updated according to
%\color{blue!90!black}
%
\begin{equation}
    \gamma_{+} = \gamma_{-} \cdot 
    \sign\!\left(f(x_{-};\gamma_{-})\right) \cdot 
    \sign\!\left(f(x_{+};\gamma_{-})\right),
    \label{eq:sigma_update}
\end{equation}
where $\sign(\cdot) \in \{-1,+1\}$ denotes the sign function and for ease of notation we used $x_{\pm} := x(t^\ast \pm h)$, $\gamma_{\pm} := \gamma(t^\ast \pm h)$. Note that $f(\cdot;\gamma_{-})$ denotes evaluation using the pre-crossing orientation factor.

In contrast to Branin's original strategy~\cite{branin1972widely}, which switches the sign at every singularity crossing, our orientation update may yield $\gamma_{+} = \gamma_{-}$ when the sign of $g'(x)$ is preserved across the singularity (e.g., at inflection points), and $\gamma_{+} = -\gamma_{-}$ when the sign $g'(x)$ is effectively changed (e.g., at local extrema).

Indeed, with this choice in Eq.~\eqref{eq:sigma_update}, since $\sign(ab)=\sign(a)\sign(b)$ and $\sign(a)\sign(a)=1$, we observe that
\begin{equation}
%\small
\begin{split}
    \sign&(f(x_{+};\gamma_{+})) 
    = \sign(\gamma_{+}) \sign\!\left(\frac{g(x_{+})}{g'(x_{+})}\right) =\\
    &=  \sign(\gamma_{-}) \sign(\gamma_{-}) \sign\!\left(\frac{g(x_{-})}{g'(x_{-})}\right)
    \sign(\gamma_{-}) \sign\!\left(\frac{g(x_{+})}{g'(x_{+})}\right)
       \sign\!\left(\frac{g(x_{+})}{g'(x_{+})}\right) =\\
    &= \sign(\gamma_{-}) \sign\!\left(\frac{g(x_{-})}{g'(x_{-})}\right)=\sign(f(x_{-};\gamma_{-})).
\end{split}
\end{equation}
%
% Since $g(x)$ does not change sign across the singularity (where $g'(x^*)=0$ and $g(x^*)\neq 0$), we have $\sign(g(x_{+})) = \sign(g(x_{-}))$, so $\sign(g(x_{-}))\sign(g(x_{+})) = 1$. Thus
% \begin{equation}
%     \sign(f(x_{+};\gamma_{+})) 
%     = \sign(\gamma_{-}) \sign(g'(x_{-})) \sign(g'(x_{+})) 
%     = \sign(f(x_{-};\gamma_{-})),
% \end{equation}
% %
%This confirms that $\sign(f(x_{+};\gamma_{+}))=\sign(f(x_{-};\gamma_{-}))$ and the orientation of the flow is preserved across the crossing.
%\color{black}
Therefore, this rule preserves the orientation of the flow when $f$ has the same sign on both sides of the singularity — as occurs for instance at an inflection point of $g$, e.g., $g(x) = x^3 + 1$ with singularity at $x = 0$ — and reverses it when $f$ changes sign, as occurs at a local minimum or maximum of $g$. Integration then resumes from $(t^\ast + h,\, x(t^\ast + h))$.

%%%%%%%%%%%%%%%%%%
\subsubsection{Illustrative example: estimation of blow-up parameters}
\label{sec:quadratic_sqrt2}

\begin{figure}[ht!]
\centering
    \subfigure[]{
    \includegraphics[width=0.3\textwidth]{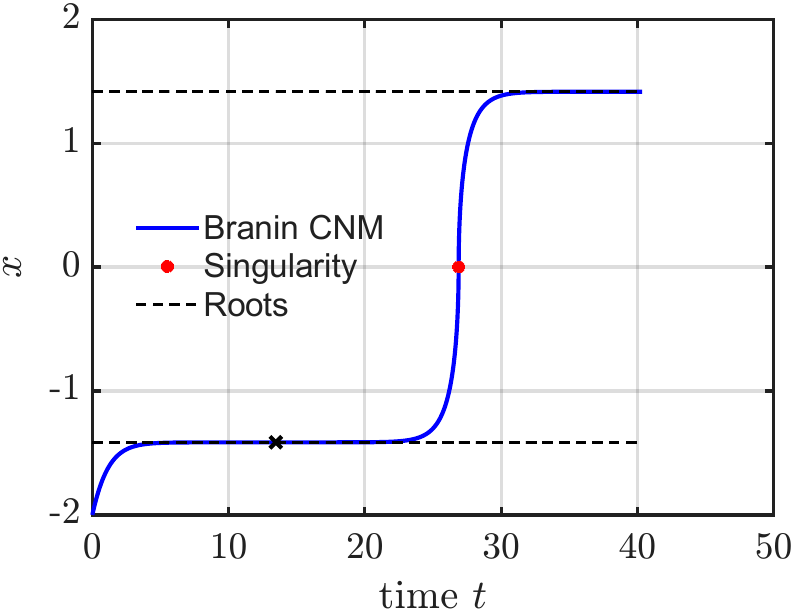}
    }
    %\hfill
    \subfigure[]{
    \includegraphics[width=0.3\textwidth]{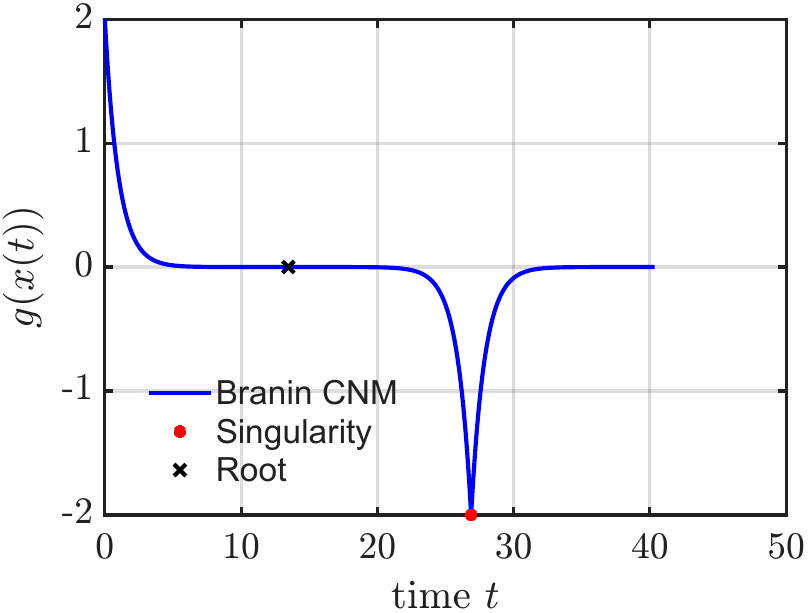}
    }
    %\hfill
    \subfigure[]{
    \includegraphics[width=0.3\textwidth]{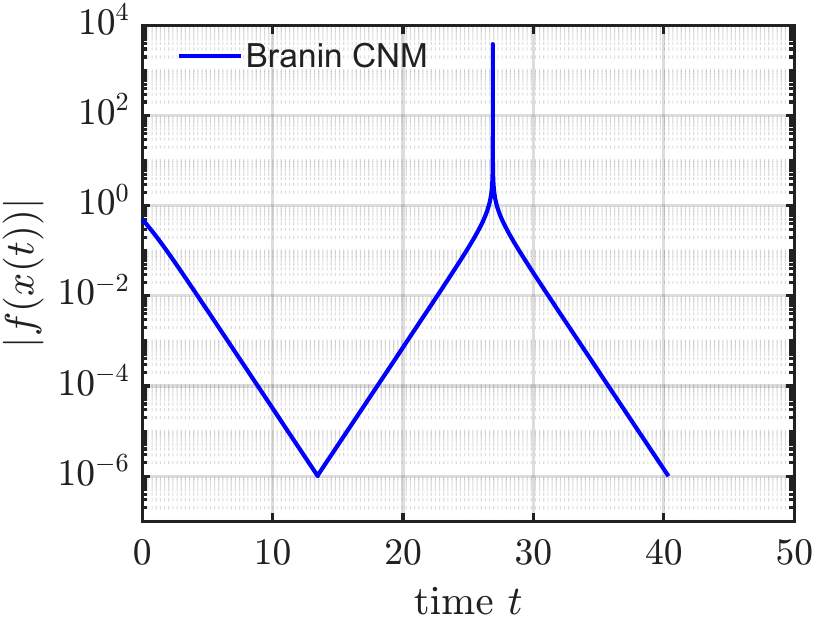}
    }
    \caption{\scriptsize Illustrative case study for $g(x) = x^2 - 2 = 0$ under the modified continuous Newton flow with initial condition $x_0 = -2$ and $\gamma(0) = -1$. (a) State trajectory $x(t)$; (b) evolution of $g(x(t))$; (c) semi-logarithmic plot of $|\dot x(t)|$, highlighting the blow-up of the vector field at the singularity crossing. With a black cross we highlight the point where the flow reaches the first root.\label{fig:trajx2}
    }
    \vspace{0.1cm}\hrule
\end{figure}

As a first didactic illustration of the proposed estimation framework, we consider the scalar equation
\begin{equation}
    g(x) = x^2 - 2 = 0,
    \label{eq:rootof2}
\end{equation}
whose solutions are $x^{*}_1 = -\sqrt{2}$ and $x^{*}_2 = \sqrt{2}$.
The associated Newton flow is given by
\begin{equation}
    \dot x(t) = \gamma(t) \frac{x(t)^2 - 2}{2x(t)}.
    \label{eq:rootof2flow}
\end{equation}
We start from $x(0) = -2$ with $\gamma(0) = -1$. The flow converges directly to the negative root $x^*_{1} = -\sqrt{2}$ (highlighted by a black cross in Fig.~\ref{fig:trajx2}). This root is reached with a tolerance $\varepsilon_f=10^{-6}$, around time $t\approx 13.47$. Having found one solution, we flip the orientation by setting $\gamma = 1$ and restart from a nearby point to the right of $x^*_{1}$, aiming to locate the remaining positive root. For ease of comparison with the analytical results below, we also run the simulation starting from $x(0) = -1$  and $\gamma(0)=+1$ (which corresponds to the post-root restart), though the figure shows the entire trajectory from $x(0) = -2$. The full trajectory passes through $x \approx -1$ at approximately $t \approx 26.19$.
%%%YGK  this corresponds to TIME exploding at finite x - that is what we use to make the two cases (x eplosion at finite t, and t exploding at finite x) into a circle in the paper with Panos - you don't have to say anything about it, or you could mention something - no pressure. 

With $\gamma = 1$, the flow evolves with positive velocity toward the singularity at $x = 0$, where the vector field blows up. At this instant, the right-hand side of \eqref{eq:rootof2flow} blows up. The integration crosses this infinity and resumes with $\gamma(t) = -1$ for $t > t^\ast$, allowing the trajectory to approach the positive root $x^*_{2} = \sqrt{2}$.

The full trajectories of $x(t)$ and $\dot x(t)$ are displayed in Fig.~\ref{fig:trajx2}. Near the infinity crossing, the adaptive framework described in Section~\ref{subsubsec:detection} estimates the self-similar parameters $(q, t^\ast, a)$ from \eqref{eq:selfsimilar}. Applying the estimation procedure to the threshold crossings, we obtain the following estimation of the self-similar parameters:
\begin{equation}
    \hat q = -0.500, \qquad \hat t^\ast = 0.693, \qquad \hat a = 0.707.
    \label{eq:estimated_params}
\end{equation}
These estimates can be verified analytically. Indeed, the exact solution of \eqref{eq:rootof2flow} for $\gamma = 1$ and $x(0)=-1$ is
\begin{equation}
    x(t) = -\sqrt{2 - e^{t}},
    \label{eq:exact_solution}
\end{equation}
from which the blow-up time is obtained by setting $x(t^\ast) = 0$, yielding $t^\ast = \log 2 \approx 0.693$. Furthermore, expanding the ODE near the singularity $x = 0$ gives the asymptotic approximation
\begin{equation}
    \dot x(t) = -\frac{1}{x(t)} + \mathcal{O}(x(t)),
    \label{eq:asymptotic_ode}
\end{equation}
whose solution near $t^\ast$ is
\begin{equation}
    x(t) \sim -\sqrt{2(t^\ast - t)},
    \qquad \text{and hence} \qquad 
    \dot x(t) \sim \frac{1}{\sqrt{2(t^\ast - t)}}.
    \label{eq:asymptotic_solution}
\end{equation}
Comparing with the self-similar ansatz in Eq.~\eqref{eq:selfsimilar}, we identify
\begin{equation}
    q = -\frac{1}{2}, \qquad a = \frac{1}{\sqrt{2}} \approx 0.707,
    \label{eq:analytical_params}
\end{equation}
which are in excellent agreement with the numerical estimates in Eq.~\eqref{eq:estimated_params}. This simple example thus demonstrates the accuracy of the proposed estimation framework and confirms that the self-similarity parameters can be recovered reliably from the computed trajectory data, even in the presence of a finite-time singularity.
%%%%%%%%%%%%%%%%%%%%%%%%%%%%%%

\subsubsection{Illustrative example: finite-time blow-up of the state}
\label{subsec:arctan_example}

We consider the scalar equation
\begin{equation}
    g(x) = \arctan(x) = 0,
    \label{eq:arctan_g}
\end{equation}
which has a unique real root at $x^\ast = 0$.
The corresponding continuous Newton flow is
\begin{equation}
    \dot x(t) = f(x(t)) = \gamma(t)\arctan(x(t)) \, (1 + x(t)^2).
    \label{eq:arctan_flow}
\end{equation}

Although the natural Newton/Davidenko choice $\gamma(0) = -1$ would lead directly to the root, we deliberately and illustratively set $\gamma(0) = 1$ to expose the method to a finite-time blow-up of the state.
%\color{black}
%
Starting from $x(0) = 1$, the flow initially moves toward larger positive values. However, since $g$ possesses a horizontal asymptote as $x \to +\infty$, we have $g'(x) \to 0$, causing the vector field to blow-up in finite time. Both $x(t)$ and $\dot x(t)$ diverge to $+\infty$; the trajectory crosses infinity via the projective compactification of $\mathbb{R}$ and reappears on the other side. After crossing, the orientation is updated to $\gamma(t) = -1$ and integration resumes, allowing the flow to converge to the root.

\begin{figure}[ht!]
    \centering
    \subfigure[]{
    \includegraphics[width=0.3\textwidth]{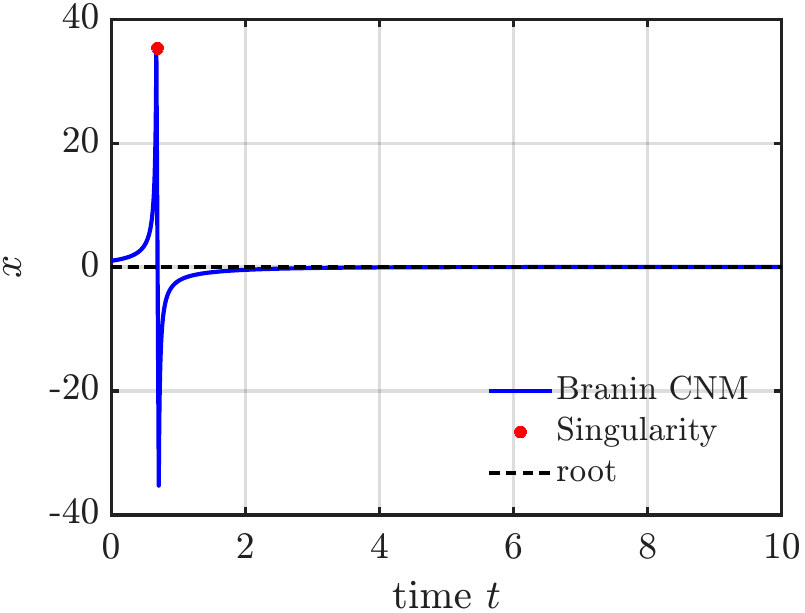}
    }
    %\hfill
    \subfigure[]{
    \includegraphics[width=0.3\textwidth]{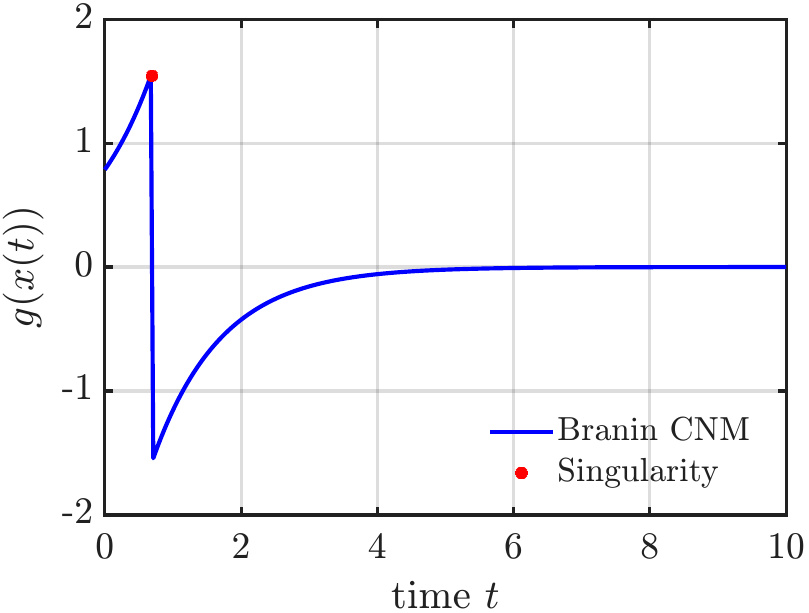}
    }
    %\hfill
    \subfigure[]{
    \includegraphics[width=0.3\textwidth]{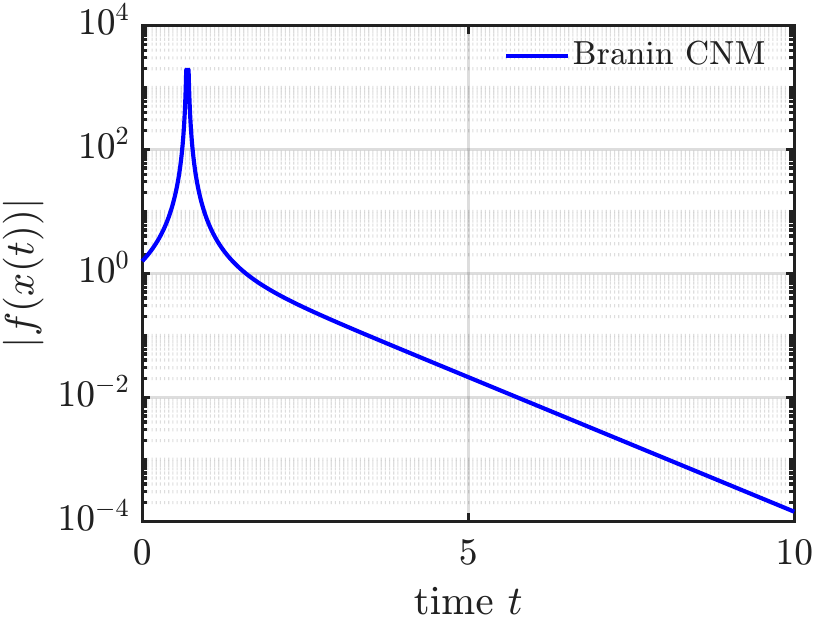}
    }
    \caption{\scriptsize Illustrative example for $g(x) = \arctan(x)$ with $x_0 = 1$ and $\gamma(0) = 1$. (a) State trajectory $x(t)$ showing finite-time blow-up; (b) evolution of $g(x(t))$; (c) semi-logarithmic plot of $|\dot x(t)|$.\label{fig:arctan}}
    \vspace{0.1cm}\hrule
\end{figure}

The trajectory and the blow-up of the vector field are shown in Fig.~\ref{fig:arctan}. Using the estimation framework of Section~\ref{subsubsec:detection}, we obtain the self-similar parameters
\begin{equation}
    \hat q = -1.997, \qquad \hat t^\ast = 0.693, \qquad \hat a = 0.642.
    \label{eq:arctan_estimates}
\end{equation}

These estimates can be verified analytically. For $\gamma = 1$, the ODE in Eq.~\eqref{eq:arctan_flow} admits the exact solution
\begin{equation}
    x(t) = \tan\left( \frac{\pi}{4} e^{t} \right),
    \label{eq:arctan_exact}
\end{equation}
which follows from the substitution $u = \arctan(x)$. Setting $x(t^\ast) = +\infty$ gives $\arctan(x(t^\ast)) = \pi/2$, so that
\begin{equation}
    t^\ast = \log 2 \approx 0.693,
    \label{eq:arctan_tstar}
\end{equation}
in agreement with the estimated value. For large $x$, we have $\arctan(x) \to \pi/2$ and $1 + x^2 \sim x^2$, so that the ODE asymptotically behaves as
\begin{equation}
    \dot x \sim \frac{\pi}{2} x^2.
    \label{eq:arctan_asymptotic_ode}
\end{equation}
Solving this asymptotic ODE gives
\begin{equation}
    x(t) \sim \frac{2}{\pi}(t^\ast - t)^{-1},
    \label{eq:arctan_x_asymptotic}
\end{equation}
and consequently, substituting back into the asymptotic expression for the vector field,
\begin{equation}
    f(x(t)) = \dot x(t) \sim \frac{2}{\pi}(t^\ast - t)^{-2}.
    \label{eq:arctan_f_asymptotic}
\end{equation}
Comparing with the self-similar ansatz $f(x(t)) = a (t^\ast - t)^q$ in Eq.~\eqref{eq:selfsimilar} yields
\begin{equation}
    q = -2, \qquad a = \frac{2}{\pi} \approx 0.637,
    \label{eq:arctan_analytical}
\end{equation}
which are in excellent agreement with the numerical estimates in Eq.~\eqref{eq:arctan_estimates}.

%\color{black}
\subsection{Infinity crossing for the infinite time grow-up scenario 
(\texttt{IT2})}
\label{subsec:infinity_crossing_in_infinite_time}

We now address the case in which the state $x(t)$ grows up in infinite time, that is condition \texttt{IT2}. This scenario arises when the vector field $f(x(t))$ remains of constant sign and non-negligible for all $t \geq t_1$, driving $x(t)$ toward $+\infty$ or $-\infty$ without any finite-time blow-up in $f$.
Rather than halting the integration, we preserve the orientation of the flow and allow the trajectory to pass through infinity, continuing on the other side of the compactified real line — in the spirit of the projective compactification discussed in~\cite{kevrekidis2017infinity}. This is consistent with the orientation-preserving philosophy outlined in Section~\ref{sec:intro}. %: reversing the flow direction upon detecting large values of $|x(t)|$ would risk undoing progress already made toward a potential solution.

\subsubsection{Detection}
\label{subsubsec:detection_type2}
The approach of $x(t)$ toward infinity is detected via the event function
\begin{equation}
    E_{\textrm{state}}(x) = |x| - r_{\max},
    \label{eq:event_type2}
\end{equation}
where $r_{\max} \gg 1$ is a user-defined threshold considered a practical proxy for infinity.
The use of the absolute value captures both $x \to +\infty$ and $x \to -\infty$ within a 
single event, with the sign of $x$ at the trigger time determining the direction of the divergence.
Detection is triggered when $E_{\textrm{state}}(x)$ transitions from negative to positive, i.e., when $|x(t)|$ first exceeds $r_{\max}$, while simultaneously $|f(x(t))|$ remains below $f_{\max}$, distinguishing this scenario from the case of \texttt{IT1}.

\subsubsection{Time reparameterization}
\label{subsubsec:time_reparam}
The fundamental difficulty of \texttt{IT2} is that the crossing must occur at $t = +\infty$, which is not directly accessible to a finite-time ODE integrator.
To overcome this, we apply the time reparameterization~\cite{wuerth2009dynamic}
\begin{equation}
    t(s) = t_1 + \frac{s}{1-s},
    \label{eq:time_transform}
\end{equation}
where $t_1$ is the time at which $|x(t_1)| = r_{\max}$. This transformation satisfies
\begin{equation}
    t(0) = t_1 \qquad \text{and} \qquad 
    \lim_{s \uparrow 1} t(s) = +\infty,
\end{equation}
mapping the semi-infinite interval $[t_1, +\infty)$ bijectively onto the finite interval $[0, 1)$. Under this reparameterization, the ODE in Eq.~\eqref{eq:continuous_newton} transforms to
\begin{equation}
    \frac{\dd x}{\dd s} = f(x(t(s))) \cdot \frac{\dd t}{\dd s} = 
    \frac{f(x(t(s)))}{(1-s)^2},
    \label{eq:ode_transformed}
\end{equation}
where we used $\dd t/\dd s = (1-s)^{-2}$. The factor $(1-s)^{-2}$ amplifies the vector field as $s \to 1$, effectively accelerating the trajectory toward infinity in the finite pseudo-time $s$.

To verify that this transformation reduces the problem to a finite-time blow-up, suppose that for some $\hat{t}_1 \geq t_1$ the vector field satisfies
\begin{equation}
    f(x(t)) \geq \frac{c}{(t - t_1)^r}, \qquad t > \hat{t}_1,
    \label{eq:f_lower_bound}
\end{equation}
for some $c > 0$ and $0 < r \leq 1$. This condition is mild and is satisfied, for instance, whenever $f$ is bounded away from zero or decays no faster than $1/t$.
Under the transformation in Eq.~\eqref{eq:time_transform}, Eq.~\eqref{eq:f_lower_bound} implies
\begin{equation}
    \frac{\dd x}{\dd s} \geq c s^{-r} (1-s)^{r-2}, \qquad s>\hat s_1,
\end{equation}
where $\hat{s}_1$ is defined by $t(\hat{s}_1) = \hat{t}_1$. 
Since $r - 2 \leq -1$ for $0 < r \leq 1$, the right-hand side is not integrable on $[\hat{s}_1, 1)$, confirming that both $x(s)$ and $\dd x/\dd s$ diverge in finite pseudo-time as $s \to 1$.
The problem has thus been reduced to a finite-time blow-up in the $s$-variable, to which the machinery of Subsection~\ref{subsec:infinity1} for \texttt{IT1} case applies directly.

\subsubsection{Illustrative example: state growing-up in infinite time}
\label{subsec:infinity_infinite_time}
We consider a case where the state variable grows-up in infinite time due to the presence of a horizontal asymptote of $g$, corresponding to a root at infinity. Let
\begin{equation}
    g(x) = e^{-x} \bigl((x^2 - 4)x + 4\bigr),
    \label{eq:rootofscalar}
\end{equation}
which has a unique real root $x^\ast \approx -2.3830$, corresponding to the single real root of the cubic polynomial $p(x) = x^3 - 4x + 4$.

\begin{figure}[ht!]
\centering
    \subfigure[]{
    \includegraphics[width=0.4\textwidth]{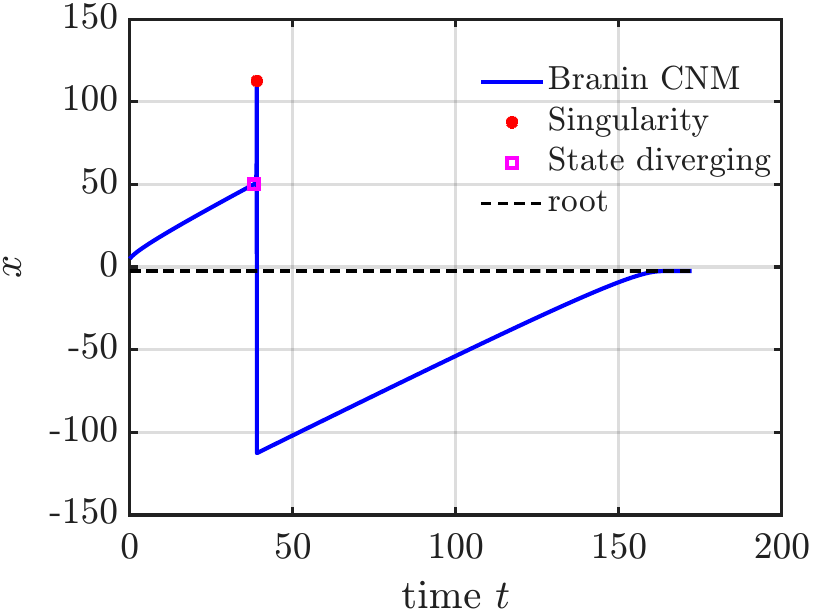}
    }
    \subfigure[]{\includegraphics[width=0.4\textwidth]{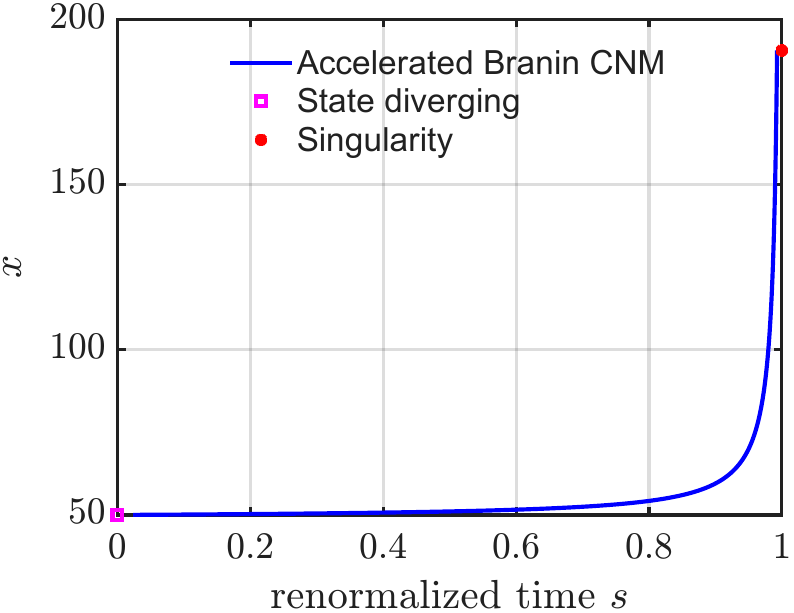}
    }\\
    %\hfill
    \subfigure[]{
    \includegraphics[width=0.4\textwidth]{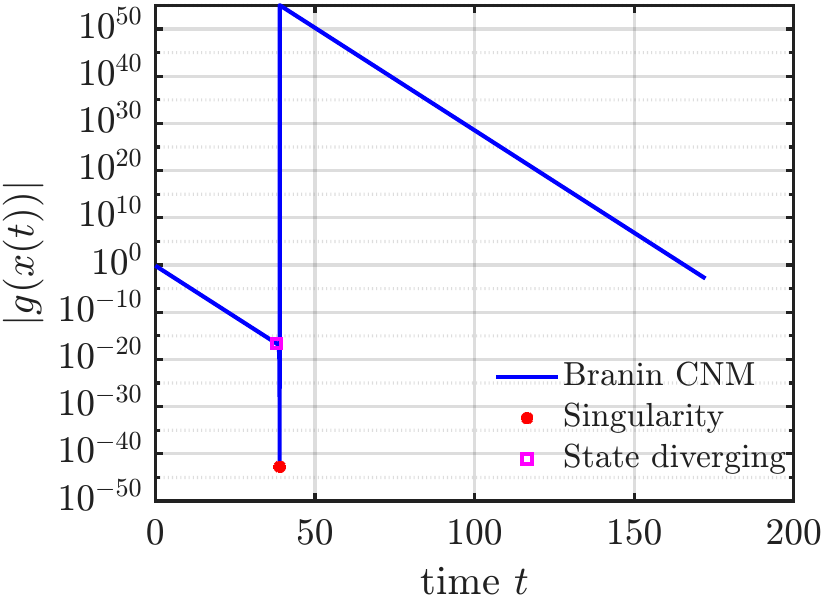}
    }
    %\hfill
    \subfigure[]{
    \includegraphics[width=0.4\textwidth]{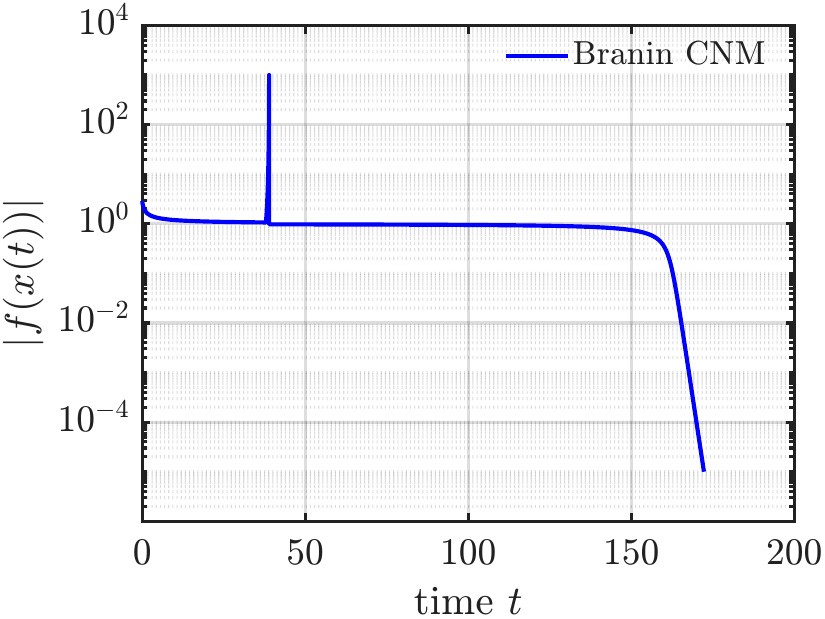}
    }
    \caption{\scriptsize Illustrative example for $g(x) = e^{-x}(x^3 - 4x + 4)$ with $x_0 = 5$ and $\gamma(0) = -1$. (a) State trajectory $x(t)$; (b) zoom-in of the accelerated dynamics in renormalized time $s\in[0,1)$; (c) semi-logarithmic plot of $|G(x(t))|$; (d) semi-logarithmic plot of $|\dot x(t)|$. Note that as $x \to +\infty$, $|G(x(t))|$ becomes very small, but $|\dot x(t)|$ approaches a finite nonzero value (approximately 1), indicating a horizontal asymptote rather than a root.\label{fig:case2b}}
\vspace{0.1cm}\hrule
\end{figure}

The continuous Newton flow is simulated with initial condition $x(0) = 5$ and orientation $\gamma(0) = -1$. In this configuration, the trajectory moves toward larger positive values, approaching the horizontal asymptote of $g$ as $x \to +\infty$. Since the asymptote is approached only in infinite time, the vector field does not blow-up; instead, $|f(x(t))|$ tends to a finite nonzero value. To handle this case numerically, we declare $|x(t)| > r_{\max} = 50$ (reached around $t=37.95$) as a practical proxy for infinity and apply the time transformation in Eq.~\eqref{eq:time_transform} to accelerate the dynamics. This transformation introduces the visible spikes in the state and right-hand side at $t^* \approx 38.947$ and $q\approx-1.988$ (see Fig.~\ref{fig:case2b}).
{A zoom-in of the accelerated dynamics is showed in panel (b).}
After crossing infinity via the projective compactification of $\mathbb{R}$, the sign of $x(t)$ is flipped, the orientation is updated, {and integration resumes in the original time variable $t$}, allowing the trajectory to converge to the unique root $x^\ast \approx -2.3830$.

The trajectories are shown in Fig.~\ref{fig:case2b}. In panel (c), $|G(x(t))|$ decays exponentially as $x \to +\infty$, which might suggest the presence of a root at infinity. However, panel (d) reveals that $|\dot{x}(t)|$ tends to a finite positive value rather than zero, confirming that this is indeed a horizontal asymptote and not a finite root in the usual sense. In a certain compactified sense, one could view $x = +\infty$ as a root at infinity, since $\lim_{x\to\infty} g(x) = 0$. However, this is not a zero crossing: the function remains strictly positive for all finite $x>x^\ast$, and there is no point where $g(x)$ changes sign. 

\section{Infinity crossings for multidimensional problems}
\label{sec:multidimensionalodes}
The extension of the infinity-crossing procedure to $n$-dimensional systems follows the same structure as the scalar case developed in Section~\ref{sec:scalarodes}. The key difference is that the blow-up may now involve multiple components of the vector field $\bm{f}(t, \bm{x}(t))$ simultaneously, requiring a unified detection strategy that avoids the computational complexity of monitoring each component independently.
%
%%%%%%%%%%%%%%%%%%%%%%%%%%%%%%%%%%%
\subsection{Singularity structure of the continuous Newton flow}
\label{subsec:singularity_structure_nd}

To analyze the singularities of the CNM flow in higher dimensions, we recall the vector field definition. With orientation factor $\gamma(t)\in\{-1,+1\}$, the flow is
\begin{equation}
\dot{\bm{x}}(t)
:= \bm{f}(\bm{x}(t),t)
= \gamma(t)\, \mathbf{J}(\bm{x})^{-1}\bm{g}(\bm{x}),
\qquad
\mathbf{J}(\bm{x}) := \nabla \bm{g}(\bm{x}).
\label{eq:nd_newton_flow}
\end{equation}
This is well-defined on $\mathbb{R}^n \setminus \mathcal{S}$, where $\mathcal{S}$ is the set of singularities, i.e.,
\begin{equation}
\mathcal{S} := \{\bm{x}\in\mathbb{R}^n : \det(\mathbf{J}(\bm{x})) = 0\}.
\label{eq:singular_set}
\end{equation}

Using the adjugate identity, the vector field can be written as
\begin{equation}
\bm{f}(\bm{x})
=
\gamma(t)\,
\frac{\adj(\mathbf{J}(\bm{x}))\,\bm{g}(\bm{x})}
{\det(\mathbf{J}(\bm{x}))}.
\label{eq:nd_adjugate_form}
\end{equation}
%
%\begin{remark}[$\adj(\mathbf{J})\bm g=\bm 0$ requires $\det(\mathbf{J})=0$

%
Following Branin's terminology~\cite{branin1972widely}, a singular point $\bm{x}^\ast \in \mathcal{S}$ is defined as a \emph{extraneous singularity} when both the Jacobian is singular and the adjugate action on the residual vanishes, i.e.
\begin{equation}
\adj(\mathbf{J}(\bm{x}^\ast))\,\bm{g}(\bm{x}^\ast) = \bm{0}.
\label{eq:branin_extraneous_def}
\end{equation}
%so that the leading-order pole in Eq.~\eqref{eq:nd_adjugate_form} is canceled and the vector field does not exhibit a first-order blow-up despite $\det(\mathbf{J}(\bm{x}^\ast)) = 0$.
%In this case, the local behavior of the CNM flow is determined by higher-order terms in the expansions of $\mathbf{J}$ and $\bm{g}$.
%
Note that, if $\det(\mathbf{J}(\bm x))\neq0$, then $\mathbf{J}(\bm x)^{-1}$ exists, is invertible and $\adj(\mathbf{J}(\bm x))=\det(\mathbf{J}(\bm x))\,\mathbf{J}(\bm x)^{-1}$, so that
\begin{equation}
\adj(\mathbf{J}(\bm x))\bm g(\bm x)=\bm 0
\;\Longleftrightarrow\;
\mathbf{J}(\bm x)^{-1}\bm g(\bm x)=\bm 0
\;\Longleftrightarrow\;
\bm g(\bm x)=\bm 0.
\end{equation}
Hence, away from $\mathcal S$, the condition $\adj(\mathbf{J})\bm g=\bm 0$ only occurs at a regular root of $\bm g$, not at a singularity. Consequently, the Branin \emph{extraneous-singularity} condition~\eqref{eq:branin_extraneous_def} is nontrivial only on $\mathcal S$. %=\{\det(\mathbf{J})=0\}$.

We now express the Jacobian using the singular value decomposition
\begin{equation}
\mathbf{J}(\bm{x})
=
\mathbf{U}(\bm{x})\,\mathbf{\Sigma}(\bm{x})\,\mathbf{V}(\bm{x})^\top,
\qquad
\mathbf{\Sigma}(\bm{x})=\mathrm{diag}(\sigma_1,\dots,\sigma_n),
\label{eq:nd_svd}
\end{equation}
with $\sigma_1 \geq \cdots \geq \sigma_n \geq 0$.
Substituting these expressions into the vector field yields
\begin{equation}
\bm{f}(\bm{x})
=
\gamma(t)\,
\sum_{i=1}^n
\frac{\bm{u}_i(\bm{x})^\top \bm{g}(\bm{x})}{\sigma_i(\bm{x})}
\,\bm{v}_i(\bm{x}),
\label{eq:nd_flow_svd}
\end{equation}
where $\bm{u}_i(\bm{x})$ and $\bm{v}_i(\bm{x})$ denote the $i$-th left and right singular vectors.

Using the SVD representation, we can also rewrite the determinant and the adjugate. The determinant reads
\begin{equation}
\det(\mathbf{J}(\bm{x})) = \det(\mathbf{U}(\bm{x}))\det(\mathbf{V}(\bm{x}))\prod_{i=1}^n \sigma_i(\bm{x}),
\label{eq:det_svd}
\end{equation}
while the adjugate admits the expansion
\begin{equation}
\begin{split}
\adj(\mathbf{J}(\bm{x}))
&=\det(\mathbf{U}(\bm{x}))\det(\mathbf{V}(\bm{x}))\, \mathbf{V}(\bm{x})\,\adj(\mathbf{\Sigma}(\bm{x}))\,\mathbf{U}(\bm{x})^\top \\
&=
\pm \sum_{i=1}^n
\left(
\prod_{j\neq i}\sigma_j(\bm{x})
\right)
\bm{v}_i(\bm{x})\,\bm{u}_i(\bm{x})^\top.
\end{split}
\label{eq:adj_svd}
\end{equation}

If $\operatorname{rank}(\mathbf{J}(\bm{x}^\ast)) \le n-2$, then at least two singular values vanish at $\bm{x}^\ast$, implying $\adj(\mathbf{J}(\bm{x}^\ast))=\bm{0}$. Consequently, the singular point is automatically a Branin extraneous singularity (see Eq.~\eqref{eq:branin_extraneous_def}).
%%%YGK  what is an EXTRANEOUS singularity ? and, one line below as above, what is a SIMPLE singularity ? 
%
Whereas, if $\operatorname{rank}(\mathbf{J}(\bm{x}^\ast))=n-1$ and the singularity has a single vanishing singular value (a corank-one, single-direction singularity), i.e. $\sigma_n(\bm{x}^\ast)=0$ and $\sigma_{n-1}(\bm{x}^\ast)>0$---we call such a singularity \emph{simple}---then
\begin{equation}
    \det(\mathbf{J}(\bm{x}))
    \sim
    C\,\sigma_n(\bm{x}),
    \qquad
    C=\pm\prod_{i=1}^{n-1}\sigma_i(\bm{x}^\ast)\neq 0,
    \label{eq:det_svd_asymptote}
\end{equation}
locally along the trajectory as $\bm{x}\to\bm{x}^\ast$, with $C$ evaluated at $\bm{x}^\ast$.
Additionally, the adjugate admits the asymptotic expansion
\begin{equation}
    \adj(\mathbf{J}(\bm{x}))
    =
    C \,\bm{v}_n(\bm{x}^\ast) \bm{u}_n(\bm{x}^\ast)^\top
    +O(\sigma_n(\bm{x})),
    \label{eq:adj_simplesing}
\end{equation}
where $C$ is as in Eq.~\eqref{eq:det_svd_asymptote}.
Hence, near a \emph{simple} singularity,
\begin{equation}
    \bm{f}
    =
    \gamma(t)\frac{\adj(\mathbf{J})\,\bm{g}}
    {\det(\mathbf{J})}
    =
    \gamma(t)\frac{\bm{u}_n^\top\bm{g}}
    {\sigma_n}\,
    \bm{v}_n
    +O(1),
\end{equation}
showing that the local dynamics are completely determined by two scalar quantities: the smallest singular value $\sigma_n$ and the projected residual
\begin{equation}
    \phi(\bm{x})
    :=
    \bm{u}_n(\bm{x})^\top\bm{g}(\bm{x}).
\end{equation}
These play the role of the scalar derivative $g'$ and the scalar residual $g$, respectively, in the one-dimensional case.

Moreover, for \emph{non-extraneous} singularity (see Eq.~\eqref{eq:branin_extraneous_def}), 
%%%YGK same question - what is non-extraneous ?? 
the direction of the flow $\bm v_f$ is asymptotically aligned with the right singular vector $\bm{v}_n(\bm{x}^\ast)$,
\begin{equation}
    \bm v_{f}:=\frac{\bm{f}(\bm{x})}{\|\bm{f}(\bm{x})\|_2} = \gamma(t)\sign(\phi(\bm{x}^\ast))\,\bm{v}_n(\bm{x}^\ast) + O(\sigma_n(\bm{x})),
    \label{eq:v_f}
\end{equation}

\noindent demonstrating that the CNM flow escapes to infinity along the one-dimensional subspace spanned by $\bm{v}_n(\bm{x}^\ast)$ at a rate proportional to $\sigma_n^{-1}$.

\begin{remark}
We emphasize that the method does not assume the trajectory is necessarily attracted to a singularity or that it will pass through one. In the analysis above, we considered--as a working assumption--the case where the trajectory does approach a singularity at a certain point, which is one of the possible scenarios that may arise during integration. In the CNM flow, singularities generically form hypersurfaces defined by $\det(J(x)) = 0$. Such a hypersurface may separate the high-dimensional space into different regions, and a trajectory may be forced to cross it to move from one region to another--as occurs, for instance, when crossing separatrices between solution branches.
In this sense, the 1D intuition of being ``forced through" a singularity extends naturally to higher dimensions when singularities are hypersurfaces.
However, isolated singularities (points) in higher dimensions are more unlikely to be encountered. Indeed, in higher dimensions a trajectory can bypass an isolated point by moving in a transverse direction. For such a point to be encountered, the singularity must possess some notion of local attractivity in the flow. The method does not attempt to force encounters with isolated singularities; it merely detects and crosses them when they occur, whether they are points on hypersurfaces or isolated points.
\end{remark}

\subsection{Directional self-similar blow-up detection}
\label{subsec:directional_detection}
The central observation is that, near a blow-up, the vector field $\bm{f}(t, \bm{x}(t))$ grows large along a dominant direction. 
We have defined in Eq.~\eqref{eq:v_f} the instantaneous flow direction $\bm{v}_f(t)$. Note that the scalar projection of $\bm{f}$ onto $\bm{v}_f$ satisfies
\begin{equation}
    \bm{v}_f^\top \bm{f}(t, \bm{x}(t)) = 
    \|\bm{f}(t, \bm{x}(t))\|_2.
    \label{eq:projection}
\end{equation}
This reduces the $n$-dimensional detection problem to a single scalar quantity: the norm of the vector field. Detection of an upcoming blow-up is therefore implemented via the single event function
\begin{equation}
    E_{\textrm{fun}}(\bm{x}) = \|\bm{f}(t, \bm{x})\|_2 - f_{\max},
    \label{eq:event_nd}
\end{equation}
which triggers when $\|\bm{f}\|_2$ first exceeds $f_{\max}$, regardless of which components are responsible for the growth. 

\begin{assumption}
\label{ass:vf_constant}
In a neighborhood of the blow-up time $t^\ast$, the blow-up direction $\bm{v}_f(t)$ remains approximately constant, i.e., $\bm{v}_f(t) \approx \bm{v}_f^\ast$ for $t \in (t^\ast - \varepsilon, t^\ast)$ and some $\varepsilon > 0$.
\end{assumption}

Under Assumption~\ref{ass:vf_constant}, the scalar projection 
$\|\bm{f}(t, \bm{x}(t))\|_2$ inherits the self-similar structure 
of the scalar case, and the directional self-similar ansatz reads
\begin{equation}
    \|\bm{f}(t, \bm{x}(t))\|_2 = a\,|t^\ast - t|^q,
    \label{eq:selfsimilar_nd}
\end{equation}
with $a > 0$ and $q < 0$. The parameters $(q, t^\ast, a)$ are 
estimated from the sequence of threshold crossings of 
$E_{\textrm{fun}}(\bm{x})$ exactly as in Section~\ref{subsubsec:detection}, 
with $|f(x(t))|$ replaced by $\|\bm{f}(t,\bm{x}(t))\|_2$ 
throughout. In particular, the linear regression 
formula in Eq.~\eqref{eq:q_estimate} for $q$, the blow-up time estimate in Eq.~\eqref{eq:tstar_estimate}, and the amplitude estimate in Eq.~\eqref{eq:a_estimate} all carry over unchanged.

%\color{gray!30!red}
\begin{remark}
Assumption~\ref{ass:vf_constant} excludes rotating or spiraling blow-up, in which the direction $\bm{v}_f(t)$ rotates while $\|\bm{f}\|_2$ diverges. For general ODE systems, such behavior cannot be ruled out {\em a priori}. However, for the CNM flow prior to a sign switch ($\gamma \equiv -1$), the dynamics are gradient-like with respect to the Lyapunov function $L(\bm{x}) = \|\bm{g}(\bm{x})\|_2^2 = \bm{g}(\bm{x})^\top \bm{g}(\bm{x})$. Indeed, along trajectories of the CNM flow,
\begin{equation}
\begin{split}
    \dot{L}(\bm{x}(t)) = \frac{\dd}{\dd t} \left(\bm{g}^\top \bm{g}\right)
    = 2 \bm{g}^\top \dot{\bm{g}}
    &= 2 \bm{g}^\top \left(\nabla \bm{g} \, \dot{\bm{x}}\right)
    = 2 \bm{g}^\top \left(\nabla \bm{g} \, (-\nabla \bm{g}^{-1} \bm{g})\right)
    =\\
    &= -2 \|\bm{g}\|_2^2
    = -2 L(\bm{x}(t)) \leq 0,
\end{split}
\end{equation}
where we used the chain rule, the definition of the CNM flow $\dot{\bm{x}} = -\nabla \bm{g}^{-1} \bm{g}$, and the fact that $\nabla \bm{g} \nabla \bm{g}^{-1} = \mathbf{I}_n$. Thus $L$ decreases monotonically along trajectories for $\gamma \equiv -1$, and no periodic or chaotic solutions can exist~\cite{zufiria2002role}.
Note that the Lyapunov function $L$ considered above is also the basis of the squared-operator method (SOM)~\cite{subramanian2021exploring}, which reformulates root-finding as the minimization of $\|g(x)\|^2$.
Moreover, as shown above, for non-extraneous singularities, the vector field asymptotically aligns with the right singular vector corresponding to the smallest singular value of the Jacobian, which stabilizes the blow-up direction. Thus, Assumption~\ref{ass:vf_constant} is justified for the class of problems considered in this work.
\end{remark}

%\color{black}
%%%%%%%%%%%%%%%%%%%%%%%%%%%%%%%%%%%%%
\subsection{Crossing strategy}
\label{subsec:crossing_nd}
With $(q, t^\ast, a)$ estimated along the blow-up direction 
$\bm{v}_f^\ast$, the crossing strategy mirrors that of 
Section~\ref{subsubsec:crossing}, applied directionally.

\textit{Case $-1 < q < 0$: vector field blows-up, 
state bounded.}
The state $\bm{x}(t)$ remains bounded across the singularity 
and the crossing correction is applied along the blow-up 
direction:
\begin{equation}
    \bm{x}(t^\ast + h) = \bm{x}(t^\ast - h) + 
    \bm{v}_f^\ast \cdot \frac{2ah^{q+1}}{q+1}.
    \label{eq:crossing_nd}
\end{equation}
This is a single vector correction along $\bm{v}_f^\ast$, 
consistent with the geometry of the flow near the singularity.

\textit{Case $q \leq -1$: both vector field and 
state blow-up.}
The state $\bm{x}(t)$ itself diverges in finite time along 
the blow-up direction. As in the scalar case, motivated by 
the projective compactification discussed in~\cite{kevrekidis2017infinity}, we reset the state via
\begin{equation}
    \bm{x}(t^\ast + h) = -\bm{x}(t^\ast - h).
    \label{eq:sign_flip_nd}
\end{equation}
The underlying idea is to compactify $\mathbb{R}^n$ by identifying antipodal points on the sphere at infinity,  so that a trajectory that diverges to infinity along a given ray re‑emerges from the opposite direction. The natural reset is therefore $\bm{x} \leftarrow -\bm{x}$, which maps the state to the antipodal point on the sphere.

In practice, the blow‑up may not be perfectly radial; the state may be large in only a few components while others remain small. In such cases, one may choose to flip only those components that contribute most significantly to $\|\bm{x}\|_2$, e.g., those for which $|x_i|$ exceeds a fraction of $\|\bm{x}\|_2$. This component‑wise reset reflects the observation that the singularity is effectively one‑dimensional along the dominant direction, and the transverse components are not involved in the blow‑up. 
%The global flip $\bm{x} \leftarrow -\bm{x}$ is a simpler alternative that treats all components equally. Both heuristics are defensible; they may lead to slightly different subsequent trajectories, but in practice the flow typically realigns quickly after crossing, and the choice does not affect the ultimate convergence to a root. In this work we adopt the global flip for its simplicity and conceptual clarity.

\subsubsection{Orientation update}
%\color{gray!30!red}
Consistent with Branin's formulation, we consider only the case where the orientation factor takes the values $\gamma(t) = \pm 1$. %, corresponding to a global reversal of the vector field.

% case in which we can change the entire sign of the orientation $\gamma$ from $-1$ to $+1$ or viceversa. Another option (as done by Ralf) is to change only a few particular directions of the flow, thus only some elements of the diagonal entries of $\gamma$. However this creates a mixed behavior, for which there does not exist a clean solution as in the Branin solution in Eq.~\eqref{eq:branin_solution}; instead one obtains
% \begin{equation}
% \dot{\bm{g}}=\nabla \bm{g} \, \gamma \, \big(\nabla \bm{g}\big)^{-1} \bm{g}.
% \end{equation}
% which is not equivalent to $\gamma \bm{g}$ unless $\nabla \bm{g}$ and $\gamma$ commute. Consequently, there is no straightforward interpretation of $\gamma$ in terms of the Morse index of the $\bm{g}$; the orientation update is instead dictated by the sign change of the projection onto the blow‑up direction, as derived in the scalar case and extended directionally.
%
The sign of the orientation factor $\gamma$ is updated according to Eq.~\eqref{eq:sigma_update}, with the sign of the 1d ODE $f$, replaced by the sign of the projection $\bm{v}_f^{*\top} \bm{f}$:
\begin{equation}
    \gamma(t^\ast + h) = \gamma(t^\ast - h) \cdot 
    \sign\!\left(\bm{v}_f^{*\top} 
    \bm{f}(\bm{x}(t^\ast - h))\right) \cdot 
    \sign\!\left(\bm{v}_f^{*\top} 
    \bm{f}(\bm{x}(t^\ast + h))\right),
    \label{eq:sigma_update_nd}
\end{equation}
where the direction of the flow $\bm{v}_f^{*\top}$ is computed at time $t^\ast-h$.
%
%Our choice of a global sign flip for $\gamma$ is deliberately pragmatic: it aims to continue the flow along the current direction of motion, projected onto the blow‑up direction.
%
The scalar projection $\phi = \mathbf{v}_f^\top \mathbf{f}$ reduces the dynamics to an effective one‑dimensional problem along the direction $\mathbf{v}_f$. Along this direction, the possible local behaviors are exactly those of the scalar case: local minima, local maxima, and inflection points of the projected dynamics.
%
% By flipping all entries of $\gamma$, we preserve the alignment of the flow with $\mathbf{v}_f$ after crossing, which is consistent with the one‑dimensional reduction.
% %
% This approach does not attempt to capture higher‑order saddle points or degeneracies associated with the full Hessian of $\bm{g}$; such features are not resolved by the scalar projection and would require a more elaborate treatment beyond the scope of this work.
%\color{black}

\subsection{Infinity crossing in infinite time}
\label{subsec:infinity2_nd}
For \texttt{IT2} case in $n$ dimensions, the state $\bm{x}(t)$ itself grows-up to infinity in infinite time. Detection is implemented via the event function
\begin{equation}
    E_{\textrm{state}}(\bm{x}) = \|\bm{x}\|_2 - r_{\max},
    \label{eq:event_state}
\end{equation}
which triggers when $\|\bm{x}\|_2$ first exceeds a prescribed threshold $r_{\max} \gg 1$. This threshold defines a large sphere $\{\bm{x} : \|\bm{x}\|_2 = r_{\max}\}$ as a practical proxy for infinity.

When the grow‑up occurs in infinite time, the vector field does not diverge, but the state grows without bound as $t \to \infty$. In this situation, the time reparameterization introduced in Section~\ref{subsubsec:time_reparam} applies without modification, with $\bm{f}(t, \bm{x}(t))$ now vector-valued. The transformed ODE reads
\begin{equation}
    \frac{\dd\bm{x}}{\dd s} = \frac{\bm{f}(\bm{x}(t(s)))}{(1-s)^2},
    \label{eq:ode_transformed_nd}
\end{equation}
where $t(s) = t_1 + s/(1-s)$ maps the semi‑infinite interval $[t_1, \infty)$ onto $[0,1)$. This transformation reduces the infinite‑time grow‑up to a finite‑time blow‑up in the pseudo‑time $s$, allowing the detection and crossing machinery of Section~\ref{subsec:crossing_nd} to be applied directly, with all scalar quantities replaced by their $n$-dimensional counterparts and detection performed via the norm‑based event in Eq.~\eqref{eq:event_nd}.

%%%%%%%%%%%%%%%%%%%
\subsection{Passing a root and continuing integration}
\label{subsec:passing_root}

When the integration converges to a root $\bm{x}^\ast$ (detected when $\|\bm{g}(\bm{x})\|$ or $\|\bm{f}\|$ falls below a prescribed tolerance), the algorithm can continue the search for additional solutions by passing beyond the root. Following Branin's strategy, we flip the orientation factor to $\gamma = +1$ and restart from a point just beyond the root along the last direction of approach.
This step is natural once one recognizes that reaching a steady state is, in a precise sense, the mirror image of finite-time blow-up: the former corresponds to \emph{time blow-up at finite state} ($x\to x^\ast$ as $t\to\infty$), the latter to state blow-up at finite time ($x\to\infty$ as $t\to t^*$). In the projective compactification of \cite{kevrekidis2017infinity} these are the two ends of a single circle (cf.\ Eqs.~(5)--(6) therein), and crossing the root in this sense is the dual of the infinity crossing.
%%%YGK  this corresponds again to "time going to infinity for finite x
%% as discussed in the paper with Panos - we can then do "dt/dx instead %% of dx/dt"

Let $\bm{d} = \bm{x}_k - \bm{x}_{k-1}$ be the last direction of integration. We generate a set of candidate points along this direction:
\begin{equation}
    \bm{x}_{\text{cand}}^{(j)} = \bm{x}_k + j \, \bm{d}, \qquad j = 1, \ldots, M,
    \label{eq:candidates}
\end{equation}
where $M$ is a fixed number of trial steps (e.g., $M = 100$). We then evaluate $\|\bm{g}(\bm{x}_{\text{cand}}^{(j)})\|$ and select the first candidate for which the residual exceeds both the residual at the root and a prescribed tolerance:
\begin{equation}
    \|\bm{g}(\bm{x}_{\text{cand}}^{(j)})\| > \|\bm{g}(\bm{x}_k)\| 
    \quad \text{and} \quad 
    \|\bm{g}(\bm{x}_{\text{cand}}^{(j)})\| > \varepsilon_g.
    \label{eq:candidate_selection}
\end{equation}
%
%\color{blue!90!black}
The tolerance $\varepsilon_g$ is chosen empirically; in our numerical experiments, values in the range $10^{-6}$ to $10^{-4}$ (relative to the scale of $\bm{g}$) were found to work reliably.
%\color{black}
The orientation is then set to $\gamma = +1$ (necessarily, since convergence to roots occurs only with $\gamma = -1$), and integration resumes. This heuristic allows the trajectory to escape the basin of attraction of the current root and explore other regions of the solution landscape.

%\color{blue!90!black}
We emphasize that the algorithm does not attempt to guarantee the discovery of all roots. In our implemented code, the number of roots sought is a user-specified parameter, as is the maximum integration time. In practice, the procedure may find fewer roots than requested, may encounter the same root multiple times (if the trajectory revisits it after going to infinity), or may find no further roots after a given point—for instance, when the flow diverges monotonically after setting $\gamma = +1$ and never switches back. Thus, the method is best understood as an exploration tool that can locate multiple roots along a trajectory, rather than a complete root-counting procedure. Conditions under which all roots lie on a single connected trajectory of the Branin Newton flow are known only under restrictive assumptions~\cite{diener1987global, zufiria2002role}.
%%%%%%%%%%%%%%%%

\section{Algorithm: continuous Newton's method with infinity crossing}
\label{sec:algorithm}
%%%%%%%%%%%%%%%%%%%%%
\begin{algorithm}[t!]
\caption{Continuous Newton with Blow-up detections}
\label{alg:contNewton}
\begin{algorithmic}[1]
\STATE \textbf{Input:} $@(\bm x) \bm g(\bm x), @(\bm x)\mathbf{J}(\bm x),\bm x_0,t_{\rm end}$
\STATE \textbf{Detection params:} $f_{\max},\,r_{\max},\,c_{\rm grow},\,N_{\rm detc},\,\varepsilon_g,\,\varepsilon_f$ ; \;
\STATE Define flow: $\dot{\bm{x}}=@(t,\bm x) \bm f(t,\bm x)=\gamma\,\mathbf{J}(\bm x)^{-1}\bm g(\bm x)$; \;
\STATE Define Events: \(E_1=\|\bm g(\bm x)\|-\varepsilon_g;\; E_2=\|\bm f\|-\varepsilon_f,\); \;
\STATE \(E_{3}=\|\bm f\|-f_{\max};\; E_{4}=\|\bm x\|-r_{\max}.\);\; %$(i=1:\text{dimension}(x))$
\STATE $t_{\rm curr}\leftarrow0,\;\bm x_{\rm curr}\leftarrow \bm x_0$.
\WHILE{$t_{\rm curr}<t_{\rm end}$}
    \STATE $[t_{all},x_{all}]\leftarrow$ Integrate $ \dot{\bm{x}}=\bm f(t,\bm x)$ %until $t_{end}$
    , with Events $\{E_1,E_2,E_3,E_4\}$;
    \IF{$E_1\le0$ and $E_2\le0$ events triggered}
        \RETURN trajectories $[t_{all},x_{all}]$.
    \ENDIF
    \IF{$E_{3}\ge0$ triggered}
        \STATE $[q,t^{\ast},a]$=Estimate\_BlowUpRate($t_{\rm curr}, \bm x_{\rm curr},\bm f, f_{\rm max},c_{\rm grow},N_{\rm detc}$);\; (Alg.~\ref{alg:estimateBU})
        \IF{$(-1<q<0)$}
            \STATE $h\leftarrow t^{\ast}-t_{\rm curr}$; \; $\bm x_{prev}\leftarrow \bm x_{\rm curr}$; \; $t_{\rm curr}\leftarrow t_{\rm curr}+2h$; \;
            \STATE $\bm{v}_f^{\ast}\leftarrow\bm{f}(\bm{x}_{prev})/\|\bm{f}(\bm{x}_{prev})\|$;
            \STATE $\bm x_{curr}\leftarrow \bm{x}_{curr}+\bm{v}_f^{\ast}\cdot (2ah^{(q+1)})/(q+1);$
        \ELSE
            \STATE $\bm x_{curr}\leftarrow -\bm x_{curr}$; \;
        \ENDIF
        \STATE %update $\gamma$ to preserve the direction of $\dot{\bm{x}}$:
        $\gamma\leftarrow \gamma \cdot \sign((\bm{v}_f^{\ast})^\top \bm f(\bm x_{\rm curr}))\cdot \sign((\bm{v}_f^{\ast})^\top\bm f(\bm x_{prev})) $; \;
    \ENDIF
    \IF{$E_{4}$ triggered}
        \STATE Reparameterize time: $t(s)=t_{\rm curr}+\dfrac{s}{1-s}$; \;
        \STATE Define flow: $\dot{\bm{x}} \leftarrow @(s,\bm x) \bm f(t(s),\bm x)/(1-s)^2$; \;
    \ENDIF
\ENDWHILE
\RETURN trajectory and singularity list: $[t_{all},x_{all},t^{\ast}]$.
\end{algorithmic}
\end{algorithm}

Algorithm~\ref{alg:contNewton} integrates the continuous Newton flow
\(\dot{x}=\gamma \mathbf{J}(x)^{-1}g(x)\) while monitoring norm-based thresholds.
Algorithm~\ref{alg:estimateBU} processes a short sequence of monotonic, scaled threshold hits and estimates the triplet \((q,t^\ast,a)\) characterizing the self-similar blow-up law \(\|f(t)\|\sim a(t^\ast-t)^q\).
In the implementation, function handles are denoted by the \texttt{@(x)} notation, consistent with standard MATLAB syntax.

%%%%%%%%%%%%%%%%%%
\begin{algorithm}
\caption{Estimate\_BlowUpRate $(\|\bm f\|\approx a|t-t^{\ast}|^q)$}
\label{alg:estimateBU}
\begin{algorithmic}[1]
\STATE \textbf{Input:} $t_{\rm curr}, \bm x_{\rm curr},@(t,\bm x) \bm f(t,\bm x), \bm f_{\rm max},c,N$.
%\STATE \textbf{Output:} $q,t^\ast,a,x^\ast$
\STATE $t_1\leftarrow t_{\rm curr};\quad \bm x_1\leftarrow \bm x_{\rm curr}$; \;
\FOR{($j=2:N$)}
    \STATE $f_{\rm detc} \leftarrow f_{\rm max} c^{(j-1)}$; \;
    \STATE Define Event $E_j=\|\bm f\|-f_{\rm detc}$; \;
    \STATE $[t_{j},\bm x_{j}]\leftarrow$ Integrate $ \dot{\bm{x}}=\bm f(t,\bm x)$ with Event checking $\{E_j\}$; \;
\ENDFOR
\STATE Compute $\Delta_j\leftarrow t_{j+1}-t_j$; \; (with $j=1:N-1$)
\STATE Estimate:
\(
q\leftarrow\frac{\log(c)(N-1)}{\log\Delta_{N}-\log\Delta_{1}}; \quad
t^{\ast}=t_{N+1}+\frac{c^{1/q}}{1-c^{1/q}}\,\Delta_{N},
\)
\STATE
\(
a=\dfrac{f_{\rm detc}}{(t^{\ast}-t_{N+1})^{q}};
\) \;
%\STATE (Optional) refine $(\log|a|,q,t^\ast)$ by minimizing \( \sum_{j=1}^{N}\Big(\log|f_i(t_j)|-\log|a|-q\log(t^\ast-t_j)\Big)^2, \).
\RETURN $(q,t^\ast,a,x^{N+1})$
\end{algorithmic}
\end{algorithm}

%\newpage
\section{Illustrative case studies} \label{sec:casestudies}
In this section, we demonstrate the practical performance of the proposed framework on a set of representative problems. The aim is to illustrate its ability to detect and cross finite-time blow-ups, to handle grow-up scenarios (infinite-time escape to infinity), and to navigate multiple sequential singularities. %
%and, in certain cases, escape spurious loops.
For each example, we report the CNM trajectory evolution, the continuous time-evolution of the residual equations and/or its $L^2$-norm $\|\bm{g}(\bm{x})\|_2$, the continuous time-evolution of the vector field and/or its $L^2$-norm $\|\bm{f}(\bm{x})\|_2$.
Importantly we highlight the number of singularity crossings, and the case in which the state grows-up without observing a blow-up in the vector field.
All integrations are carried out using an adaptive Runge–Kutta solver (Dormand–Prince 5(4), as implemented in \texttt{MATLAB R2025b}'s \texttt{ode45}) with absolute and relative tolerances set to $10^{-10}$ and $10^{-8}$, respectively. Blow-up detection is triggered when the norm of the vector field exceeds $f_{\max}$, and the state norm exceeds $r_{\max}$. The self-similar parameters are estimated from a sequence of $5$ successive detections, with growth factor $c=1.4$.
Integration is terminated when either $\|g(\bm{x})\|_2 < \varepsilon_g = 10^{-10}$ or $\|\bm{f}(\bm{x})\|_2 < \varepsilon_f = 10^{-8}$, indicating convergence to a root or to a stationary point of $G$, respectively.
The computations were performed on a laptop with an \texttt{Intel Core Ultra 9 275HX} processor and 32 GB of RAM.

\begin{figure}[ht!]
\centering
    % \subfigure[]{
    % \includegraphics[width=0.3\textwidth]{new_figures/case1z_z2_2_rootsqrt2_xoft.pdf}
    % }
    \subfigure[]{
    \includegraphics[width=0.4\textwidth]{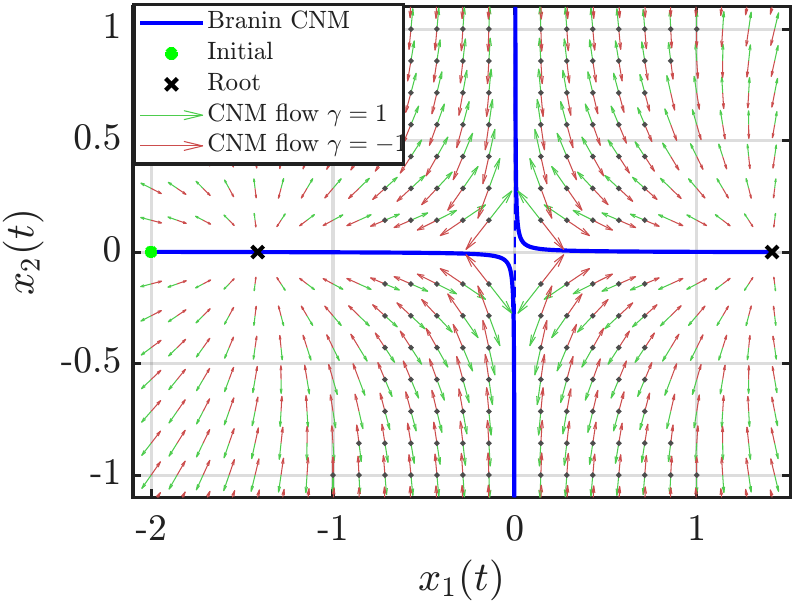}
    }
    \subfigure[]{
    \includegraphics[width=0.45\textwidth]{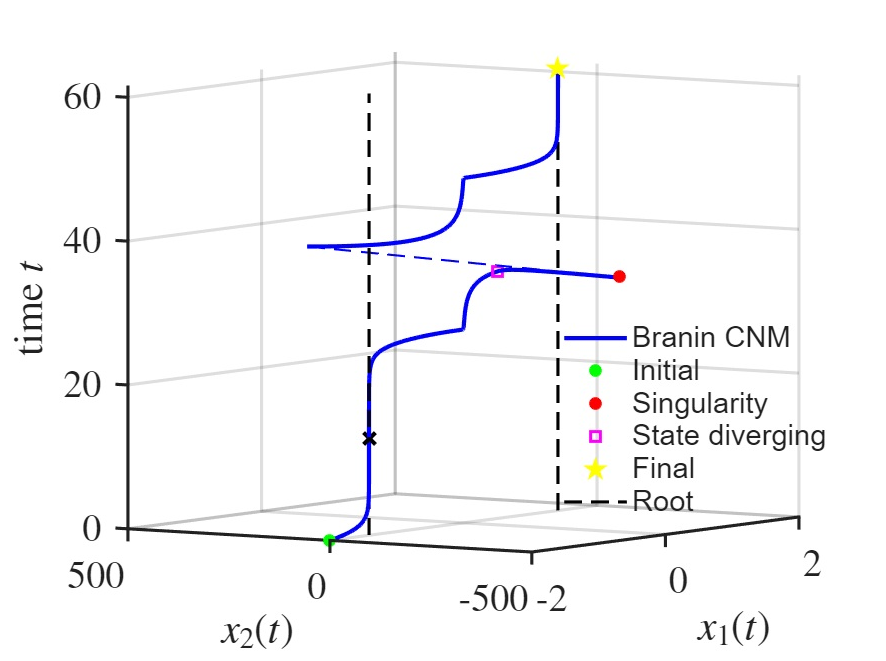}
    }
    \\
    %\hfill
    \subfigure[]{
    \includegraphics[width=0.4\textwidth]{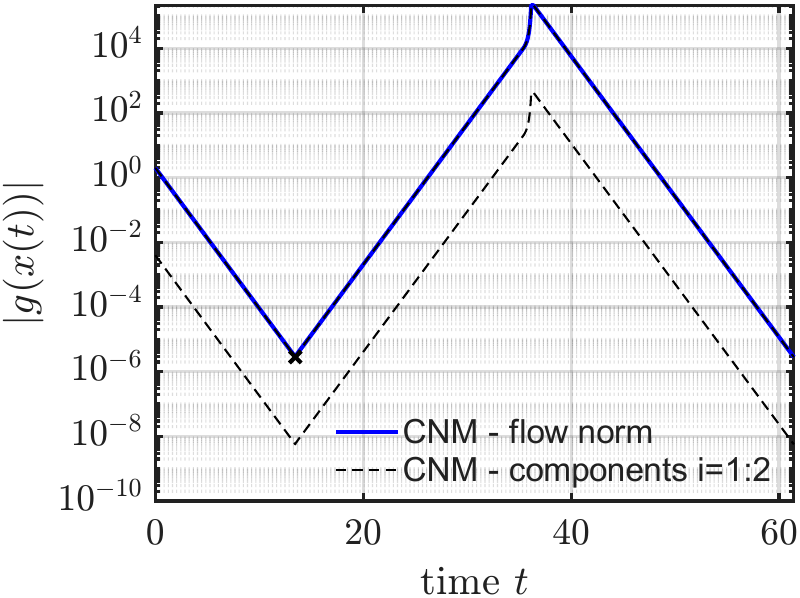}
    }
    %\hfill
    \subfigure[]{
    \includegraphics[width=0.4\textwidth]{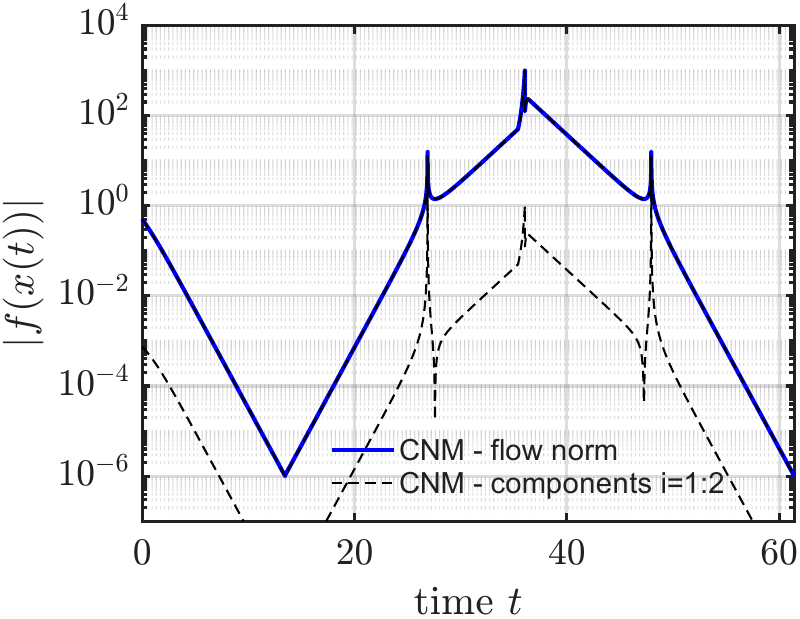}
    }
    %%\hfill
    
    \caption{\scriptsize Illustrative case study 1: complexification of $x^2-2$. Continuous Newton flow for $z^{2}-2 = 0$, $z=x_1+i\,x_2\in\mathbb{C}$. Initial condition $x_1(0)=-1$, $x_2(0)=1E-3$, $\gamma(0)=-1$.
    (a) phase portrait in the $(x_1,x_2)$-plane showing the trajectory (thick blue line) overlaid on the vector field for $\gamma=1$ (green arrows) and $\gamma=-1$ (red arrows).
    (b) Trajectory in the $(x_1, x_2, t)$-space; $x_1(t) = \Re(z(t))$ and $x_2(t) = \Im(z(t)))$. The dashed blue line indicates the compactification-based reset ($\bm z \to -\bm z$).
    (c-d) Time evolution of (c) the residual $g(z(t))$ and (d) $f(z(t))$ in a semi-logarithmic representation.
    \label{fig:trajx2complex}.}
    \vspace{0.1cm}\hrule
\end{figure}

%%%
\subsection{Case study 1: Complexification of a quadratic equation}
To illustrate the extension of the proposed framework to higher dimensions, we consider the complexification of the scalar equation
\begin{equation}
g(z) = z^2 -2, \quad z=x_1+i\, x_2 \in \mathbb{C}, \quad x_1,x_2\in \R.
\end{equation}
The roots are $z^\ast = \pm \sqrt{2} \in \mathbb{R}$. The corresponding continuous Newton flow, written in terms of the real and imaginary parts, is given by
\begin{align}
    \dot x_1(t) &= \gamma(t) \left( \frac{x_1(t)}{2} - \frac{x_1(t)}{x_1(t)^2 + x_2(t)^2} \right), \label{eq:complex_flow_x} \\
    \dot x_2(t) &= \gamma(t) \left( \frac{x_2(t)}{2} + \frac{x_2(t)}{x_1(t)^2 + x_2(t)^2} \right), \label{eq:complex_flow_y}
\end{align}
%with $\gamma(0) = 1$.
We initialize the trajectory at $z(0) = -2 + i 10^{-3}$, slightly perturbed off the real axis and $\gamma(0) = -1$. The resulting trajectories are shown in Fig.~\ref{fig:trajx2complex}.
The CNM converge first to the root $z=-\sqrt{2}$ around time $13.47$, highlighted by a black cross in Fig.~\ref{fig:trajx2complex}. Having found one solution, we flip
the orientation by setting $\gamma = +1$ and restart from a nearby point after the root extrapolating in the same direction that approached the root.

Then the real part $x_1(t)$ initially decreases toward the Jacobian singularity at $z = 0$, while the imaginary part $x_2(t)$, though initially small, grows as the trajectory approaches the singularity (Fig.~\ref{fig:trajx2complex}, panels (a-b)).

The detection framework is triggered when the state norm exceeds $r_{\max} = 100$ at $t \approx 35.40$, 
indicating that $x_2(t)$ is growing-up to infinity (\texttt{IT2}). The time transformation in Eq.~\eqref{eq:time_transform} is then applied, inducing a singularity of type $q \leq -1$. From a sequence of threshold crossings, the estimated self-similar parameters are
\begin{equation}
    \hat q = -2.539, \qquad \hat t^\ast = 36.265, \qquad \hat a = 13.377.
    \label{eq:complex_estimates}
\end{equation}
The estimated exponent $\hat q \approx -2.54$ confirms that both the state and the vector field blow-up in finite time after the pseudo-time reparametrization (\texttt{IT1}). The singularity is crossed via projective compactification: the state, which had reached approximately $(-0.462, -461.67)$ on one side of infinity, re-emerges at $(0.462, 461.67)$ on the opposite side, the orientation is updated to $\gamma = -1$, and integration resumes.
Following the crossing, the flow is initially weakly attracted toward the positive root $\sqrt{2}$, with $y(t)$ rapidly decaying toward zero. Around time $t \approx 47.95$, a transient occurs as can be seen in panel (d): the trajectory rapidly transitions from moving primarily in the $x_2$-direction to moving primarily in the $x_1$-direction. %This transient is not a singularity; the residual $|g(z(t))|$ remains finite and evolves monotonically (increasing for $\gamma=1$, decreasing for $\gamma=-1$).
%\color{black}
%
The trajectory then converges to the positive real root $z^\ast = \sqrt{2}$, with integration terminating at $t \approx 61.43$ when $\|\dot z\|_2 < 10^{-6}$ and $\|g(z)\|_2 = 2.83 \times 10^{-6}$.
The entire computation completes in approximately $0.063$ seconds.
The phase portrait in panel~(a) is particularly instructive for understanding the CNM dynamics. The Newton flow has a pole at the origin (a saddle-type equilibrium of the desingularized field),
and two additional fixed points, a source and a sink, located at $z = \pm \sqrt{2}$. The role of $\gamma$ is then fully transparent: it simply \emph{exchanges} which of the two fixed points is the source and which is the sink. %Consequently, after approaching the singularity near the origin along the unstable manifold 
%%%YGK  along whose unstable manifold ? or whose stable manifold ? 
%and performing an infinity crossing, the flow is guided from one fixed point to the other.
Consequently, the flow first approaches the origin along the real axis—the stable direction of the desingularized saddle-type equilibrium at the origin. As it gets too close, it is repelled along the imaginary axis—the corresponding unstable direction. It is along this vertical excursion that the infinity crossing is performed; this crossing is therefore situated on the path along which the flow is guided from one fixed point to the other.

\begin{figure}[ht!]
\centering
    \subfigure[]{
    \includegraphics[width=0.3\textwidth]{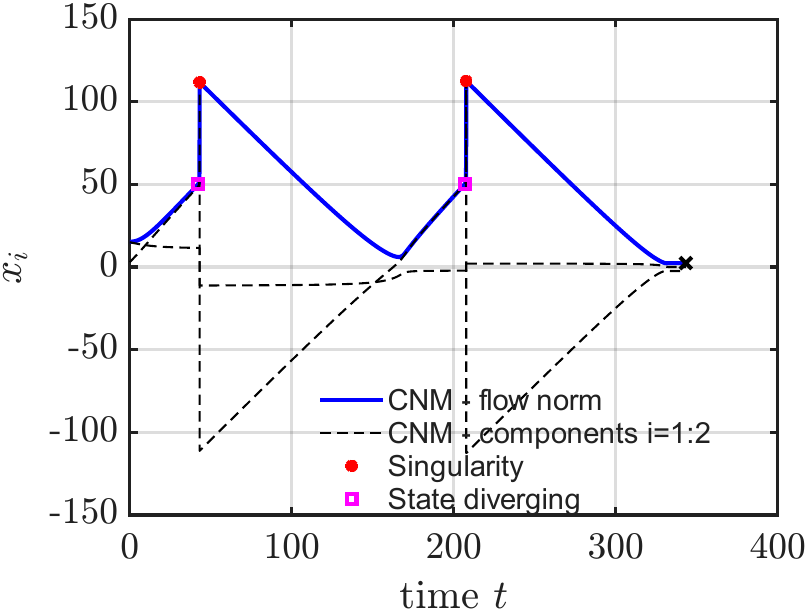}
    }
    \subfigure[]{
    \includegraphics[width=0.3\textwidth]{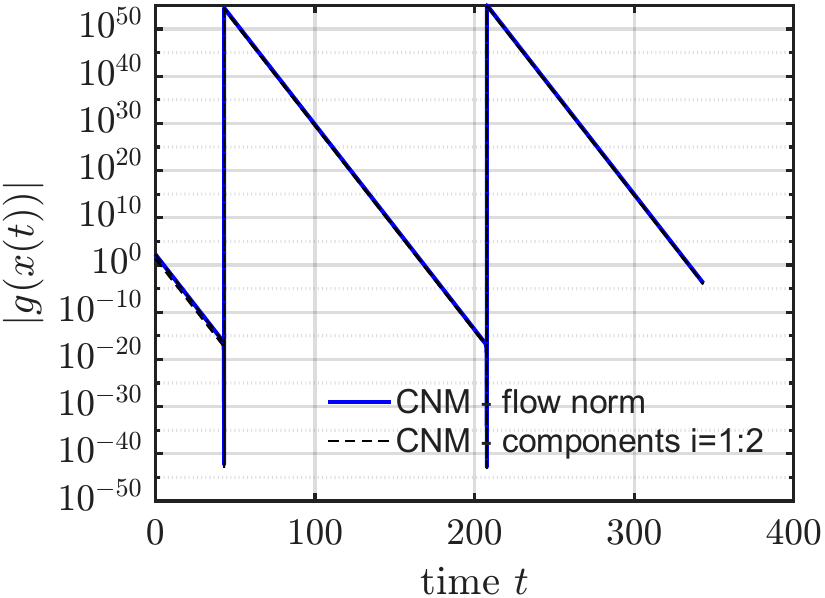}
    }
    \subfigure[]{
    \includegraphics[width=0.3\textwidth]{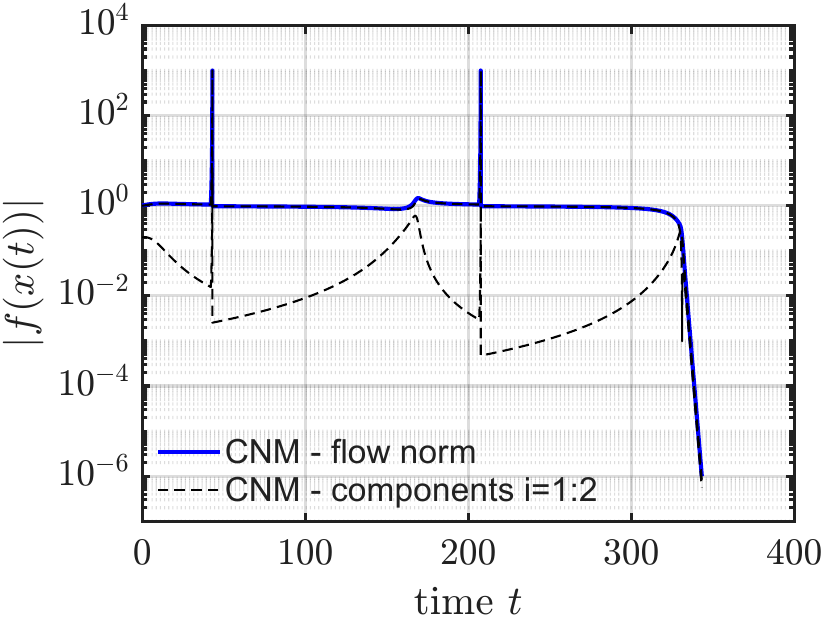}
    }
    \\
    \subfigure[]{
    \includegraphics[width=0.45\textwidth]{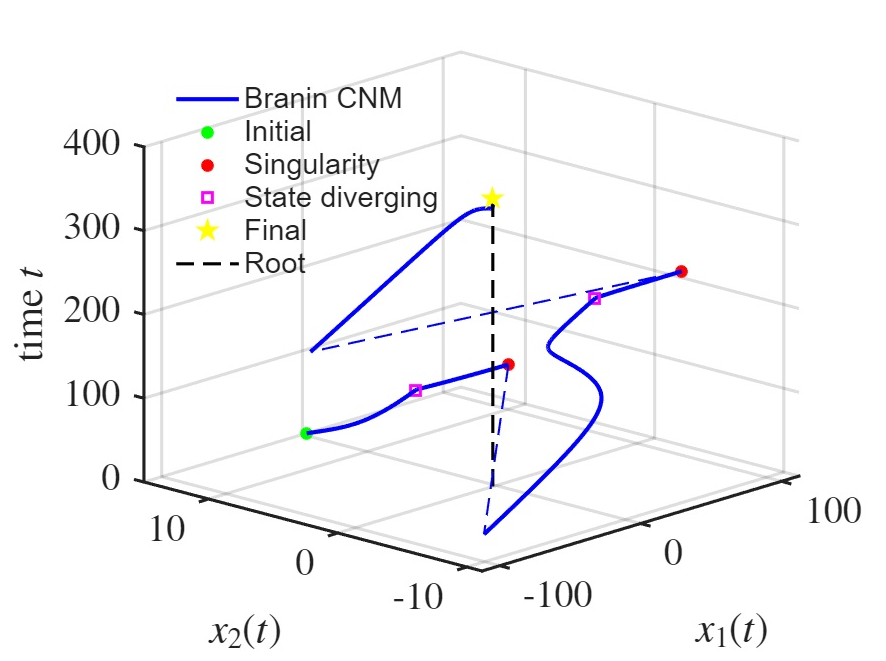}
    }
    \subfigure[]{
    \includegraphics[width=0.4\textwidth]{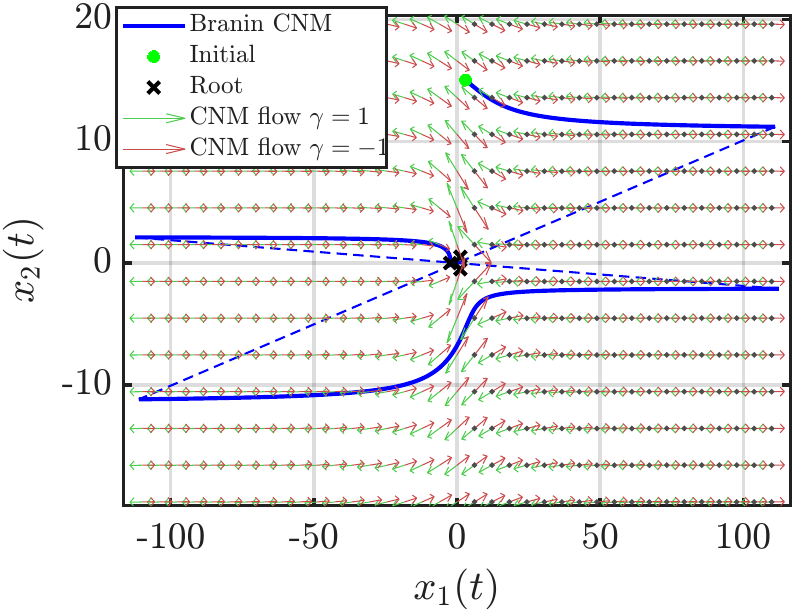}
    }
    
    %\hfill
    
    \caption{\scriptsize Illustrative case study 2: Complexification of nonlinear equation with root at infinity. Continuous Newton flow for $g(z) = e^{-z}((z^2-4)z+4)$, $z = x_1 + i x_2 \in \mathbb{C}$. Initial condition $x_1(0) = 3$, $x_2(0) = 15$, $\gamma(0) = -1$.
    (a) Evolution of the state norm $\|\bm x(t)\|_2$; (b-c) evolution of (b) $|g(z(t))|$ and (c) $|f(z(t))|$ in semi-logarithmic scale.
    (d) three-dimensional trajectory $(x_1(t), x_2(t), t)$; (e) phase portrait in the $(x_1,x_2)$-plane showing the trajectory (thick blue line) overlaid on the vector field for $\gamma=1$ (green arrows) and $\gamma=-1$ (red arrows). %(d) evolution of $\|\dot z(t)\|_2$ in semi-logarithmic scale.
    \label{fig:traj_complex2}}
    \vspace{0.1cm}\hrule
\end{figure}
\subsection{Case study 2: Complexification of nonlinear equation with root at infinity}
\label{subsec:complexification_infinity}

We now consider the complexification of the scalar equation
\begin{equation}
g(z) = e^{-z} \bigl((z^2 - 4)z + 4\bigr), \qquad z = x + i y \in \mathbb{C}, \quad x_1,x_2 \in \mathbb{R},
\label{eq:complex_expx}
\end{equation}
which in the real scalar case has a unique real root at $x^\ast \approx -2.38$ and a horizontal asymptote as $x \to +\infty$. In the complex plane, the cubic polynomial $p(z) = z^3 - 4z + 4$ has, in addition to the real root, two complex conjugate roots at $z_1^\ast \approx 1.19 + 0.51 i$ and $z_2^\ast \approx 1.19 - 0.51 i$; more generally for an $n$-th degree polynomial, there will be
$n$ such roots in the complexified problem. The complexified Newton flow is given by
\begin{align}
    \dot x_1(t) &= \gamma(t) \, \Re\left(\frac{g(x_1+ix_2)}{g'(x_1+ix_2)}\right), \label{eq:complex_expx_flow_x} \\
    \dot x_2(t) &= \gamma(t) \, \Im\left(\frac{g(x_1+ix_2)}{g'(x_1+ix_2)}\right), \label{eq:complex_expx_flow_y}
\end{align}

\noindent with $\gamma(0) = -1$.

We initialize the trajectory at $z(0) = 3 + 15 \,i $. The resulting trajectories are shown in Fig.~\ref{fig:traj_complex2}.
The trajectory crosses infinity twice, at $t \approx 43.1$ and again at $t \approx 207.7$, with estimated self-similar exponents $\hat q \approx -1.985$ in both cases, confirming that both the state and the vector field blow up in finite time ($q \leq -1$).
The orientation factor $\gamma$ remains $-1$ throughout, consistent with the preservation of the flow direction. The trajectory then converges to the real root at $t \approx 343.3$, with residual $\|f(z)\|_2 \approx 1 \times 10^{-6}$.

Due to the exponential factor $e^{-z}$, which dominates the Jacobian, the flow is predominantly horizontal: for $\gamma = -1$, it moves from left to right, and for $\gamma = 1$, from right to left (see panel (e)). In the central region near (for $|x_1|<10$), there is a contraction in the imaginary direction, which, after multiple crossings, brings the trajectory close to one of the three roots (all located near the center of the phase space). Notably, the orientation factor $\gamma$ remains $-1$ throughout the crossings, as the positivity and magnitude of $e^{-z}$ strongly affect the vector field (see panel (e)).

However, if one flips the orientation to $\gamma = 1$ in an attempt to find additional roots, the flow reverses direction, moving from right to left. In the central region, the imaginary component now grows, eventually driving the trajectory away from the roots and toward infinity. After repeated crossings, the trajectory tends to grow in the $x_2$-direction as well.

%%%%%%%%%%%%%%%%%%%%%%
\begin{figure}[ht!]
\centering
    \subfigure[]{
    \includegraphics[width=0.4\textwidth]{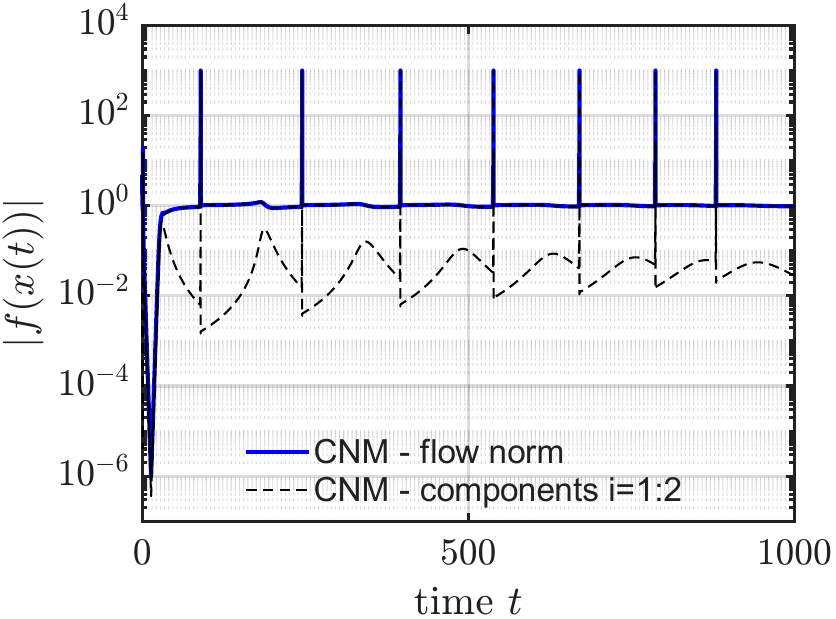}
    }
    \subfigure[]{
    \includegraphics[width=0.4\textwidth]{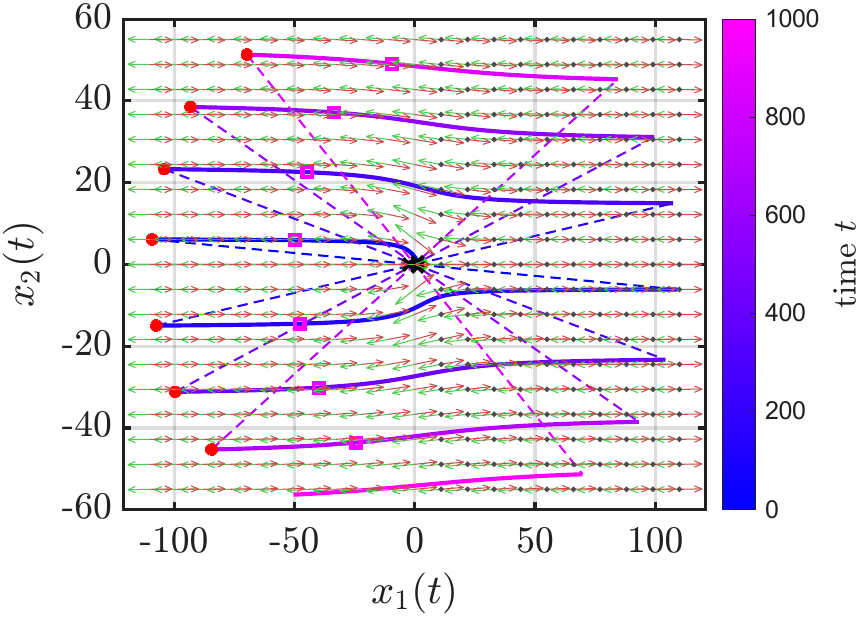}
    }
    \caption{\scriptsize Failed attempt to find additional roots after flipping $\gamma = 1$. The trajectory escapes to infinity rather than converging to a complex root. (a) $|f(x(t))|$ in semi-logarithmic scale; (b) phase portrait showing the escape path (thick blue line) overlaid on the vector field. For ease of visualization, the trajectory is colored by time (blue to magenta; see colorbar), tracing its cyclic evolution and resets at infinity. \label{fig:complex_expx_escape}}
    \vspace{0.1cm}\hrule
\end{figure}
%%%%%%%%%%%%%%%%%%%%%%%

To illustrate this behavior, and the possibility to converge to a different root, we run a separate simulation starting from $z(0) = 3 + 0.001 \, i$. With $\gamma = -1$, the CNM converges first to the complex root $z_1^* \approx 1.19 + 0.51 i$ at time $t\approx 13.65$, without encountering singularities. We then flip the orientation to $\gamma = 1$; the resulting trajectory escapes to infinity, as shown in Fig.~\ref{fig:complex_expx_escape}.
Whenever the system exceeds the maximum norm $r_{\max} = 50$, this event triggers a time reparametrization that accelerates the trajectory toward infinity and detects a singularity when $\|f\|_2 > 10^3$. As can be observed Fig.~\ref{fig:complex_expx_escape}(a), many of such events occurs, as reflected in the multiple peaks observed.
The phase portrait in panel~(b) illustrates the trajectory exiting the root $z_1^*$ and proceeding toward $x_1 \to -\infty$, with resets occurring near $x_1 = -100$. Upon reset, the trajectory reappears at $x_1 \approx 100$, while the $x_2$ coordinate is flipped from approximately $+5$ to $-5$. Notably, the sign of $\gamma$ remains positive at the crossing to preserve the horizontal direction of the vector field. As the trajectory evolves back toward $x_1 = -100$, $x_2$ decreases to roughly $-15$, is subsequently reset to $+15$, and continues to grow to about $+25$ before the next crossing. This process repeats indefinitely, with $x_2$ alternating sign and growing monotonically in magnitude with each cycle.
Contrary to the case $\gamma=-1$, the contraction in the $x_2$ direction becomes expansion.
This behavior highlights a limitation of the Branin-Newton flow: reversing the orientation does not guarantee convergence to another solution, as the flow may instead be pushed away from the roots and toward infinity.
\begin{remark}
    Alternative, complementary strategies exist for identifying multiple roots without flipping $\gamma$ to $+1$. Deflation techniques~\cite{farrell2015deflation} remove already-discovered roots by modifying the system, allowing the flow---with $\gamma \equiv -1$---to converge to previously undiscovered solutions. Note that deflation alters the system itself, whereas the Branin-Newton preserve the original system and change the direction of integration.
    Note also, that in the example $g(z) = e^{-z}p(z)$ considered above, the finite roots are determined entirely by the polynomial $p(z)$; the exponential factor introduces roots at infinity but does not affect the finite roots. In such cases, one could simplify the system ``by hand", symbolically removing the exponential factor and solving $p(z) = 0$. More generally, whenever a multiplicative factor that does not vanish at finite roots can be identified, it may be advantageous to factor it out before applying the CNM.
\end{remark}

%However, it also highlights that reversing the orientation does not guarantee convergence to another root, as the flow may escape to infinity instead.

\begin{figure}[ht!]
\centering
    \subfigure[]{
    \includegraphics[width=0.3\textwidth]{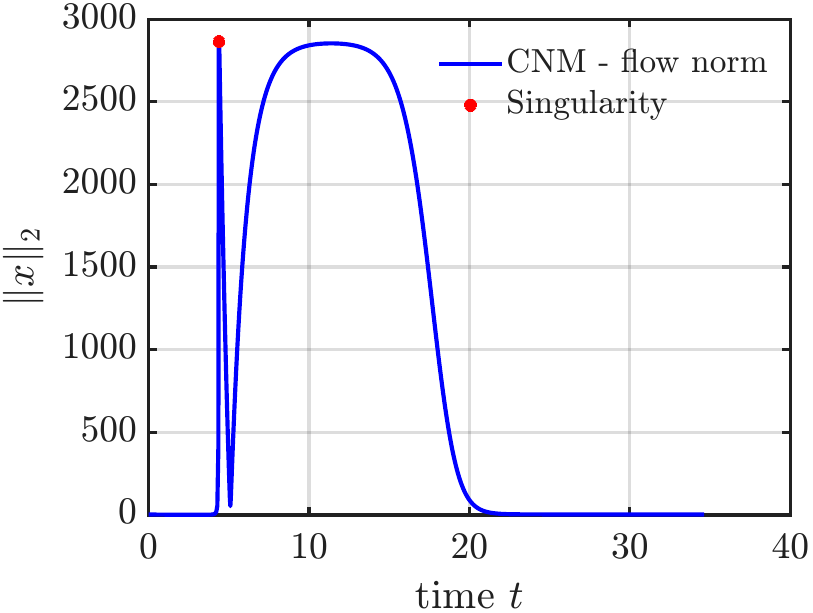}
    }
    \subfigure[]{
    \includegraphics[width=0.3\textwidth]{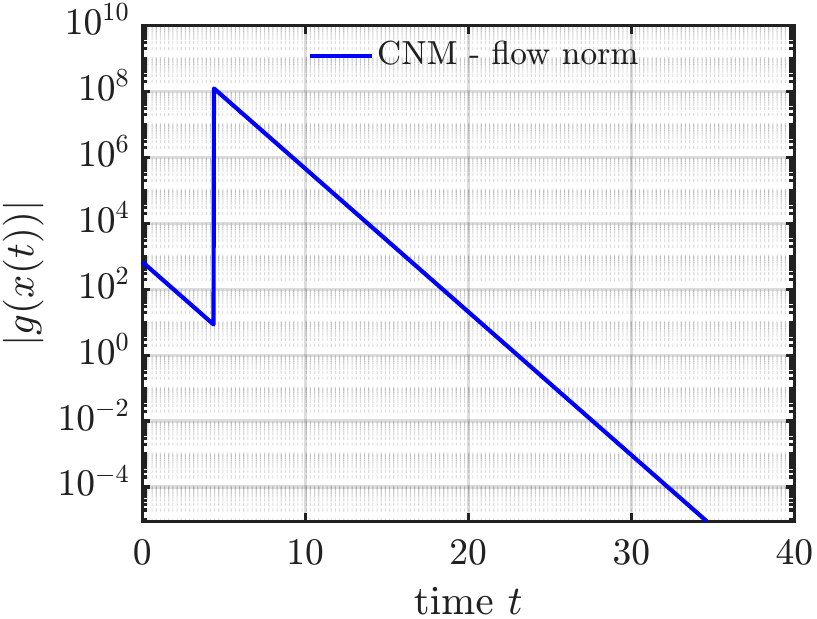}
    }
    %\hfill
    \subfigure[]{
    \includegraphics[width=0.3\textwidth]{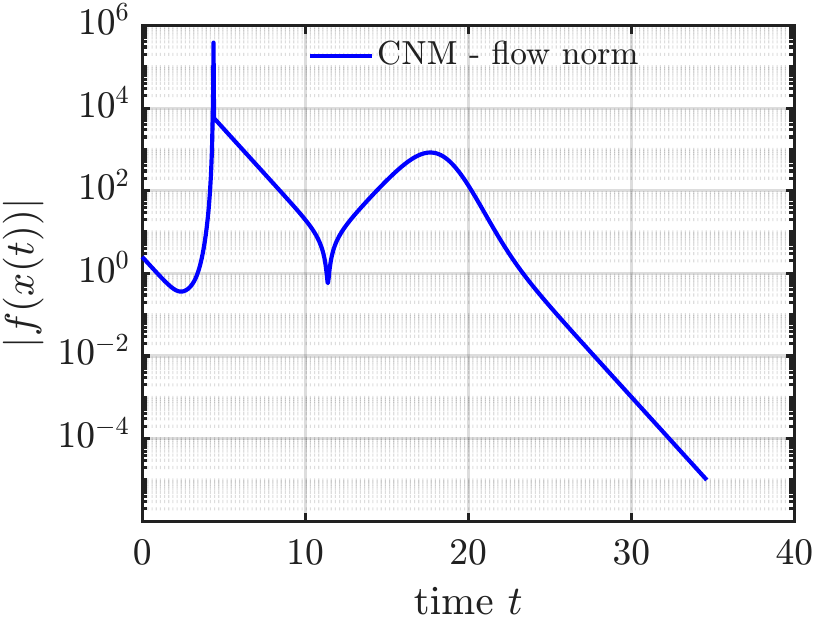}
    }
    %\hfill
    \\
    \subfigure[]{
    \includegraphics[ width=0.45\textwidth]{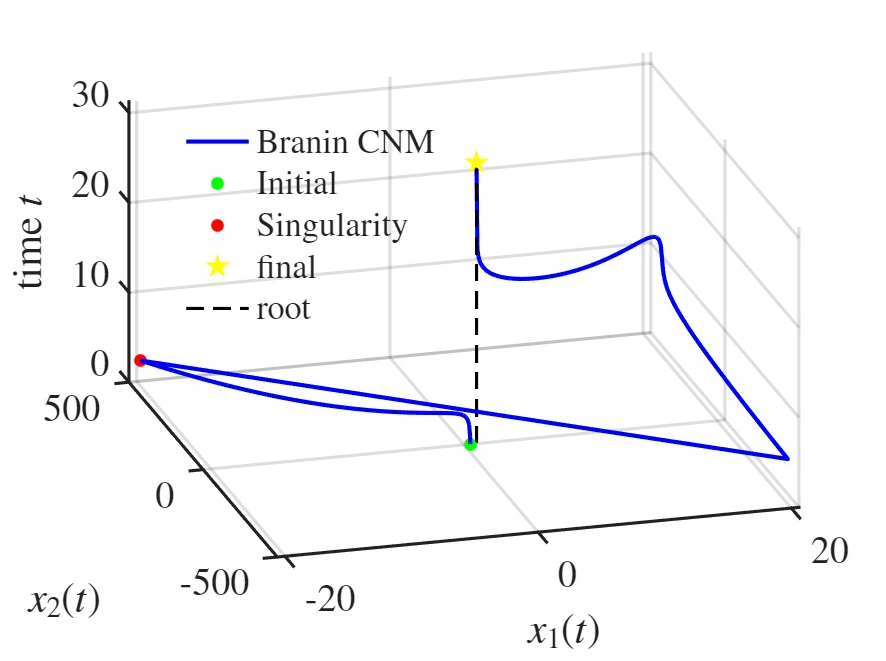}
    }
    \subfigure[]{
    \includegraphics[ width=0.4\textwidth]{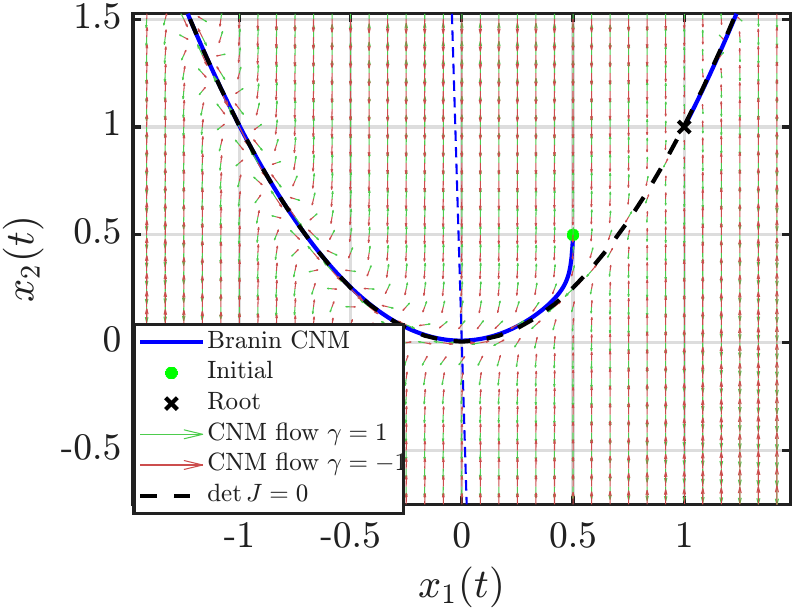}
    }
    \caption{\scriptsize Illustrative case study 3: 10-dimensional Rosenbrock problem. Continuous Newton flow for $\mathbf{g}(\bm{x}) = \nabla f(\bm{x}) $ with $d = 10$, initial condition $\bm{x}(0) = (1/2,\ldots,1/2)$, and $\gamma(0) = -1$.
    (a) Evolution of the state norm $\|\bm{x}(t)\|_2$; 
    (b) evolution of the residual norm $\|\mathbf{g}(\bm{x}(t))\|_2$ in logarithmic scale; (c) evolution of the vector field norm $\|\dot{\bm{x}}(t)\|_2$ in logarithmic scale, showing the finite-time blow-up at $t \approx 4.37$.
    (d) three-dimensional trajectory $(x_1(t), x_2(t), t)$, for indicatively $x_1$ and $x_2$ components;
    (e) zoom-in phase portrait for indicatively $x_1$ and $x_2$ components.
    \label{fig:rosenbrock}}
    \vspace{0.1cm}\hrule
\end{figure}
\subsection{Case study 3: extended Rosenbrock problem}
\label{subsec:rosenbrock}

To demonstrate the applicability of the proposed framework to higher-dimensional problems, we consider the extended Rosenbrock function~\cite{more1981testing}, a classic benchmark in optimization. % in its standard (non-overlapping) formulation. 
The function is defined as
\begin{equation}
    G(\bm{x}) = \sum_{i=1}^{d/2} \left[ (1 - x_{2i-1})^2 + 100 (x_{2i} - x_{2i-1}^2)^2 \right],
    \label{eq:rosenbrock}
\end{equation}
which we treat as an optimization problem whose solution is the global minimum. Reformulating this optimization problem as a system of nonlinear equations, we take the gradient $\mathbf{g}(\bm{x}) = \nabla G(\bm{x}) = \mathbf{0}$, whose unique solution is the global minimum $\bm{x}^\ast = \mathbf{1} \in \mathbb{R}^{d}$. The Jacobian $\mathbf{J}(\bm{x}) = \nabla^2 G(\bm{x})$ is block-diagonal with $d/2$ blocks of size $2 \times 2$, each corresponding to a pair of variables $(x_{2i-1}, x_{2i})$. Here we select $d=10$.

The continuous Newton flow is simulated with initial condition $\bm{x}(0) = (1/2,\dots,1/2)$ and orientation $\gamma(0) = -1$.
From this initial guess, the trajectory is initially attracted toward the narrow curved valley characteristic of the Rosenbrock function, where the Jacobian becomes ill-conditioned.
As the trajectory progresses, the Jacobian becomes singular at a finite time, causing the vector field to blow-up. The detection framework is triggered when the vector field norm exceeds $f_{\max} = 10^5$ at $t \approx 4.37$, at which point the state norm has already grown to $\|\bm{x}\|_2 \approx 1.17 \times 10^3$. The estimated self-similar parameters are
\begin{equation}
    \hat q = -2.808, \qquad \hat t^\ast = 4.389, \qquad \hat a = 2.253.
    \label{eq:rosenbrock_estimates}
\end{equation}
The estimated exponent $\hat q \approx -2.81$ places this singularity in the regime $q \leq -1$, indicating that both the state and the vector field blow-up in finite time. The singularity is crossed using the projective compactification strategy: the state is reset by sign flipping, and integration resumes.
%In this case, the orientation factor $\gamma$ is not updated after crossing, as the sign of the vector field remains unchanged across the singularity; that is, $\sign(\mathbf{v}_f^\top \mathbf{f}_{\text{before}}) \cdot \sign(\mathbf{v}_f^\top \mathbf{f}_{\text{after}}) > 0$.
Integration then continues, and the trajectory converges to the global minimum $\bm{x}^\ast = \mathbf{1}$, with the residual satisfying $\|\mathbf{g}(\bm{x})\|_2 = 9.12 \times 10^{-6}$ and $\|\dot{\bm{x}}\|_2 < 10^{-5}$ at termination. The integration terminates at $t \approx 34.63$, with total execution time of approximately $0.788$ seconds.

The trajectories are shown in Fig.~\ref{fig:rosenbrock}. Panel (a) displays the evolution of the state norm $\|\bm{x}(t)\|_2$, which grows rapidly as the trajectory approaches the singularity and then converges to $\bm{x}^\ast=\mathbf{1}$ ($\|\bm{x}^\ast\|_2 = \sqrt{10}$) after crossing. Panel (d) provides an indicative three-dimensional view of the trajectory ($x_1,x_2,t$). 
%\color{blue!90!black}
The transient growth of $\|x(t)\|$ in panel~(a) is caused by the trajectory crossing infinity: the $x_2$-component transitions from approximately $+500$ to $-500$. At the reset point, the vector field is nearly vertical (aligned with $x_2$), as shown in panel~(e). This drives $x_2$ rapidly from large negative values through zero and then upward, while $x_1$ remains nonzero. From this state, the trajectory gradually approaches the Rosenbrock parabolic manifold, along which it converges to the root. The dynamics in this regime are particularly stiff. A zoom of the phase portrait near the root is shown in panel~(e).
Notably, The trajectory is only weakly attracted to the Rosenbrock parabola, which lies close to, but not exactly on, the singular locus $\det(\mathbf{J}(\bm{x}))=0$, given by $y=x^{2}+1/200$, and eventually arrives at a detected singularity.
This example illustrates that the heuristic resetting mechanism can generate new trajectories that bear little relation to the path followed prior to the growth of the state.
%\color{black}
Panel (b) shows the residual norm $\|\mathbf{g}(\bm{x}(t))\|_2$ in logarithmic scale: the residual initially decreases, then spikes as the trajectory passes through the singularity, and finally decays monotonically to below $10^{-5}$ upon convergence to the root. Panel (c) displays the vector field norm $\|\dot{\bm{x}}(t)\|_2$ in logarithmic scale, where the vertical asymptote at $t \approx 4.37$ confirms the finite-time blow-up of the vector field. %This example demonstrates that the proposed framework scales to higher-dimensional problems and handles singularities arising from ill-conditioned Jacobians in the Rosenbrock function, successfully guiding the trajectory through the singularity and to the global minimum.

\subsection{Case study 4: Steady-state of the Liouville--Bratu--Gelfand PDE}
\label{subsec:bratu}

To demonstrate the applicability of the proposed framework to large-scale systems arising from the discretization of partial differential equations, we consider the one-dimensional Liouville--Bratu--Gelfand problem
\begin{equation}
    \frac{\partial u}{\partial t} = \frac{\partial^2 u}{\partial x^2} + \lambda e^{u}, \qquad x \in (0,1),
    \label{eq:bratu_pde}
\end{equation}
subject to homogeneous Dirichlet boundary conditions $u(0) = u(1) = 0$.
\begin{figure}[ht!]
\centering
    \subfigure[]{
    \includegraphics[width=0.4\textwidth]{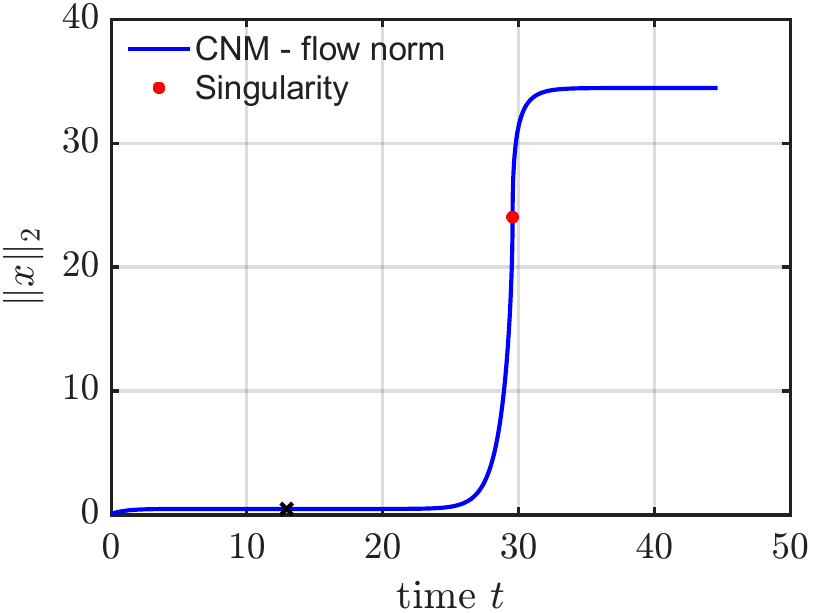}
    }
    %\hfill
    \subfigure[]{
    \includegraphics[width=0.45\textwidth]{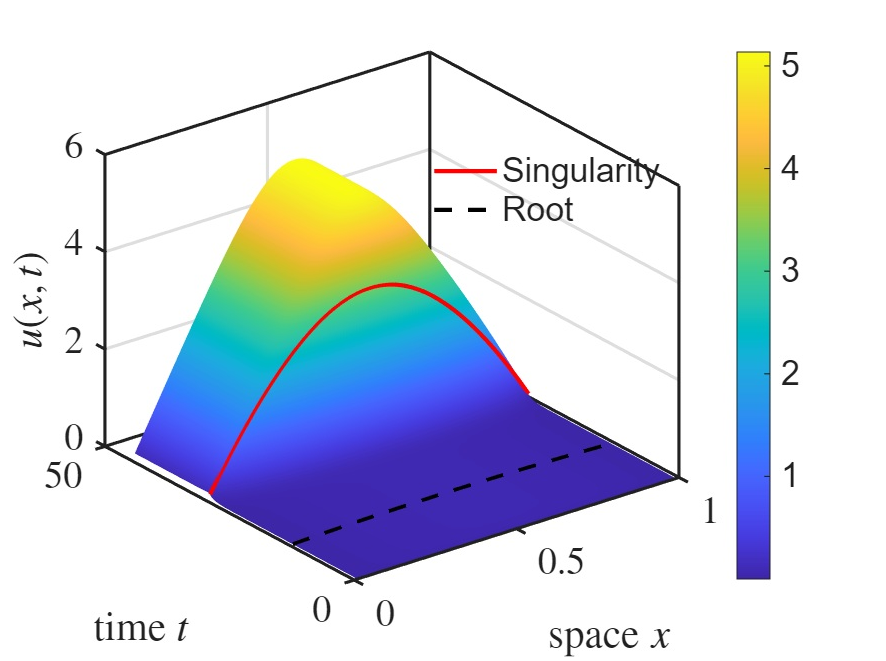}
    }
    \\
    \subfigure[]{
    \includegraphics[width=0.4\textwidth]{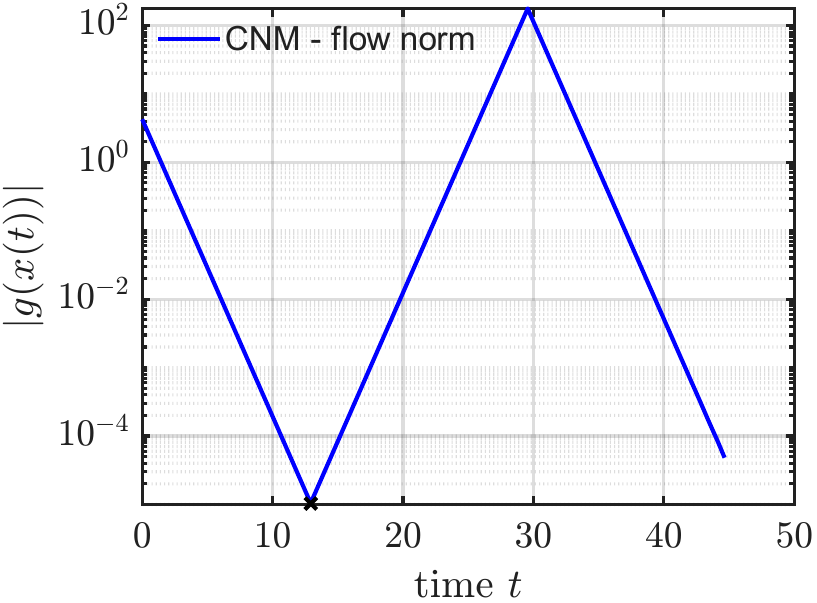}
    }
    %\hfill
    \subfigure[]{
    \includegraphics[width=0.4\textwidth]{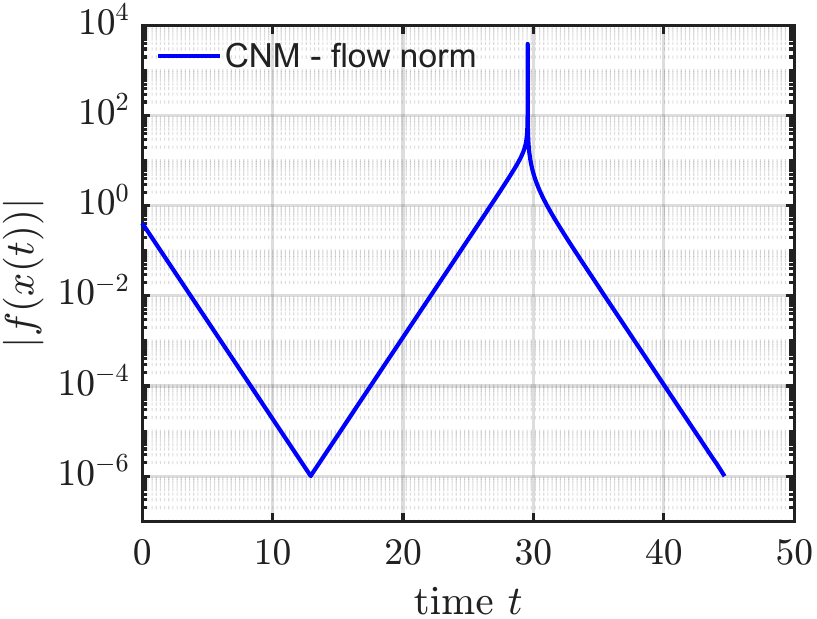}
    }
    \caption{\scriptsize Illustrative case study 4: Liouville--Bratu--Gelfand PDE discretized with $n = 101$ interior points. Continuous Newton flow for the system in Eq.~\eqref{eq:bratu_discrete} with $\lambda = 0.5$, initial profile $u(x,0) = 0.04x(1-x)$, and $\gamma(0) = -1$.
    (a) Evolution of the state norm $\|\mathbf{u}(t)\|_2$; (b) surface plot of the PDE profile $u(x,t)$ under the continuous Newton flow, with the singularity profile (red line) marking the crossing of the separatrix; (c) evolution of the residual norm $\|\mathbf{g}(\mathbf{u}(t))\|_2$ in logarithmic scale; (d) evolution of the vector field norm $\|\dot{\mathbf{u}}(t)\|_2$ in logarithmic scale.
    \label{fig:bratu}}
    \hrule
\end{figure}
The steady-state equation reads
\begin{equation}
    u_{xx} + \lambda e^{u} = 0,
    \label{eq:bratu_steady}
\end{equation}
which, for $\lambda < \lambda^* \approx 3.51$, admits two solutions: a stable lower branch and an unstable upper branch. We choose $\lambda = 0.5$ and discretize the steady-state equation using second-order centered finite differences on a uniform grid with $n = 101$ interior points. Denoting $u_i = u(x_i)$ for $x_i = i/(n+1)$, the resulting system of nonlinear algebraic equations is
\begin{equation}
    g_i(\mathbf{u}) = \frac{u_{i-1} - 2u_i + u_{i+1}}{h^2} + \lambda e^{u_i} = 0, \qquad i = 1, \ldots, n,
    \label{eq:bratu_discrete}
\end{equation}
with $u_0 = u_{n+1} = 0$ and $h = 1/(n+1)$.

The continuous Newton flow is initialized with the smooth profile $u(x,0) = 0.04\, x (1-x)$ and orientation $\gamma(0) = -1$.
The resulting trajectories are shown in Fig.~\ref{fig:bratu}.
The trajectory of the PDE profile first evolves toward the stable solution, that it is reached around time $t=12.93$ with a tolerance on $\|f\|=10^{-6}$, denoted by dashed black line in panel (b).
We then surpass the root along the direction of approach, increasing the magnitude above the stable branch, and set $\gamma = 1$.
With $\gamma = 1$, the root becomes unstable and the profile grows, eventually crossing the separatrix between the stable and unstable solution branches (red line in panel (b)). The detection framework is triggered when $\|f\|_2$ exceeds $f_{\max} = 10^3$ at $t \approx 29.56$ (\texttt{IT1}), with state norm $\|\mathbf{u}\|_2 \approx 23.94$. The estimated self-similar parameters are
\begin{equation}
    \hat q = -0.498, \qquad \hat t^\ast = 29.558, \qquad \hat a = 7.946.
    \label{eq:bratu_estimates}
\end{equation}
The estimated exponent $\hat q \approx -0.5$ places this singularity in the regime $-1 < q < 0$, indicating that the vector field blows up in finite time while the state remains bounded. Integration then resumes and the trajectory converges to the unstable steady-state solution, with integration terminating at $t \approx 44.66$ when $\|f\|_2 < 10^{-6}$ and $\|g(\mathbf{u})\|_2 = 4.82 \times 10^{-5}$.
%\color{black}

This example demonstrates that the proposed framework scales to large-scale systems arising from PDE discretizations and successfully handles singularities encountered when crossing separatrices between stable and unstable solution branches.

\subsection{Case study 5: Brent counterexample}
\label{subsec:brent}

We now consider the two-dimensional system proposed by Brent \cite{brent1972davidenko} as a counterexample to Branin's conjecture that the CNM with sign switching is globally convergent. The system is defined by
\begin{equation}
    \mathbf{g}(\bm{x}) = 
    \begin{bmatrix}
        4 (x_1 + x_2) \\
        4 (x_1 + x_2) + (x_1 - x_2)\left((x_1 - 2)^2 + x_2^2 - 1\right)
    \end{bmatrix}.
    \label{eq:brent}
\end{equation}
\begin{figure}[ht!]
\centering
    % \subfigure[]{
    % \includegraphics[width=0.4\textwidth]{new_figures/case9_brent_xoft.pdf}
    % }
    \subfigure[]{
    \includegraphics[width=0.45\textwidth]{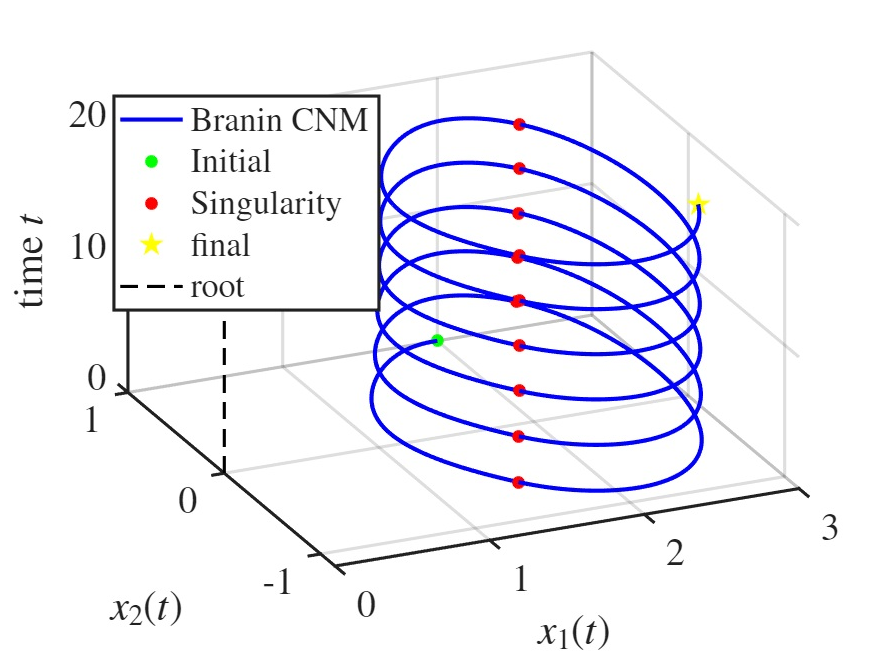}
    }
    %\hfill
    \subfigure[]{
    \includegraphics[width=0.4\textwidth]{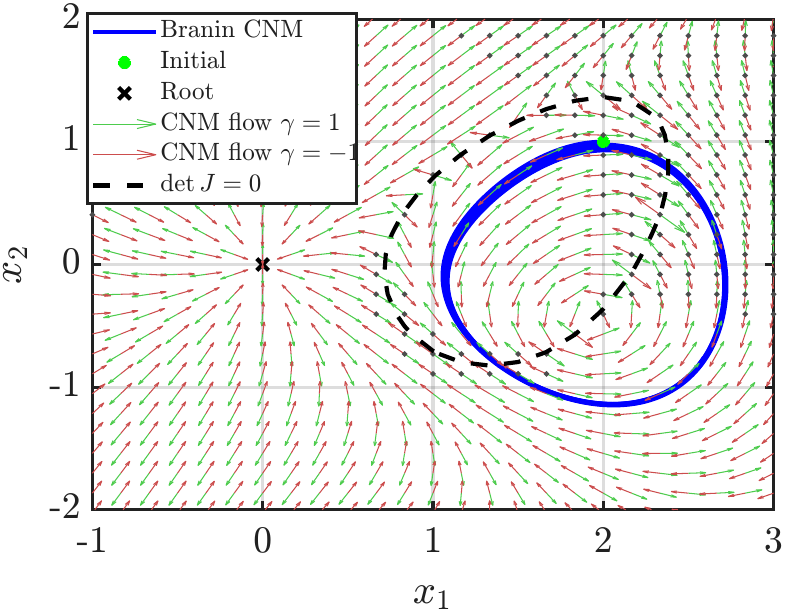}
    }
    \\
    \subfigure[]{
    \includegraphics[width=0.4\textwidth]{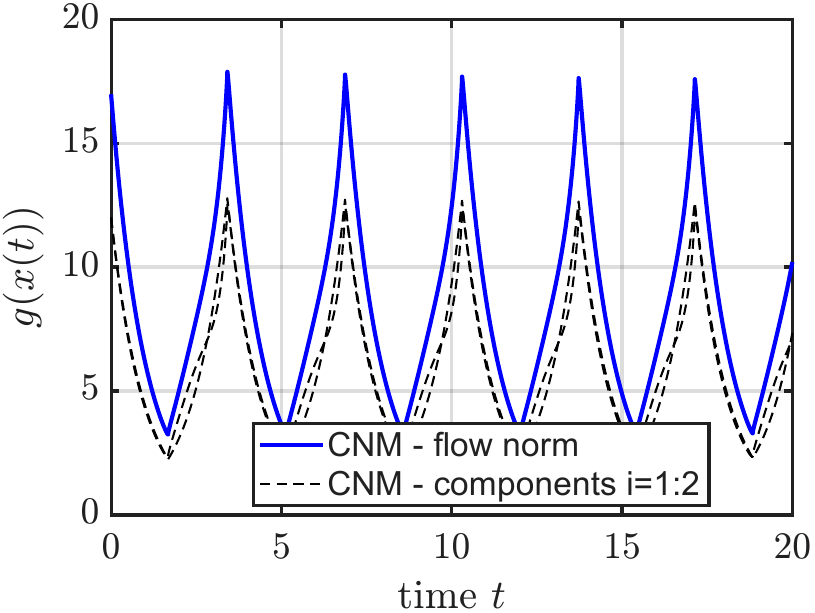}
    }
    %\hfill
    \subfigure[]{
    \includegraphics[width=0.4\textwidth]{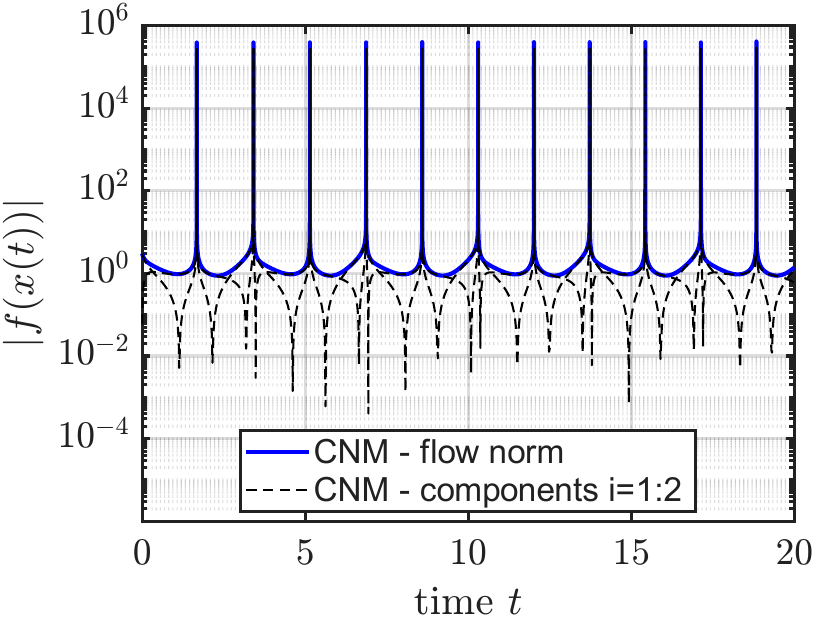}
    }
    \caption{\scriptsize Illustrative case study 5: Brent counterexample. Continuous Newton flow for $\mathbf{g}(\bm{x})$ in Eq.~\eqref{eq:brent} with $a=4$, $b=2$, $c=1$, initial condition $\bm{x}(0)=(2,1)$, and $\gamma(0)=-1$.
    %(a0) Time evolution of $x_1(t)$ and $x_2(t)$;
    (a) three-dimensional trajectory $(x_1(t), x_2(t), t)$;
    (b) phase portrait in the $(x_1,x_2)$-plane showing the trajectory (thick blue line) overlaid on the vector field for $\gamma=1$ (green arrows) and $\gamma=-1$ (red arrows).
    (c) evolution of $\|\mathbf{g}(\bm{x}(t))\|_2$ in logarithmic scale; (d) evolution of $\|\dot{\bm{x}}(t)\|_2$ in logarithmic scale, showing repeated finite-time blow-ups of the vector field.
    \label{fig:brent}}
    \vspace{0.1cm}\hrule
\end{figure}
The system has a unique root at $\bm{x}^\ast = (0,0)$. From the initial condition $\bm{x}(0) = (2,1)$ and orientation $\gamma(0) = -1$, Branin's original sign-switching strategy produces a spurious periodic orbit that cycles indefinitely without ever reaching the root. This failure is a direct consequence of the saddle-type structure induced in the desingularized flow by an extraneous singularity located near the trajectory: at each crossing, the trajectory is redirected by this saddle in a way that traps it in a bounded region containing no root.

The proposed framework is applied to this problem with $f_{\max} = 10^4$ and integrated up to $t_{\text{end}} = 20$. The results are shown in Fig.~\ref{fig:brent}. The trajectory encounters a sequence of finite-time blow-ups, each corresponding to a singularity where the vector field diverges while the state remains bounded ($-1 < q < 0$). At each detected singularity, the self-similar parameters are estimated and the crossing is performed via the integral correction in Eq.~\eqref{eq:crossing_formula}, with the orientation $\gamma$ updated according to \eqref{eq:sigma_update}. In all cases, the estimated exponent is consistently $\hat q \approx -0.5$, confirming that the singularities are of the type where the vector field blows-up but the state remains finite.

%Panel (a0)  shows the time evolution of $x_1(t)$ and $x_2(t)$; the trajectory oscillates between bounded regions, with each oscillation terminated by a finite-time blow-up of the vector field.
Panel (a) of Fig.~\ref{fig:brent} provides a three-dimensional view of the trajectory $(x_1(t), x_2(t), t)$, illustrating the periodic structure of the flow.
Panel (b) shows the phase portrait in the $(x_1,x_2)$-plane, where the trajectory (thick blue line) is overlaid on the vector field for $\gamma=1$ (green arrows) and $\gamma=-1$ (red arrows), together with the curve $\det J = 0$ (black line) marking the singularities. The trajectory repeatedly crosses this singularity curve, alternating between the two regions of the plane and producing a piecewise continuous periodic orbit.
Panel (c) displays the residual norm $\|\mathbf{g}(\bm{x}(t))\|_2$ in logarithmic scale, which spikes at each singularity crossing. Panel (d) shows the vector field norm $\|\dot{\bm{x}}(t)\|_2$ in logarithmic scale, where the vertical asymptotes %at $t \approx 1.67, 3.41, 5.14, \ldots$
confirm the repeated finite-time blow-ups.

A total of eleven singularities are detected and crossed before the integration reaches the final time $t_{\text{end}} = 20$. At each crossing, the estimated exponent is $\hat q \approx -0.5$, placing the singularity in the regime $-1 < q < 0$ where the vector field blows-up while the state remains bounded. %The integral correction in Eq.~\eqref{eq:crossing_formula} is applied to jump across each singularity, and the orientation is updated according to \eqref{eq:sigma_update}.
Despite successfully crossing all singularities, the trajectory remains trapped in the periodic orbit and does not converge to the root by $t = 20$. %the final residual is $\|\mathbf{g}\|_2 \approx 10.2$, far above the tolerance of $10^{-7}$.
This behavior is consistent with Brent's original observation: the extraneous singularities induce a spurious periodic orbit from which the flow cannot escape, regardless of the crossing strategy.
%
%\color{blue!90!black}
The phase portrait in panel~(b) offers insight into the trapping mechanism. The curve $\det J = 0$ forms a closed loop in the $(x_1,x_2)$-plane. The trajectory repeatedly intersects this curve; each intersection triggers a singularity crossing, after which the flow re-emerges on the other side of the singular set. However, rather than escaping the trapping region, the trajectory is attracted back to a different point on the same singular curve, forming a closed periodic orbit. Thus, while our infinity-crossing framework successfully navigates each individual singularity, it does not alter the underlying topology of the flow: the singular set constitutes a trapping structure from which the Branin Newton flow cannot escape, regardless of the crossing strategy.
%\color{black}

This example illustrates both the capabilities and the limitations of the proposed framework. On one hand, the method successfully detects and crosses every singularity encountered along the trajectory, demonstrating its robustness in handling repeated finite-time blow-ups of the vector field. On the other hand, the trajectory remains confined to the periodic orbit, confirming that the failure to converge is not due to an inability to cross singularities, but rather to the intrinsic geometry of the flow: the extraneous singularities form a trapping structure from which the continuous Newton flow cannot escape.  

\section{Conclusions}
\label{sec:conclusions}
We have developed a systematic framework for detecting and crossing singularities in CNM. The central idea is to equip the Branin Newton flow with an on-the-fly self-similar blow-up detection and crossing procedure that operates directly on the real ODE flow.
For the scalar continuous Newton flow, we prove under the stated hypotheses (Lemma~\ref{lem:m_zero_taylor} and Lemma~\ref{lem:asymptote}) that the vector field exhibits asymptotic self-similar behavior near any singularity encountered—whether at a local extremum or inflection point of $g$ (yielding $q \in (-1,0)$), or at a horizontal asymptote (yielding $q \leq -1$). In higher dimensions, the directional self-similar ansatz is justified by the SVD analysis of Section~\ref{subsec:singularity_structure_nd}, which shows that near a simple (single-direction) Jacobian singularity the flow asymptotically aligns with the right singular vector corresponding to the smallest singular value, reducing the dynamics to an effective one-dimensional problem. The self-similarity exponent $q$, the blow-up time $t^\ast$, and the amplitude $a$ are estimated from a short sequence of threshold crossings via linear regression.
The estimated exponent distinguishes two qualitatively different scenarios: one in which only the vector field diverges while the state remains finite ($-1 < q < 0$), and one in which the state itself blows-up in finite time ($q \leq -1$). For the former, we apply an explicit integral correction; for the latter, we employ projective compactification and sign flipping. In addition, we handle infinite-time blow-up through a suitable time reparameterization. The extension to higher dimensions is achieved through a directional self-similar detection strategy that monitors the norm of the vector field along its instantaneous direction, reducing the detection and estimation to a single scalar problem regardless of dimension.

This projective reset operation finds a natural interpretation within the broader geometric framework of projective compactification, where the real line is extended by adding a point at infinity to form the projective space $\mathbb{RP}^1$~\cite{wright1985finding}. 

For scalar ODEs such as $\dot{x}=x^2$, this compactification naturally allows the trajectory to pass through infinity and re-emerge from the opposite direction, consistent with the projective identification of $+\infty$ and $-\infty$ as the same point~\cite{kevrekidis2017infinity}. While the present work adopts a self-similar, data-driven approach to detect and cross singularities on the fly, this geometric perspective provides a rigorous underpinning for the projective reset operation used in our algorithm when the state itself diverges ($q \leq -1$). The connection also highlights a key limitation: for ODEs whose blow-up leads to complex-valued solutions beyond the singularity, as in $\dot{x}=x^3$, the real projective compactification alone is insufficient, and one must either extend to the complex domain or resort to alternative crossing strategies.
A promising direction for future work is to integrate techniques from numerical algebraic geometry, such as those implemented in Bertini~\cite{bates2013numerically}, to reduce the heuristic nature of the projective reset and to systematically handle complex blow-up regimes where the real trajectory ceases to exist.
We also note that for PDEs discretized on a spatial grid, the self-similarity observed in the discretized ODE system may be accompanied by spatial rescaling in the underlying PDE~\cite{kavousanakis2026going}. While our method detects and handles singularities at the ODE level, whether these correspond to genuine PDE self-similarity or are artifacts of discretization is problem-dependent and remains an interesting direction for future investigation.
%\color{black}

The methodology is demonstrated on a range of problems of increasing complexity, from scalar nonlinear equations to the steady-state solution of the Liouville--Bratu--Gelfand PDE discretized on a fine spatial grid. The method successfully identifies and crosses multiple sequential singularities, handles complexified problems where the trajectory must pass through a singularity, and recovers both stable and unstable solution branches in PDE applications. The Brent counterexample illustrates both the capabilities and the limitations of the approach: while the framework successfully detects and crosses every singularity encountered along the trajectory, the flow might remain trapped in the spurious periodic orbit, confirming that the failure to converge is not due to an inability to cross singularities, but rather to the intrinsic topology of the vector field.
More generally, the method cannot overcome fundamental trapping structures from which the continuous Newton flow cannot escape, regardless of the crossing strategy. Seeking to bypass this limitation, potentially through considering continuum variants of Nesterov's method~\cite{weijie2016differential,ward2017toolkitsteadystatesnonlinear,subramanian2021exploring} and its associated dissipative dynamics may be an interesting direction for future studies.
Alternative strategies for systematically exploring multiple solutions have been developed in related contexts, including, for instance, the generalized high-index saddle dynamics framework~\cite{yin2021searching}, which addresses the problem of constructing solution landscapes and connecting stationary points.

Despite this limitation, the proposed framework significantly extends the practical range of applicability of continuous Newton methods, enabling them to cross singularities that would otherwise cause numerical breakdown and to converge to solutions that are unreachable by standard integration.

\subsection*{Acknowledgments}
G.F. and I.G.K acknowledge the Department of Energy (DOE) support under grant no. DE-SC0024162. I.G.K. also acknowledges partial support by the National Science Foundation under grant nos. CPS2223987 and FDT2436738.
J.D.H. acknowledge the National Science Foundation under grant no. CCF2331400.

\subsection*{Code and Data Availability}
The MATLAB code implementing the proposed algorithms, along with the scripts used to generate the numerical results presented in this work, will be made available upon publication.

%\bibliographystyle{plainnat}
%\bibliography{infinity}	

%\section*{Acknowledgments}
%We would like to acknowledge ....

\end{document}